\documentclass[a4paper]{amsart}
\usepackage{bbm} %,ulem}
\usepackage{graphicx}
\usepackage{subcaption}
\usepackage{amssymb}
\usepackage[utf8]{inputenc}
\usepackage{amsfonts}
\usepackage{amsmath}
\usepackage{amssymb}
\usepackage{amsthm}
\usepackage{url}
\usepackage[colorlinks,citecolor=blue]{hyperref}
\usepackage{mathabx}
\usepackage{wrapfig}
\usepackage{dsfont}
\usepackage{resizegather}
\usepackage{yfonts}
\usepackage{geometry}
\usepackage[dvipsnames]{xcolor}
\usepackage{algpseudocode}
\usepackage{algorithmicx}
\usepackage{algorithm}
\usepackage[pagewise]{lineno}\nolinenumbers

\DeclareMathOperator*{\argmin}{arg\,min}

\providecommand{\U}[1]{\protect\rule{.1in}{.1in}}

\newtheorem{prop}{Proposition}[section]
\newtheorem{cor}[prop]{Corollary}

\newtheorem{rmk}[prop]{Remark}
\newtheorem{lem}[prop]{Lemma}

\newtheorem{theo}[prop]{Theorem}

\newcommand{\tr}{\mbox{\rm Tr}}

\newcommand{\EE}{\mathbb{E}}

\newcommand{\LL}{\mathbb{L}}

\newcommand{\PP}{\mathbb{P}}

\newcommand{\RR}{\mathbb{R}}

\newcommand{\UU}{\mathbb{U}}
\newcommand{\VV}{\mathbb{V}}
\newcommand{\XX}{\mathbb{X}}
\newcommand{\WW}{\mathbb{W}}
\newcommand{\YY}{\mathbb{Y}}

\newcommand{\Ba}{ {\mathcal B }}

\newcommand{\Da}{ {\mathcal D }}
\newcommand{\La}{ {\mathcal L }}

\newcommand{\Ka}{ {\mathcal K }}

\newcommand{\Sa}{ {\mathcal S }}
\newcommand{\Ra}{ {\mathcal R }}
\newcommand{\Va}{ {\mathcal V }}
\newcommand{\Ua}{ {\mathcal U }}

\newcommand{\Ga}{ {\mathcal G }}
\newcommand{\Qa}{ {\mathcal Q }}

\newcommand{\Xa}{ {\mathcal X }}
\newcommand{\Ma}{ {\mathcal M }}

\newcommand{\Ha}{ {\mathcal H }}

\newcommand{\Pa}{ {\mathcal P }}
\newcommand{\Za}{ {\mathcal Z }}
\newcommand{\Ya}{ {\mathcal Y }}
\newcommand{\Wa}{ {\mathcal W }}

\newcommand{\point}{\mbox{\LARGE .}}
\newcommand{\cqfd}{\hfill\blbx \\}
\def\blbx{\hbox{\vrule height 5pt width 5pt depth 0pt}\medskip}

\def \PP{\mathbb{P}}
\def \RR{\mathbb{R}}

\def \EE{\mathbb{E}}

\def \LL{\mathbb{L}}

\def \WW{\mathbb{W}}

\newcommand{\cchi}{\protect\raisebox{2pt}{$\chi$}}

\usepackage[textwidth=2.5cm, textsize=scriptsize]{todonotes}

\newcommand{\vertiii}[1]{{\left\vert\kern-0.25ex\left\vert\kern-0.25ex\left\vert #1
    \right\vert\kern-0.25ex\right\vert\kern-0.25ex\right\vert}}

\usepackage{amsopn}

\numberwithin{equation}{section}

\title[On the contraction properties of Sinkhorn semigroups]{On the contraction properties of Sinkhorn semigroups}

\author{O. Deniz Akyildiz}
\address{Department of Mathematics, Imperial College London, UK}
\email{\textcolor{blue}{\footnotesize \texttt{deniz.akyildiz@imperial.ac.uk}}}

\author{Pierre Del Moral}
\address{Centre de Recherche Inria Bordeaux Sud-Ouest, Talence, France}
\email{\textcolor{blue}{\footnotesize \texttt{pierre.del-moral@inria.fr}}}

\author{Joaqu\'in Miguez}
\address{Department of Signal Theory \& Communications, Universidad Carlos III de Madrid, Spain}
\email{\textcolor{blue}{\footnotesize \texttt{joaquin.miguez@uc3m.es}}}


\begin{document}

\maketitle
\begin{abstract}

We develop a novel stability theory for Sinkhorn semigroups based on Lyapunov techniques and quantitative contraction coefficients, and establish exponential convergence of Sinkhorn iterations on weighted Banach spaces. This operator-theoretic framework yields explicit exponential decay rates of Sinkhorn iterates toward Schrödinger bridges with respect to a broad class of  $\phi$-divergences and Kantorovich-type distances, including relative entropy, squared Hellinger integrals, $\alpha$-divergences, weighted total variation norms, and Wasserstein distances. To the best of our knowledge, these results provide the first systematic contraction inequalities of this kind for entropic transport and the Sinkhorn algorithm.

We further introduce Lyapunov contraction principles under minimal regularity assumptions, leading to quantitative exponential stability estimates for a large family of Sinkhorn semigroups. The framework applies to models with polynomially growing potentials and heavy-tailed marginals on general normed spaces, as well as to more structured boundary state-space models, including semicircle transitions and Beta, Weibull, and exponential marginals, together with semi-compact settings. Finally, our approach extends naturally to statistical finite mixtures of such models, including kernel-based density estimators arising in modern generative modeling.

\end{abstract}

% \\
% \\ 
  % \\
% \textbf{Keywords:}  Entropic optimal transport, Sinkhorn's algorithm, iterative proportional fitting procedure, Gaussian processes, Schr\"odinger bridges, Monge maps, Riccati matrix difference equations.\\
% \\
% \noindent\textbf{Mathematics Subject Classification:} Primary 49N05, 49Q22, 94A17, 62J99, 60J20; secondary 62C10, 35Q49.

% \newpage
% \tableofcontents

%%%%%%%%%%%%%
%			%
%%%%%%%%%%%%%
\section{Introduction}

\subsection{Entropic optimal transport}

Let $(\XX,\lambda)$ and $(\YY,\nu)$ be locally compact Polish spaces equipped with locally bounded positive measures $\lambda$ and $\nu$. We consider probability measures of the form
\begin{equation}\label{ref-intro-UV}
\lambda_U(dx):=e^{-U(x)}~\lambda(dx)\quad \mbox{\rm and}\quad
\nu_V(dy):=e^{-V(y)}~\nu(dy)
\end{equation}
where $U:\XX\to\RR$ and $V:\YY\to\RR$ are measurable potential functions chosen so that $\lambda_U$ and $\nu_V$ are normalized.
Let $\Pi_{U,V}$ be the set of probability measures $\Pa$ on  $(\XX\times\YY)$, with 
prescribed marginals $(\lambda_U,\nu_V)$. 
We introduce a reference Markov kernel $\Qa$ from $\XX$ to $\YY$ of the form
 \begin{equation}\label{def-Qa}
\Qa(x,dy)=e^{-W(x,y)}~\nu(dy)\quad\mbox{\rm and set}\quad
\Ra(y,dx):=e^{-W(x,y)}~\lambda(dx).
\end{equation}
where $W:\XX\times\YY\to\RR$ is a measurable transition potential. The associated reference probability measure on $\XX\times\YY$ is defined by
\begin{equation}\label{ref-ref}
\Pa_0(d(x,y)):=
(\lambda_U\times \Qa)(d(x,y)):=\lambda_U(dx)~\Qa(x,dy).
\end{equation}
The entropic optimal transport problem, also known as the Schrödinger bridge problem, with respect to the reference measure $\Pa_0$ consists in finding
\begin{equation}\label{def-entropy-pb-v2}
\PP=\argmin_{\Pa\,\in\, \Pi_{U,V}}\mbox{\rm Ent}(\Pa~|~\Pa_0),\quad \mbox{\rm with the relative entropy}\quad
\text{Ent}(\Pa~|~\Pa_0):=\int \log\frac{d\Pa}{d\Pa_0}~d\Pa,
\end{equation}
where $d\Pa/d\Pa_0$ denotes the Radon-Nikodym derivative of $\Pa$ with respect to $\Pa_0$.
We implicitly assume that there exists at least one $\Pa\in\Pi_{U,V}$ such that
$\mbox{\rm Ent}(\Pa|\Pa_0)<\infty$. This condition ensures  the existence of a solution  $\PP$ to (\ref{def-entropy-pb-v2})  (cf. the seminal article by Csisz\'ar~\cite{csiszar-2}, as well as Section 6 in the Lecture Notes by Nutz~\cite{nutz}, and also the survey article by L\'eonard~\cite{leonard} and references therein). 
Note that the product distribution
 $$
( \lambda_U\otimes \nu_V)(d(x,y)):= \lambda_U(dx)~\nu_V(dy)
$$
trivially belongs to $\Pi_{U,V}$. Moreover, if $\Pa=(\lambda_U\times\Ka)\,\in\, \Pi_{U,V}$ for some Markov kernel $\Ka(x,dy)$, then
 $$
 \displaystyle\frac{d\Pa}{d\Pa_0}(x,y)=\frac{d(\lambda_U\otimes \nu_V)}{d\Pa_0}(x,y)\times 
  \frac{d\Pa}{d(\lambda_U\otimes \nu_V)}(x,y)
 $$
 which yields the decomposition
 $$
 \mbox{\rm Ent}(\Pa~|~\Pa_0)=\int~W(x,y)~\Pa(d(x,y))+
\mbox{\rm Ent}(\Pa~|~\lambda_U\otimes \nu_V)-\int V(y)\nu_V(dy).
 $$
Choosing a transition potential  of the form $W(x,y)=\Wa(x,y)/t$,  for some parameter $t>0$, problem (\ref{def-entropy-pb-v2}) becomes equivalent to the regularized optimal transport problem
\begin{equation}\label{reg-version}
\PP=\argmin_{\Pa\,\in\, \Pi_{U,V}}\left(\int~\Wa(x,y)~\Pa(d(x,y))+t~
\mbox{\rm Ent}(\Pa~|~\lambda_U\otimes \nu_V)\right).
\end{equation}
In the optimal transport literature, $\Wa$ is referred to as the cost function and
$t>0$ as the entropic regularization parameter. The limiting case $t=0$ formally
corresponds to the classical Kantorovich optimal transport problem with cost $\Wa$.
For further background and developments, we refer to
\cite{dobrushin2006perturbation,rachev1985monge,villani2008optimal}
and the references therein.

%%%%%%%%%%
%
%%%%%%%%%%%

\subsection{Background} % and contributions}

The Schr\"odinger bridge problem \eqref{def-entropy-pb-v2} consists in constructing a random evolution that minimizes the relative entropy with respect to a given reference model while matching prescribed initial and final marginal distributions. This convex optimization problem can be interpreted as an entropic regularization of the Monge--Kantorovich optimal transport problem; we refer to \cite{leonard,nutz,peyre} for recent surveys and comprehensive bibliographies.

Entropic optimal transport and Schr\"odinger bridges have become central tools in a wide range of applications, including generative modeling and machine learning \cite{arjovsky,doucet-bortoli,kolouri,peyre,rioux2023semi,rioux2024entropic}, statistical barycenters \cite{agueh,andersen,bigot,cramer-2,cuturi-w,loubes}, economics \cite{bojitov}, computer vision \cite{dominitz,solomon}, and control theory \cite{chen,chen-2}. 

Under mild integrability conditions, the Schr\"odinger bridge problem admits a unique solution that can be characterized as the fixed point of the Sinkhorn equations, also known as the iterative proportional fitting procedure. We refer to the pioneering works \cite{deming,fortet,kruithof,sinkhorn-2,sinkhorn-3} and to more recent developments in \cite{chen,cuturi,leonard-2} for historical background and modern perspectives.
Sinkhorn bridges involve nonlinear integral operators acting on infinite-dimensional function spaces, and explicit solutions are generally unavailable. In finite state spaces, the Sinkhorn algorithm reduces to alternating row and column normalizations of matrices. For Gaussian models, however, Sinkhorn equations can be solved analytically using sequential Bayesian conjugacy arguments \cite{adm-24,cramer-2,cramer-3,janadi}.

Sub-linear convergence rates of Sinkhorn equations have also been developed in the articles~\cite{alts-2017,chak-2018,dvu-2017}. Linear rates  on non-necessarily compact spaces were first obtained by L\'eger in~\cite{leger} using elegant gradient descent and Bregman divergence techniques (see also the recent articles~\cite{doucet-bortoli,karimi-2024}). Refined convergence rates at least one order faster for a general class of models including sub-Gaussian target marginals has also been developed in~\cite{promit-2022} under some exponential moments estimates on Schr\" odinger potentials. Most existing results on the exponential convergence of Sinkhorn iterations concern finite or compact state spaces, or bounded cost functions, typically relying on Hilbert’s projective metric and related techniques \cite{borwein,fienberg,sinkhorn-2,sinkhorn-3,soules,chen,deligiannidis,franklin,marino}. See also \cite{greco2023non} for refined analyses on the torus and \cite{chizat2025sharper} for recent advances in compact settings.

As emphasized in \cite{cohen}, contraction in the Hilbert metric is a very strong property, as this metric dominates all $\phi$-divergences. While Hilbert metric contraction implies exponential decay for a broad class of divergences (including total variation, relative entropy, and Wasserstein distances), it is generally difficult to obtain sharp contraction estimates for specific divergence criteria. Moreover, even on totally bounded spaces such as $\XX=]0,1[=\YY$, the cost function and Schr\"odinger potentials may be unbounded, and the associated Markov operators may have unit Birkhoff contraction coefficient; see Sections~\ref{birk-sec}, \ref{bst-mod-sec}, and \ref{beta-w-sec}.

Quadratic costs on Euclidean spaces provide another important class of unbounded models. In this setting, \cite{adm-24} develops a refined stability analysis for Gaussian Schr\"odinger bridges based on matrix Riccati equations. 
Beyond Gaussian models, only a few recent works address exponential convergence of Sinkhorn iterations for non-compact spaces and unbounded costs, notably \cite{durmus,chiarini,eckstein2025hilbert}.
The article \cite{durmus} establishes quantitative exponential convergence for symmetric quadratic costs and marginal distributions with asymptotically positive log-concavity, using diffusion reflection couplings and convexity propagation techniques \cite{conforti2024weak,eberle2016reflection}. These results were further refined in \cite{chiarini}, which applies to all regularization parameters and certain perturbations of log-concave models. A central tool in \cite{chiarini} is a stability theorem controlling the entropy of Schr\"odinger bridges by the entropy of their marginals, under uniform semi-concavity assumptions on unknown dual Schr\" odinger potentials.
This approach is also limited to symmetric quadratic-type costs and does not cover important classes of nonsymmetric reference transitions, such as those arising in linear Gaussian models for Ornstein--Uhlenbeck processes and denoising diffusion models. Recent work \cite{del2025stability} partially overcomes these limitations by establishing entropic continuity results for nonsymmetric transitions and strongly convex-at-infinity potentials.

An alternative approach based on generalized Hilbert metrics for weighted function spaces is developed in \cite{eckstein2025hilbert}.  Extending results in~\cite{chen,deligiannidis} to unbounded cost functions, the recent article ~\cite{eckstein2025hilbert} presents exponential decays of these new classes of
Hilbert’s projective metrics with respect to some cone of functions
 in settings where the marginal distributions have sufficiently light tails compared to the growth of the cost function.   These results also imply the exponential convergence of the total variation distance. On normed finite dimensional spaces, the main condition in ~\cite{eckstein2025hilbert} is that $W$ grows faster that a polynomial of order $p$ and the tails of both marginals decay faster that $e^{-r^{p+q}}$ for some $q>0$. This condition doesn't cover Gaussian models.

 Finally, Sinkhorn semigroups belong to the broader class of time-inhomogeneous Markov processes that share a common invariant measure at each step \cite{adm-24}. Stability analysis for such processes is notoriously delicate, as transition operators may vary substantially over time. 
 These classes of time varying Markov semigroups arise in a variety of areas in applied mathematics including on nonlinear filtering, physics and molecular chemistry, see for instance~\cite{dpa,dm-04,dg-01,horton-dp,saloff-zuniga-3}  and references therein. 
 While general tools exist—such as coupling methods, spectral techniques, and functional inequalities \cite{douc,saloff-zuniga,saloff-zuniga-2}-their direct application to Sinkhorn semigroups is challenging due to the nonlinear conjugate structure of the updates.
 This motivates the development of new analytical tools. In this work, we introduce a Lyapunov-based, operator-theoretic framework for the quantitative stability analysis of Sinkhorn semigroups on weighted Banach spaces.

%%%%%%%%%%
%
%%%%%%%%%%
\subsection{Contributions}

We investigate Lyapunov approaches and related weighted-norm contraction techniques to analyze the convergence of Sinkhorn semigroups. These semigroup techniques are basic quantitative methods for studying the stability of nonlinear and time-varying stochastic models, see for instance~\cite{horton-dp} as well as Section 2 in~\cite{dpa}. In the context of time-homogeneous models, we also refer to the pioneering articles by Hasminskii~\cite{hasminskii} and Meyn and Tweedie~\cite{meyn-1}, see also the book~\cite{meyn-2} and references therein. 

These powerful off-the-shelf semigroup techniques are mainly based on two conditions: 
\begin{enumerate}
 \item The existence of a Lyapunov function satisfying a drift-type condition.
 \item  A local minorization condition on the sub-level sets of the Lyapunov function.
 \end{enumerate}
Such drift-minorization approaches are well-established techniques to localize the minorization conditions on the compact sub-level sets of the Lyapunov function. In this context, Markov operators act on weighted Banach spaces and operator norms coincide with the  contraction coefficients of the Markov integral operator. For a more detailed discussion on this contraction coefficients we refer to~\cite{dpa,delgerber2025,horton-dp,penev}.

This Lyapunov approach and the operator-theoretic framework  presented in this article strongly differ from the approaches discussed earlier and taken in the series of recent articles~\cite{chiarini,durmus,del2025stability,eckstein2025hilbert}.
 Our main contributions can be summarized in the following way:
\begin{itemize}
\item {\em Strict contraction for bounded costs:} For bounded transition potentials, we establish strict contraction inequalities for a broad class of $\phi$-divergences, including $\LL_p$ norms, Havrda–Charvat and $\alpha$-divergences, Jensen–Shannon divergence, squared Hellinger and Kakutani–Hellinger distances, as well as Jeffreys and R\'enyi divergences (see Theorem~\ref{theo1-intt} and Section~\ref{strong-mixing-sec}).\\

\item  {\em Weighted-norm stability for unbounded models:} For locally bounded costs and marginal potentials with compact sub-level sets, we introduce weighted total variation and Kantorovich–Wasserstein contraction estimates, together with relative entropy bounds for Schrödinger bridges (Sections~\ref{sec-lyap-i}, \ref{sec-wass}, and \ref{sec-entrr}). The choice of Lyapunov functions and Wasserstein metrics is explicitly adapted to the marginal distributions.\\

\item  {\em General Lyapunov contraction principles:} We formulate practical drift–minorization conditions yielding quantitative exponential convergence rates for a large class of Sinkhorn semigroups. These results apply to models with polynomial growth, heavy-tailed marginals, boundary state spaces, semi-compact domains, and finite mixture models, including kernel-based density estimators used in generative modeling.\\

\end{itemize}

To the best of our knowledge, these results provide the first general exponential stability estimates for Sinkhorn semigroups with respect to weighted norms, $\phi$-divergences, and Wasserstein-type distances in non-necessarily compact and unbounded settings.

%%%%%%
%
%%%%%%
\subsection{Outline of the article}

After introducing notation in Section~\ref{ssBasic}, Section~\ref{sec-models} recalls the Schr\"odinger bridge problem and the associated Sinkhorn semigroups, emphasizing their interpretation as time-inhomogeneous nonlinear Markov evolutions (see also~\cite{adm-24}).

Section~\ref{stab-th-intro} presents a concise overview of the main stability and contraction results established in this work. Motivating applications are already discussed in Section~\ref{sec-app-com}, which highlights the relevance of our approach for non-compact state spaces, unbounded costs, and heavy-tailed marginals.

Sections~\ref{contract-Markov-sec}–\ref{ccoef-sec} develop the operator-theoretic foundations underlying our analysis. We review Lyapunov drift–minorization techniques and weighted Banach space contractions for Markov integral operators (Section~\ref{contract-Markov-sec}, with complements in Appendix~\ref{phi-psinorm-sec}), recall the Hilbert projective metric and $\phi$-divergences (Section~\ref{metrics-sec}), and introduce the associated contraction coefficients (Section~\ref{ccoef-sec}).

Building on this framework, Section~\ref{sec-sinkhorn-contract} adapts weighted-norm and Lyapunov contraction techniques to Sinkhorn semigroups. Strong mixing regimes are analyzed in Section~\ref{strong-mixing-sec}, where we establish strict contraction properties for a broad class of $\phi$-divergences, including relative entropy, total variation, $\LL_p$ norms, Hellinger-type divergences, Jensen–Shannon, Jeffreys, and R\'enyi divergences.
General drift–minorization conditions ensuring quantitative contraction bounds are developed in Section~\ref{sec-sinkhorn-contract-reg}. The main exponential stability theorems are presented in Section~\ref{main-theo-sec}. Section~\ref{sec-lyap-i} introduces a simple and practically verifiable Lyapunov condition yielding exponential convergence on weighted Banach spaces. Extensions to Wasserstein-type semi-metrics are derived in Section~\ref{sec-wass}, while Section~\ref{sec-entrr} is devoted to relative-entropy decay and the convergence of Sinkhorn bridges.

Section~\ref{sec-illustrations} illustrates the scope of the theory through a wide range of examples, including polynomial growth models (Section~\ref{sec-poly}), heavy-tailed marginals (Section~\ref{heavy-tail-sec}), bounded Lipschitz domains (Section~\ref{bst-mod-sec}), Beta and Weibull models (Sections~\ref{beta-w-sec} and~\ref{w-sec}), and multivariate linear Gaussian models (Section~\ref{sec-lin-gauss}). Section~\ref{finite-mix-sec} extends the analysis to finite mixture models, including non-Gaussian and kernel-based constructions relevant to generative modeling.

Detailed proofs and complementary material are collected in the appendices. Appendix~\ref{GS-sec} focuses on Gaussian Sinkhorn algorithms, Appendices~\ref{apHilbert} and~\ref{phi-psinorm-sec} provide additional results on Hilbert metrics and weighted-norm contractions, Appendix~\ref{appSchrodinger} reviews Schr\"odinger potentials, and Appendix~\ref{tech-proof-ap} contains further technical arguments.

\subsection{Some basic notation} \label{ssBasic}

We let  $\Ba(\XX)$ be the algebra  of locally bounded  measurable functions on 
a locally compact Polish space $\XX$. We denote by $\Ba_b(\XX)\subset \Ba(\XX)$ the sub-algebra  of bounded measurable functions $f$ on $\XX$ endowed with the supremum  norm $\Vert f\Vert:=\sup_{x\in \XX}|f(x)|$. For a given uniformly positive function $\varphi\in \Ba(\XX)$, we let $\Ba_{\varphi}(\XX)\subset \Ba(\XX)$ be the Banach space of functions $f\in \Ba(\XX)$ equipped with the norm $\Vert f\Vert_{\varphi}:=\Vert f/\varphi\Vert$.  
 
Let $\Ma(\XX)$ be the set of locally bounded signed measures on $\XX$, and let $\Ma_b(\XX)\subset \Ma(\XX)$ the subset of bounded measures. We denote the convex subset of probability measures on $\XX$ as $\Ma_1(\XX)\subset \Ma_b(\XX)$.

The Lebesgue integral of an integrable function $f\in \Ba(\XX)$ w.r.t. $\nu\in \Ma(\XX)$ is denoted by
$$
\nu(f):=\int_{\XX} f(x)~\nu(dx).
$$ 
For indicator functions $f=1_{A}$ of a measurable subset $A\subset \XX$ sometimes we write $\nu(A)$ instead of $\nu(1_A)$.

Consider  a bounded integral operator  $\Ka(x,dy)$ from $\XX$ into $\YY$, 
a measure $\mu \in\Ma(\XX)$  and a function $f \in\Ba(\YY)$. We denote by $\mu \Ka\in \Ma(\YY)$ and $\Ka(f)\in \Ba(\XX)$ the measure and the function defined by  
 \begin{equation}\label{basic-int-op}
(\mu \Ka)(dy):=\int~\mu(dx) \Ka(x,dy)\quad \mbox{\rm and}\quad \Ka(f)(x):=\int \Ka(x,dy) f(y).
 \end{equation}
For indicator functions $f=1_{A}$ sometimes we write $\Ka(x,A)$ instead of $\Ka(1_A)(x)$.
For bounded integral operator  $\Ka_1(x,dy)$ from $\XX$ into $\YY$ and  $\Ka_2(y,dx)$ from $\YY$ into $\XX$,  we also denote by $(\Ka_1\Ka_2)$ the bounded integral operator from $\XX$ into $\XX$ defined by the integral composition
$$
(\Ka_1\Ka_2)(x_1,dx_2):=\int\Ka_1(x_1,dy)
\Ka_2(y,dx_2).
$$
Given a  Markov transition $\Ka(x,dy)$ from $\XX$ into $\YY$ and
a measure $\mu \in\Ma(\XX)$ we set
$$
(\mu\times \Ka)(d(x,y)):=\mu(dx)\, \Ka(x,dy)
\quad\mbox{\rm and}\quad
(\mu\times \Ka)^{\flat}(d(y,x)):=\mu(dx)\, \Ka(x,dy). 
$$

A  measure $\nu_1\in \Ma(\XX)$ is said to be absolutely continuous with respect to another measure $\nu_2\in \Ma(\XX)$ and we write $\nu_1\ll \nu_2$ if $\nu_2(A)=0$ implies $\nu_1(A)=0$ for any $A\in \XX$. Whenever $\nu_1\ll \nu_2$, we denote by ${d\nu_1}/{d\nu_2}$ the Radon-Nikodym derivative of $\nu_1$ w.r.t. $\nu_2$. We also consider the partial order, denoted $\nu_1\leq  \nu_2$, when $\nu_1(A)\leq \nu_2(A)$ for any measurable subset $A\subset \XX$. The measures  are said comparable when we have $c~\nu_1\leq  \nu_2\leq  c^{-1}~\nu_1$ for some scalar constant $c>0$.

The positive and negative part of $a\in\RR$ are defined respectively by $a_+=\max(a,0)$
and  $a_-=\max(-a,0)$. Given $a,b\in\RR$ we also set $a\vee b=\max(a,b)$ and
$a\wedge b=\min(a,b)$.

\subsubsection*{Kantorovich and Wasserstein criteria}

We consider the weighted total variation norm associated with
some positive function  $\varphi\in \Ba(\XX)$ defined as
\begin{equation}\label{eq-kanto}
 \vertiii{\mu_1-\mu_2}_{\varphi}:=\inf_{\Pa\in \Pi(\mu_1,\mu_2)}\int~(\varphi(x_1)+
\varphi(x_2))~1_{x_1\not=x_2}~\Pa(d(x_1,x_2)),
\end{equation}
where $\Pi(\mu_1,\mu_2)$ stands for the set of coupling probability measures $\Pa$ on the product space $(\XX\times\XX)$, with 
prescribed first and second coordinate marginals $(\mu_1,\mu_2)$ such that $\mu_i(\varphi)<\infty$.  For a given semi-metric $\varrho(x_1,x_2)$ on $\XX$ and $p\geq 1$ we also consider the Wasserstein semi-metric $\WW_{p,\varrho}$ defined by
\begin{equation}\label{def-WW-varrho}
\WW_{p,\varrho}(\mu_1,\mu_2):= \inf_{\Pa\in \Pi(\mu_1,\mu_2)}\left(\int~\varrho(x_1,x_2)^p~\Pa(d(x_1,x_2))\right)^{1/p}.
\end{equation}
When $p=1$ we simplify notation and we write $\WW_{\varrho}$ instead of 
$\WW_{1,\varrho}$. When $\varrho_{\Vert \point\Vert}(x_1,x_2):=\Vert x_1-x_2\Vert$ is associated with some norm $\Vert \point\Vert$ on some vector spaces $\XX$, when there are no possible  confusions we write $\WW_p$ instead of $\WW_{p,\varrho_{\Vert \point\Vert}}$.

When $\varrho(x_1,x_2)$ is a metric on $\XX$ (possibly different from the metric of the Polish space) we recall (see for instance~\cite{edwards2011kantorovich,fernique-81,kellerer1982duality,Szulga-83,villani2008optimal}) the Kantorovich-Rubinstein duality
$$
\WW_{\varrho}(\mu_1,\mu_2)=\sup{\left\{(\mu_1-\mu_2)(f)~;~f~:~\XX\mapsto \RR\quad s.t.\quad \vert f(x_1)-f(x_2)\vert\leq \varrho(x_1,x_2) \right\}}.
$$

Note that
$$
\varrho_{\varphi}
(x_1,x_2)=(\varphi(x_1)+
\varphi(x_2))~1_{x_1\not=x_2}\quad
\Longrightarrow
\quad\WW_{\varrho_{\varphi}}(\mu_1,\mu_2)= \vertiii{\mu_1-\mu_2}_{\varphi}.
$$
More generally, for any given semi-metric $\varrho(x_1,x_2)$ on $\XX$ such that
\begin{equation}\label{gtow-0}
\varphi(x_1)+\varphi(x_2)\geq \varrho(x_1,x_2)\quad \mbox{\rm and}\quad \varphi(x)\geq c>0
\end{equation}
we have the estimates
\begin{equation}\label{gtow}
\WW_{\varrho}(\mu_1,\mu_2)\vee \left(c~\Vert \mu_1-\mu_2\Vert_{\tiny tv}\right)\leq \vertiii{\mu_1-\mu_2}_{\varphi}\,,\quad\mbox{\rm with the total variation norm $\Vert \point\Vert_{\tiny tv}$.}
\end{equation}
For instance on a given normed vector space $(\XX,\Vert\point\Vert)$, for any $b,c>0$
and $p,q\geq 1$, if
$$
%\begin{array}{l}
\displaystyle \varphi(x):=c~\exp{\left(b\Vert x\Vert^p\right)} %\\
$$
then
%\\
%\displaystyle\Longrightarrow
$$
\varphi(x_1)+\varphi(x_2)\geq \varrho_{\tiny exp}(x_1,x_2):=2c~\left(
\exp{\left(\frac{b}{2^{1+(p-1)_+}}~\Vert x_1-x_2\Vert^p\right)}-1
\right).
%\end{array}
$$
We check this claim using the fact that
$u+v\geq 2\sqrt{u}\sqrt{v}$ and $(u^p+v^p)\geq 2^{-(p-1)_+}(u+v)^p$ for any $p>0$ and $u,v\geq 0$. In this context, we have
$$
\WW_{\varrho_{\tiny exp}}(\mu_1,\mu_2)\leq\vertiii{\mu_1-\mu_2}_{\varphi}.
$$

The choice of the semi-metric $\varrho$ in (\ref{gtow}) is not unique.  Indeed, developing  the exponential function for any $q\geq 1$ we also have the estimate
\begin{equation}\label{intro-w-p-q}
\varphi(x_1)+\varphi(x_2)\geq \varrho_{p,q}(x_1,x_2):=c_{p,q}~\Vert x_1-x_2\Vert^{pq}\quad\mbox{\rm with}\quad
c_{p,q}:=2c~\frac{1}{q!}\left(
\frac{b}{2^{1+(p-1)_+}}\right)^q,
\end{equation}
and
$$
c_{p,q}~\WW_{pq}(\mu_1,\mu_2)^{pq}\leq \WW_{\varrho_{p,q}}(\mu_1,\mu_2)\leq\vertiii{\mu_1-\mu_2}_{\varphi}.
$$

\subsubsection*{Entropy and divergences}

Consider a convex function 
\begin{equation}\label{Phi-ref}
\begin{array}{cccl}
    \Phi:&\RR_+^2 &\rightarrow &\RR\cup\{+\infty\}\\
         &(u,v)   &\mapsto&\Phi(u,v)\quad\mbox{\rm s.t.  $\forall a\in\RR_+$~~~ $\Phi( a u,a v)=a~ \Phi(u,v)$ and  $\Phi(1,1)=0$.}
\end{array}
\end{equation}
The $\Phi$-entropy between two probability measures $\nu_1,\nu_2\in \Ma_1(S)$ is defined for any dominating measure $\gamma$ (such that $\gamma\gg \nu_1$ and $\gamma\gg\nu_2$) by the formula
\begin{equation}\label{def-DPhi}
D_{\Phi}(\nu_1,\nu_2):=\int~\Phi\left(\frac{d\nu_1}{d\gamma}(x),\frac{d\nu_2}{d\gamma}(x)\right)~\gamma(dx).
\end{equation}
By homogeneity arguments, the definition (\ref{def-DPhi}) does not depend on the choice of the dominating measure $\gamma$. Note that
 \begin{equation}\label{tv-osc}
 \text{if} \quad \Phi(u,v)=\Phi_0(u,v):=\frac{1}{2}~\vert u-v\vert \quad\text{then}\quad D_{\Phi_0}(\nu_1,\nu_2)=\Vert \nu_1-\nu_2\Vert_{\sf tv}.
 \end{equation}

Whenever $\Phi(1,0)=\infty$, the $\Phi$-entropy coincides with the $\phi$-divergence in the sense of Csisz\'ar~\cite{csiszar}, i.e., we also have
 $$
 D_{\Phi}(\nu_1,\nu_2)=\int~\phi\left(\frac{d\nu_1}{d\nu_2}(x)\right)~\nu_2(dx)\quad \text{where} \quad  \phi(u):=\Phi(u,1)
 $$
 and we have used the convention $D_{\Phi}(\nu_1,\nu_2)=\infty$ when $\nu_1\not\ll \nu_2$. Therefore, the $\Phi$-entropy generalizes the classical relative entropy (a.k.a. Kullback-Leibler divergence or $I$-divergence). More precisely, choosing the functions
    $$
     \Phi(u,v)= \Phi_1(u,v):=u\log{(u/v)}\quad \mbox{\rm and}\quad \phi(u)= u\log{u}\
    $$
we recover the classical formulae
 $$
 D_{\Phi_1}(\nu_1,\nu_2)=
 \mbox{\rm Ent}(\nu_1~|~\nu_2):=\int~\log{\left(\frac{d\nu_1}{d\nu_2}(x)\right)}~\nu_1(dx)=\int~\phi_1\left(\frac{d\nu_1}{d\nu_2}(x)\right)~\nu_2(dx),
$$
with the convex function with unique minimum $\phi_1(1)=0= \partial\phi_1(1)$ defined by
  $$
  \phi_1(u)=u\log{u}-(u-1).
  $$

\section{Statement of main results}\label{sec-models}
 This section introduces the Sinkhorn bridge construction, the associated time-inhomogeneous Sinkhorn semigroups, and states our main exponential stability results. Detailed proofs and refinements are deferred to later sections.
\subsection{Sinkhorn bridges}
The Schr\"odinger bridge distribution $\PP$ defined  in (\ref{def-entropy-pb-v2}) can rarely be computed analytically. {However, solutions can be approximated} efficiently using the Sinkhorn algorithm, also referred to as the iterative proportional fitting procedure~\cite{cuturi,sinkhorn-2,sinkhorn-3}. This procedure is defined sequentially  for any $n\geq 0$ by a collection of probability distributions {($\Pa_n$) and transition operators ($\Sa_n$) given by}
\begin{equation}\label{def-Pa-n}
 \Pa_{2n}=\lambda_U\times\Sa_{2n}
\quad \mbox{\rm and}\quad
 \Pa_{2n+1}=(\nu_V\times\Sa_{2n+1})^{\flat}
\end{equation}
starting from $\Pa_0$ with $\Sa_0=\Qa$, where, for a measure $\mu(dy)$ and a transition operator $\Sa(y,dx)$, the notation $(\mu\times\Sa)^\flat$ denotes the `reverse' measure
$$
(\mu\times\Sa)^\flat(d(x,y)):=\mu(dy) \Sa(y,dx).
$$
The  Markov transitions $\Sa_n$ in (\ref{def-Pa-n}) are defined  by the Sinkhorn iterations
\begin{equation}
\left\{\begin{array}{l}
(\pi_{2n}\times\Sa_{2n+1})^{\flat}=\lambda_U\times\Sa_{2n}\quad \mbox{\rm and}\quad
\pi_{2n+1}\times\Sa_{2(n+1)}=(\nu_V\times\Sa_{2n+1})^{\flat}\\
\\
\mbox{\rm with the distributions}\quad \pi_{2n}:=\lambda_U \Sa_{2n}\quad\mbox{\rm and}\quad \pi_{2n+1}:=\nu_V \Sa_{2n+1}.
\end{array}\right.
\label{s-2}
\end{equation}

In the dual formulation, the bridge distributions $\Pa_n$ are often written as
\begin{equation}\label{sinhorn-entropy-form-Sch}
\Pa_{n}(d(x,y))=e^{-W(x,y)}~\lambda_{U_n}(dx)~\nu_{V_n}(dy),
\end{equation}
with the initial potential functions $(U_0,V_0)=(U,0)$ and potential functions $(U_n,V_n)$ defined sequentially by the integral recursions
\begin{eqnarray}
U_{2n+1}&=&~U_{2n}~:=U+\log{\Qa(e^{-V_{2n}})}\nonumber\\
V_{2(n+1)}&=&V_{2n+1}:=V+\log{\Ra(e^{-U_{2n}})}.\label{prop-schp}
\end{eqnarray}
When $n\to\infty$, the sequence of distributions $\Pa_{2n}$ converges towards the Schr\"odinger bridge (\ref{def-entropy-pb-v2}) from $\lambda_U$ to $\nu_V$ (w.r.t. the reference measure $\Pa_0$) while the sequence $\Pa_{2n+1}$ converges towards the Schr\"odinger bridge from $\nu_V$ to $\lambda_U$ (w.r.t. the reference measure $\Pa_1$). 

Under mild integrability assumptions (see for instance Section 6 in~\cite{adm-24}) the Schr\"odinger bridge distribution $\PP$ solving (\ref{def-entropy-pb-v2}) can also be computed from any Sinkhorn bridge reference distribution $\Pa_{2n}$, in the sense that  for any $n\geq 0$ we have
\begin{equation}\label{def-entropy-pb-v2-2}
\PP=\argmin_{\Pa\,\in\, \Pi_{U,V}}\mbox{\rm Ent}(\Pa~|~\Pa_{2n}).
\end{equation}

The Sinkhorn transition operators can also be written in terms of the potential functions $(U_n,V_n)$, namely,
\begin{equation}\label{sinhorn-transitions-form-Sch}
\Sa_{2n}(x,dy)=\frac{\Qa(x,dy)e^{-V_{2n}(y)}}{\Qa(e^{-V_{2n}})(x)}\quad \mbox{and}\quad
\Sa_{2n+1}(y,dx)=\frac{\Ra(y,dx)~e^{-U_{2n+1}(x)}}{\Ra(e^{-U_{2n+1}})(y)},
\end{equation}
Moreover, Proposition 6.3 in~\cite{adm-24} shows that the potentials can be expressed through the marginal flows as
\begin{equation}\label{form-series} 
V_{2n}=V_0+\sum_{0\leq p< n}
\log\frac{ d\pi_{2p}}{d\nu_V}
\quad \mbox{and}\quad
U_{2n}=U_0+\sum_{0\leq p< n}\log{\frac{d \pi_{2p+1}}{d\lambda_U}}
\end{equation}
%for the potentials $(U_n,V_n)$, $n \ge 0$, expressed in terms of the distribution flow of the internal states of the Gibbs-loop process.

\subsection{Sinkhorn semigroups}\label{GL-sec}

The Sinkhorn iterations naturally induce a family of time-inhomogeneous Markov semigroups. For each $n\geq 0$ define the forward-backward Gibbs-type transitions
$$
 \Sa^{\circ}_{2n+1}:=\Sa_{2n}\Sa_{2n+1}\quad\mbox{\rm and}\quad
\Sa^{\circ}_{2(n+1)}:=\Sa_{2n+1}\Sa_{2(n+1)}.
$$
These kernels satisfy the fixed-point relations
\begin{equation}\label{fixed-points-gibbs}
\lambda_U \Sa^{\circ}_{2n+1}=\lambda_U \quad\mbox{\rm and}\quad
\nu_V \Sa^{\circ}_{2(n+1)}=\nu_V
\end{equation}
and yield the marginal evolution equations
\begin{equation}\label{gibbs-tv}
\pi_{2n}=\pi_{2(n-1)}\Sa^{\circ}_{2n}\quad\mbox{\rm and}\quad
\pi_{2n+1}=\pi_{2n-1}\Sa^{\circ}_{2n+1}.
\end{equation}

For $0\leq p\leq n$, we define the associated inhomogeneous semigroups by composition,
$$
\Sa^{\circ}_{2p,2n}=\Sa^{\circ}_{2p,2(n-1)} \Sa^{\circ}_{2n}
\quad\mbox{\rm and}\quad
\Sa^{\circ}_{2p+1,2n+1}=\Sa^{\circ}_{2p+1,2n-1} \Sa^{\circ}_{2n+1},
$$
with identity operators on the diagonal.
This representation highlights that Sinkhorn dynamics belong to the class of time-varying Markov chains sharing common invariant measures at each step.

\begin{rmk}\label{rmk-proximal}
Each kernel $ \Sa^{\circ}_{2n+1}:=\Sa_{2n}\Sa_{2n+1}$ corresponds to a two-block Gibbs sampler with target distribution
$$
\Pa_{2n}(d(x,y))=\lambda_U(dx)~\Sa_{2n}(x,dy)=\pi_{2n}(dy)~\Sa_{2n+1}(y,dx)
$$
For $n=0$, this recovers classical proximal and data-augmentation samplers; see, e.g.,
\cite{damlen1999gibbs,rendell2020global,vono2019split,vono2020asymptotically} and related works~\cite{chewi-chen,chewi-phd,jiaojiao,guan,lee,jiaming,mou}. 
\end{rmk}

\subsection{Some stability theorems}\label{stab-th-intro}

In contrast with Kantorovich criteria or related Wasserstein-type semi-distances, even when $\XX=\YY$ and the marginals coincide, $\lambda_U=\nu_V$, the identity-map deterministic coupling has infinite entropy. That is, for any reference measure $\Pa_0=(\lambda_U\times\Qa)$ of the form (\ref{ref-ref}) we have
$$
\Pa(d(x,y)):=\lambda_U(dx)~\delta_x(dy)\in \Pi_{U,V}\Longrightarrow
\mbox{\rm Ent}(\Pa~|~\Pa_0)=+\infty.
$$
This fact, together with (\ref{reg-version}), shows that the Schr\"odinger bridge  $\PP$ and the convergence analysis of Sinkhorn bridges $\Pa_{2n}$ crucially depend on the form of the reference transition $\Qa$ and the transition potential function $W$ discussed in (\ref{def-Qa}).

Markov transitions of the form (\ref{def-Qa}) with a bounded  potential $W$ typically arise on finite or compact spaces. For instance, Markov processes on compact manifolds, such as diffusions on bounded manifolds reflected at the boundary,  have Markov transitions of the form (\ref{def-Qa}) with
a uniformly bounded transition potential $W$ that depends on the time parameter~\cite{dpa,aronson,nash,varopoulos}.  The convergence analysis of Sinkhorn bridges in this context is often developed using Hilbert Hilbert projective metrics~\cite{chen,deligiannidis,franklin,marino}. 

In the context of bounded transition potentials,  the contraction principles of Markov operators presented in the article~\cite{dm-03} yield an alternative stability analysis that also
applies to any generalized $\Phi$-entropies introduced in (\ref{def-DPhi}).

\begin{theo}\label{theo1-intt}
 Assume that the transition potential $W$ is uniformly bounded. In this case, for any 
$\Phi$-entropy criteria $D_{\Phi}$ there is some parameters  $\rho\in ]0,1[$ and $c>0$ such that for any $n\geq 0$ we have the contraction inequalities
\begin{equation}\label{theo-intro-e1}
 D_{\Phi}\left(\pi_{2n},\nu_V\right)\leq \rho^n~
  D_{\Phi}\left(\pi_{0},\nu_V\right)\quad \mbox{and}\quad
  \left\Vert\log{({d \pi_{2n}}/{d\nu_V})}\right\Vert\leq c~ \rho^n.
 \end{equation}
The above  inequalities remain valid when we replace $(\pi_{2n},\nu_V)$ by 
$(\nu_V,\pi_{2n})$. Analogous bounds hold for
$(\lambda_U,\pi_{2n+1})$. Moreover,
\begin{equation}\label{theo-intro-ent1}
 \mbox{\rm Ent}(\PP~|~\Pa_{2n})  \leq   c~ \rho^{n}.
 \end{equation}
\end{theo}

The estimates (\ref{theo-intro-e1}) are a consequence of the sandwich-type contraction inequalities stated in Theorem~\ref{tx-st}. A more refined version of the entropy estimate (\ref{theo-intro-ent1}) in terms of Dobrushin contraction coefficients is provided in
Corollary~\ref{cor-entropy-wb}

We underline that even on totally bounded state spaces such as $\XX=\YY=]0,1[$ the transition potential function $W$ may be unbounded. Indeed, consider the reference Markov transition defined for any $x\in ]0,1[$ by
\begin{equation}\label{exintro}
\Qa(x,dy)=\frac{1}{q(x)}~\left(1-\vert x-y\vert\right)~1_{]0,1[}(y)~ dy\quad \mbox{\rm with}\quad \frac{1}{4}\leq q(x)=\frac{1}{2}-x(1-x)\leq\frac{1}{2}.
\end{equation}
In this context, we have
$$
\displaystyle W(x,y)=-\log{\left((x\wedge y)+(1-(x\vee y))\right)}+\log{q(x)}.
$$
The function $W(x,y)$ is locally bounded on $]0,1[^2$ but tends to $\infty$ on boundaries. For instance, we obtain
$$
W\left(\frac{1}{n},1-\frac{1}{n}\right)\geq
\log{\frac{n}{2}}-2\log{2}\stackrel{n\to\infty}{\longrightarrow}+\infty.
$$
This rather simple model  does not satisfy the boundedness condition of Theorem~\ref{theo1-intt}. 
To handle locally bounded transition potentials, consider the functions
\begin{equation}\label{UV-delta}
\Ua_{\delta}:=\delta U-W^{V}\quad\mbox{\rm and}\quad
\Va_{\delta}:=\delta V-W_{U}
\end{equation}
with the integrated costs
\begin{equation}\label{int-costs}
W_U(y):=\int~\lambda_U(dx)~W(x,y)\quad\text{and}\quad
W^V(x):=\int~\nu_V(dy)~W(x,y).
\end{equation}
Intuitively, the next theorem shows that Sinkhorn semigroups are stable as soon as the 
marginal-potentials  $(U,V)$ dominate at infinity the (integrated) transition-potentials $(W_U,W^V)$. In terms of potential energy conversion, the marginal potentials should be large enough to compensate high energy transition losses. A more thorough discussion on these regularity conditions is provided in Section~\ref{sec-lyap-i}.

We formalize this physically intuitive principle with the second main Wasserstein-type stability theorem this article, which takes basically the following form.
\begin{theo}\label{theo-intro-II}
 Assume that $W$ is locally bounded and bounded from below by a real number and $(U,V)$ are locally bounded with compact sub-level sets.
In addition, there exists some $\delta\in ]0,1[$ such that the functions $
(\Ua_{\delta/2},
\Va_{\delta/2})
$ are locally bounded with compact sub-level sets and
 we have  $\lambda_{(1-\delta)U}(1)\wedge \nu_{(1-\delta)V}(1)<\infty$.
Then, there exist some $\rho\in ]0,1[$, $n_0\geq 1$ and $c>0$ such that for any $n\geq n_0$,
\begin{equation}\label{theo-intro-e2}
\psi:=\exp{(\delta V)}\Longrightarrow
\vertiii{\pi_{2n}-\nu_V}_{\psi}\leq c~\rho^{n}\quad\mbox{and}\quad
\left\Vert\log{({d \pi_{2n}}/{d\nu_V})}\right\Vert_{\psi}\leq c~ \rho^n.
 \end{equation}
When $\psi$ satisfies the lower bound estimate (\ref{gtow-0}) for some semi-metric $\varrho$ then we also have the exponential
Wasserstein estimates
\begin{equation}\label{theo-intro-was}
\Vert \pi_{2n}-\nu_V\Vert_{\tiny tv}\vee
\WW_{\varrho}(\pi_{2n},\nu_V)\leq c~\rho^{n}.
 \end{equation}
 In addition, 
the relative entropy exponential estimates (\ref{theo-intro-ent1}) are satisfied.
Similar estimates also hold for the distributions
$(\pi_{2n+1},\lambda_U)$ with respect to the $\varphi$-norm $\vertiii{\point}_{\varphi}$ associated with  the function $\varphi=\exp{(\delta U)}$.
\end{theo}
The estimates (\ref{theo-intro-e2}) are a consequence of Corollary~\ref{cor-ref-i} and the Wasserstein formulation (\ref{eq-kanto}) of weighted norms. 
The estimates (\ref{theo-intro-was}) are a consequence of (\ref{gtow}) and (\ref{theo-intro-e2}).
More refined strict contraction estimates are established in Corollary~\ref{lcorUV} and Corollary~\ref{cor-ww-txt}.

\subsection{Illustrations and comments}\label{sec-app-com}
\subsubsection*{Some illustrations}
Section~\ref{sec-illustrations} provides a range of examples illustrating the stability results established in Sections~\ref{stab-th-intro} and~\ref{sec-sinkhorn-contract}. These include polynomial growth potentials (Section~\ref{sec-poly}), heavy-tailed marginals on general normed spaces (Section~\ref{heavy-tail-sec}), boundary state-space models (Section~\ref{bst-mod-sec}), semi-compact and semi-discrete settings (Section~\ref{sec-semi-discrete}), as well as finite mixtures of statistical models (Section~\ref{finite-mix-sec}). Section~\ref{sec-lin-gauss} is devoted to linear Gaussian models, where the regularity conditions and Lyapunov criteria admit explicit characterizations.

When the transition potential $W$ is uniformly bounded, weighted total variation norms do not provide additional insight.
In this regime, exponential convergence of Sinkhorn semigroups is well understood through Hilbert’s projective metric and related contraction techniques
\cite{borwein,chen,deligiannidis,franklin}, which already imply exponential decay in total variation.

On bounded subsets $\XX\subset\mathbb{R}^d$, total variation also dominates Wasserstein distances associated with the Euclidean metric. However, the state space $\XX$ is often non-complete and therefore non-Polish with respect to the Euclidean distance, which complicates Wasserstein convergence analysis.
In this context, it is desirable to find a complete distance on the target space that makes the corresponding Wasserstein distance in distribution space complete~\cite{bolley2008separability}.

This issue naturally motivates the introduction of alternative complete metrics. As a simple illustration, consider the Beta marginal  distributions 
 $\lambda_U=\nu_V$ on  $\XX=\YY=]0,1[$ equipped with the Euclidean metric $\varrho_e\left(x_1,x_2\right)=|x_1-x_2|$ and
  defined by
$$
\nu_V(dx)=\frac{1}{6}~x(1-x)~\nu(dx)\quad \mbox{\rm with}\quad \nu(dx):=
1_{]0,1[}(x)~dx.
$$
Applying (\ref{theo-intro-e1}) to the total variation distance (see (\ref{tv-osc})) we readily check the Wasserstein exponential decays
$$
\WW_{\varrho_e}(\pi_{2n},\nu_V)\leq \Vert\pi_{2n}-\nu_V\Vert_{\tiny tv}\leq \rho^n~
\Vert \pi_{0}-\nu_V\Vert_{\tiny tv}.
$$
Although
the open unit interval $\XX=]0,1[$ is not complete for the Euclidean metric, it becomes complete under
$$
\varrho_c\left(x_1,x_2\right):=\left\vert  \frac{1}{x_1}- \frac{1}{x_2}\right\vert+
\left\vert  \frac{1}{1-x_1}- \frac{1}{1-x_2}\right\vert\geq 2~{\varrho_e}(x_1,x_2).
$$
On the other hand, for any parameter  $\delta\in ]0,1[$, if we have
$$
%\begin{array}{l}
\displaystyle \psi(x):=\exp{(\delta V(x))}\geq \frac{6^{\delta}}{2}~\left(
 \frac{1}{x^{\delta}}+
 \frac{1}{(1-x)^{\delta}}\right)
%\\
%\\
$$
then 
$$
\displaystyle  \psi(x_1)+ \psi(x_2)\geq \varrho_{c,\delta}\left(x_1,x_2\right):=c_{\delta}~\varrho_c\left(x_1,x_2\right)^{\delta}
$$
for some $c_{\delta}>0$. In this context, the l.h.s. estimate (\ref{theo-intro-e2}) yields the $\varrho_{c,\delta}$-Wasserstein estimates
$$
2^{\delta}c_{\delta} ~\WW_{\varrho_{e}}(\pi_{2n},\nu_V)\leq \WW_{\varrho_{c,\delta}}(\pi_{2n},\nu_V)\leq
\vertiii{\pi_{2n}-\nu_V}_{\psi}\leq c~\rho^{n}
$$
as well as the weighted total variation norm decays
$$
\Vert\pi_{2n}-\nu_V\Vert_{\tiny tv}\leq c_1~
\vertiii{\pi_{2n}-\nu_V}_{\psi}\leq c_2~\rho^{n}
\quad \mbox{\rm
for some $c_1,c_2>0$.}
$$

Theorem~\ref{theo-intro-II} also applies to Gaussian marginal distributions on $\YY=\XX=\RR^d$ and a transition potential $W$ with at most quadratic growth (cf. Section~\ref{sec-poly}). In this context,   for $(p,q)=(2,1)$ in (\ref{intro-w-p-q}),
we obtain exponential convergence in the the $2$-Wasserstein distance,
$$
\WW_2(\pi_{2n},\nu_V)\leq c_2~\rho_2^n\quad \mbox{\rm for some parameters $c_2>0$ and $\rho_2\in ]0,1[$.}
$$

Section~\ref{sec-lin-gauss} provides an explicit description of the functions $(\Ua_{\delta},
\Va_{\delta})$ in the context of linear-Gaussian models. For instance, for one-dimensional centered Gaussian marginals with unit variance and a symmetric Gaussian reference transition $\Qa$ with variance $\tau>2$ we have
$$
{2}/{\tau}<\delta<1\Longrightarrow
\Ua_{\delta/2}(x)=c_{\delta}+\frac{x^2}{2}~\left(\frac{\delta}{2}-\frac{1}{\tau}\right)\stackrel{\vert x\vert\to\infty}{\longrightarrow} +\infty
\quad \mbox{\rm for some $c_{\delta}\in\RR$. }
$$
A symmetric Gaussian reference transition $\Qa$ with variance $\tau$ can be interpreted as the transition probability of a Brownian motion on a time horizon $\tau$. In this context, the exponential stability of Sinkhorn bridges w.r.t. the relative entropy for log-concave marginals at infinity, including Gaussian marginals, also holds for any time horizon $\tau>0$, see for instance~\cite{chiarini,durmus,del2025stability}. The article~\cite{adm-24} also provides Wasserstein exponential decays for any Gaussian marginals and any linear-Gaussian reference transition. For these classes of log-concave-type models, our Lyapunov approach yields new weighted total variation norm estimates for sufficiently large time horizons.

\subsubsection*{Energy barriers}
Sections~\ref{stab-th-intro} and~\ref{sec-illustrations} demonstrate that Lyapunov techniques provide a flexible tool for establishing quantitative stability of Sinkhorn semigroups. Nevertheless, several important limitations remain.
 
Firstly, note that Theorem~\ref{theo-intro-II} applies when $W$ is locally bounded.
But even well defined absolutely continuous reference transitions  on totally bounded state spaces such as $\XX=\YY=]0,1[$ may have a potential $W$ which is not well defined on the whole space. Indeed, 
 consider the reference Markov transition defined for any $x\in ]0,1[$ by
$$
\begin{array}{l}
\displaystyle\Qa(x,dy)=\frac{1}{q(x)}~\vert x-y\vert~\nu(dy)\\
\\
\displaystyle\mbox{\rm with}\quad \nu(dy):=1_{]0,1[}(y)~ dy
\quad \mbox{\rm and}\quad \frac{1}{4}\leq q(x)=\frac{1}{2}-x(1-x)\leq \frac{1}{2}.
\end{array}$$
This rather elementary example is a a slight modification of (\ref{exintro}) but
in this case the function $W$ is not defined on the diagonal and for any $x\not=y$ we have
$$
\begin{array}{l}
\displaystyle W(x,y)=-\log{\vert x-y\vert}+\log{q(x)}\Longrightarrow W(x,x)=+\infty.
\end{array}$$
To the best of our knowledge, the stability analysis of Sinkhorn bridges in the context of energy-barrier cost functions of this form remains an  open research question.
Absolutely continuous target marginals may also have an infinite potential, for instance 
$$
\begin{array}{l}
\displaystyle\lambda_U(dx):=4~\vert x-1/2\vert~\lambda(dx)\quad \mbox{\rm with}\quad \lambda(dx):=1_{]0,1[}(x)~ dx\\
\\
\displaystyle \text{~~yields~~}
U(x)=-\log{4}-\log{\vert x-1/2\vert}\Longrightarrow U(1/2)=+\infty.
\end{array}$$
In this case the potential function $U$ is not well defined on the whole state space, thus this marginal model cannot be written in terms of a Boltzmann-Gibbs measure of the form
(\ref{ref-intro-UV}). Nevertheless for uniformly bounded transition potential functions,
Theorem~\ref{theo1-intt} applies to {\it any} class of marginal distributions.

Last but not least, note that the Weibull marginal distribution with shape parameter strictly less than one is locally bounded on $\XX=]0,\infty[$ but does not have compact sub-level sets; for instance, we have
$$
\lambda_U(dx)=\frac{1}{2\sqrt{x}}~e^{-\sqrt{x}}~
1_{]0,\infty[}(x)~dx\Longrightarrow
U(x)=\sqrt{x}+\log{\sqrt{x}}+\log{2}\stackrel{x\to 0^+}{\longrightarrow} -\infty.
$$
The potential function $U$ above is well defined on the whole state space
 but does not have compact sub-level sets.  Even for bounded transition potential functions Theorem~\ref{theo-intro-II} cannot be used (but as mentioned above, Theorem~\ref{theo1-intt} applies to any marginal models).

%%%%%%%%%%%%%%%%%%%%%
%
%%%%%%%%%%%%%%%%%%%%%
\section{Contraction of Markov operators}\label{contract-Markov-sec}

Section~\ref{metrics-sec} is dedicated to a brief description of the Hilbert (projective) metric and integral operator norms. 
Section~\ref{ccoef-sec} recalls some definitions of contraction coefficients for Markov operators in the context of Hilbert projective metrics, weighted total variation norms and $\Phi$-entropies. These coefficients will be key to the analysis of contraction properties of Sinkhorn semigroups in Section \ref{sec-sinkhorn-contract}.

%%%
%
%%%
\subsection{Metrics and divergences}\label{metrics-sec}

%%%
\subsubsection{Hilbert projective metric}

 The Hilbert (projective) metric is defined, for any comparable measures $\nu_1,\nu_2\in \Ma_1(\XX)$, by the formula
\begin{equation}\label{defi-H-ineq}
 H(\nu_1,\nu_2):=\log{
    \left(
        \left\| \frac{d\nu_1}{d\nu_2} \right\|
        \left\| \frac{d\nu_2}{d\nu_1} \right\|
    \right)
} \geq (\log{3})~\Vert \nu_1-\nu_2\Vert_{\sf tv},
\end{equation}
where $\| \cdot\|$ and $\| \cdot \|_{\sf tv}$ denote the supremum and total variation norms, respectively, and $H(\nu_1,\nu_2)=\infty$ when the measures $\nu_1,\nu_2$ are not comparable. The proof of the r.h.s. estimate in (\ref{defi-H-ineq}) is provided in~\cite{atar} (cf. Lemma 1), see also~\cite{legland}. For a more thorough discussion on Hilbert projective metrics we refer the reader to the recent review and self contained article on Hilbert projective metrics~\cite{cohen}.  

%%%
\subsubsection{Integral operator norms}

 Given some functions $\varphi\in \Ba(\XX)$ such that $\varphi\geq 1/2$, let $
 \Ma_{\varphi}(\XX)\subset  \Ma_b(\XX)
$ be the subspace of measures $\nu\in  \Ma_b(\XX)$   equipped with the operator $\varphi$-norm
\begin{eqnarray}
\vertiii{\nu}_{\varphi}&:=&\sup\{|\nu(f)|~:~\Vert f\Vert_{\varphi}\leq 1\}=\sup\{|\nu(f)|~:~\mbox{\rm osc}_{\varphi}(f)\leq 1\}\nonumber\\
&=&\vert\nu|(\varphi)\quad \mbox{\rm with the total variation measure}\quad
\vert\nu \vert:=\nu_++\nu_-\label{def-phi-norm}.
\end{eqnarray}
In the above display $\nu=\nu_+-\nu_-$ stands for a Hahn-Jordan decomposition  of the measure, while $\mbox{\rm osc}_{\varphi}(f)$ stands for the $\varphi$-oscillation of a function $f\in\Ba(\XX)$, which is given by 
$$
\mbox{\rm osc}_{\varphi}(f):=\sup_{x_1,x_2}\frac{\vert f(x_1)-f(x_2)\vert}{\varphi(x_1)+\varphi(x_2)}\leq \Vert f\Vert_{\varphi}.
$$
For a detailed proof of the equivalent formulations in the latter definition \eqref{def-phi-norm}, we refer to Proposition 8.2.16 in~\cite{penev}.  
The weighted Wasserstein formulation (\ref{eq-kanto}) of the $\varphi$-norm is a consequence of
Kellerer’s extension of the Kantorovich-Rubinstein theorem~\cite{edwards2011kantorovich,kellerer1984duality,kellerer1982duality}, see also the recent article~\cite{delgerber2025}.

 The choice of condition $\varphi\geq 1/2$ in the definition of the $\varphi$-norm given in (\ref{def-phi-norm}) is imposed only to recover the conventional total variation dual distance $\vertiii{\mu_1-\mu_2}_{1/2}=\Vert\mu_1-\mu_2\Vert_{\sf tv}$ between probability measures $\mu_1,\mu_2\in \Ma_{1}(\XX)= \Ma_{1,\frac{1}{2}}(\XX)$, where $ \Ma_{1,\varphi}(\XX)$ denotes the subset of probability measures $\nu$ s.t. $\nu(\varphi)<\infty$. 
  {Recall that, given an arbitrary Polish space $S$,} for any $\mu_1,\mu_2\in \Ma_1(S)$ and $\epsilon\in ]0,1]$ we have
\begin{equation}\label{ref-coupling-tv}
\Vert \mu_1-\mu_2\Vert_{\sf tv}\leq 1-\epsilon
\quad\text{if, and only if,}\quad
\exists \nu\in \Ma_1(S)~:~\mu_1\geq \epsilon\nu\quad\mbox{\rm and}\quad \mu_2\geq \epsilon\nu.
\end{equation}
Choosing $\varphi\geq 1$ and $\kappa>0$, for any $\nu\in  \Ma_b(\XX)$ we also have  the norm equivalence formulae
\begin{equation}\label{norm-equivalence}
\begin{array}{l}
\displaystyle\varphi_{\kappa}:=1+\kappa \varphi
\\
\\
\displaystyle\Longrightarrow
\kappa~ \vertiii{\nu}_{\varphi}\leq  \vertiii{\nu}_{\varphi_{\kappa}}\leq (1+\kappa)\vertiii{\nu}_{\varphi}\quad \mbox{\rm and}\quad \vertiii{\nu}_{\varphi_{\kappa}}\geq \vertiii{\nu}_{1}\geq  \Vert\nu\Vert_{\sf tv}.
\end{array}
\end{equation}
For exponential functions of the form
$$
\varphi(x)=c~e^{g(x)} \mbox{\rm ~~with some parameter $c\geq 0$ and~~} 
g(x_1)+g(x_2)\geq \varrho_{g}(x_1,x_2)
$$
for some semi-metric $\varrho_{g}$  we also have the estimate
$$
\varphi(x_1)+\varphi(x_2)\geq 2c~ e^{g(x_1)/2}e^{g(x_2)/2}\geq \varrho_{\varphi}(x_1,x_2):=2c~\left(e^{\frac{1}{2}\varrho_{g}(x_1,x_2)}-1\right).
$$
Also note that for any $p\geq 1$ we have
$$
\varrho_{\varphi}(x_1,x_2)\geq c_p~\varrho_{g}(x_1,x_2)^p
\quad \mbox{\rm with}\quad c_p:= \frac{c}{2^{p-1}p!}
$$
and this implies that
\begin{equation}\label{eq-kanto-estimate}
 \frac{c}{2^{p-1}p!}~\WW_{p,\varrho_g}(\mu_1,\mu_2)^p
\leq \WW_{\varrho_\varphi}(\mu_1,\mu_2) \leq\vertiii{\mu_1-\mu_2}_{\varphi}.
\end{equation}
Last but not least, on normed vector spaces consider a function of the form 
$$
\varphi(x)\geq a+\Vert b(x)\Vert^p
$$
for some function $b(x)$ and $a,p>0$. In this case  we have
\begin{equation}\label{ref-poly-lyap}
\varphi(x_1)+\varphi(x_2)
\geq \varrho(x_1,x_2)\quad \mbox{\rm with}\quad
\varrho(x_1,x_2):=\frac{1}{2^{(p-1)_+}}~\Vert b(x_1)- b(x_2)\Vert^p
\end{equation}
as well as the Wasserstein estimate (\ref{gtow}).

%%%
%
%%%
\subsection{Some contraction coefficients}\label{ccoef-sec}
 
Let $\delta_x$ denote the Dirac delta measure located at $x$. Whenever
$\delta_{x_1}\Ka\ll\delta_{x_2}\Ka$ for any $(x_1,x_2)\in  \XX^2$, we denote by $\imath_{x_1,x_2}(\Ka)(y)$ the Radon-Nikodym function 
\begin{equation}\label{i-def}
\imath_{x_1,x_2}(\Ka)(y):=\frac{d\delta_{x_1}\Ka}{d\delta_{x_2}\Ka}(y).
\end{equation}
We also consider the $[0,1]$-valued parameters
\begin{eqnarray}
\jmath(\Ka)&:=&\inf_{(x_1,x_2)\in  \XX^2}~\inf_{y\in\YY}~\imath_{x_1,x_2}(\Ka)(y) \quad \text{and} \nonumber\\
\hbar(\Ka)&:=&\inf_{(x_1,x_2)\in  \XX^2}\inf_{(y_1,y_2)\in  \YY^2}
\left(\imath_{x_1,x_2}(\Ka)(y_1)~
\imath_{x_2,x_1}(\Ka)(y_2)\right)\geq \jmath(\Ka)^2,\label{it-Ka-def}
\end{eqnarray}
and we note that 
$$
 \delta_{x_1}\Ka (dy_1)~ \delta_{x_2}\Ka (dy_2)~\geq \hbar(\Ka)~\delta_{x_2}\Ka(dy_1)~ \delta_{x_1}\Ka (dy_2).
$$
 Integrating the above inequality  w.r.t. the variable $y_2$ we readily see that
\begin{equation}\label{jhj}
 \delta_{x_1}\Ka (dy_1)~\geq \hbar(\Ka)~\delta_{x_2}\Ka(dy_1),\quad \text{which readily yields}\quad
\jmath(\Ka)\geq \hbar(\Ka)\geq  \jmath(\Ka)^2.
\end{equation}

%%% 
\subsubsection{Birkhoff contraction coefficient}\label{birk-sec}

The Birkhoff contraction coefficient of a Markov transition $\Ka$ from $\XX$ into $\YY$ is defined by
\begin{equation}\label{defi-tau-intro}
\cchi_H(\Ka):=\sup_{}\frac{H(\nu_1\Ka,\nu_2\Ka)}{H(\nu_1,\nu_2)}=\frac{1-\sqrt{\hbar(\Ka)}}{1+\sqrt{\hbar(\Ka)}},
\end{equation}
where the supremum is taken over all comparable probability measures $\nu_1,\nu_2\in\Ma_1(\XX)$ such that $H(\nu_1,\nu_2)>0$ and the second equality is explicitly derived in Appendix \ref{apHilbert}. Combining \eqref{jhj} and \eqref{defi-tau-intro} we readily obtain the equivalent inequalities
$$
{\cchi_H(\Ka)<1}\Longleftrightarrow
\hbar(\Ka)>0\Longleftrightarrow \jmath(\Ka)>0.
$$
We also emphasize that even on totally bounded state spaces  there is an interesting twist because the uniform minorization condition $\jmath(\Ka)>0$ is generally not satisfied. We illustrate this assertion with the semi-circle transition on the unit interval $\XX=]0,1[=\YY$ defined by
\begin{equation}\label{semicircle-ex}
\Ka(x,dy)=\frac{1}{q(x)}~
\sqrt{1-(x-y)^2}~1_{]0,1[}(y)~dy
\end{equation}
for some normalizing constant $1/3\leq q(x)\leq 1$.
In this case, choosing $x_1=1/(2n)=1-y=1-x_2$ in (\ref{i-def}) we check that
$$
0\leq\jmath(\Ka)\leq 3~\sqrt{1-(1-1/n)^2}\stackrel{n\to\infty}{\longrightarrow}0\Longrightarrow
\jmath(\Ka)=0.
$$
The entropic transport model (\ref{def-entropy-pb-v2})  on the unit interval 
$]0,1[$ with reference transition (\ref{semicircle-ex}) is discussed in Section~\ref{beta-w-sec}. Note that in that case, the cost function is unbounded.

\subsubsection{Dobrushin contraction coefficients}

Consider functions $\varphi\in \Ba(\XX)$ and $\psi\in \Ba(\YY)$ such that $\varphi(x)\wedge\psi(y)\geq 1$ for any $(x,y)\in (\XX\times\YY)$.
The $(\varphi,\psi)$-Dobrushin coefficient $\cchi_{\varphi,\psi}(\Ka)$ of a Markov transition $\Ka(x,dy)$ from $\XX$ into $\YY$
 is defined by the norm operator
$$
 \chi_{\varphi,\psi}(\Ka)=\sup_{\nu_1,\nu_2\in \Ma_{1,\varphi}(S)}\frac{\vertiii {\nu_1\Ka-\nu_2\Ka}_{\psi}}{
\vertiii {\nu_1-\nu_2}_{\varphi}}.
$$
When $(\varphi,\XX)=(\psi,\YY)$ we may also write $ \chi_{\varphi}(\Ka)$ instead of $ \chi_{\varphi,\varphi}(\Ka)$. Given a Markov transition $\Ka_0(x,dy)$ from $\XX$ into $\YY$ and a Markov transition 
  $\Ka_1(x,dy)$ from $\YY$ into $\XX$ we have
  \begin{equation}\label{norm-product}
   \chi_{\varphi}(\Ka_0\Ka_1)=   \chi_{\varphi,\varphi}(\Ka_0\Ka_1)\leq  \chi_{\varphi,\psi}(\Ka_0)~ \chi_{\psi,\varphi}(\Ka_1).
  \end{equation}
The terminology $\varphi$-Dobrushin coefficient comes from the fact that  we recover the standard Dobrushin contraction coefficient $\cchi(\Ka)$ by choosing constant functions. Indeed, if $\varphi=1=\psi$, then
\begin{equation}
  \cchi_{1,1}(\Ka) =\sup_{\nu_1,\nu_2\,\in \Ma_1(\XX)}\frac{\Vert\nu_1\Ka-\nu_2\Ka\Vert_{\sf tv}}{\Vert\nu_1-\nu_2\Vert_{\sf tv}}=: \cchi(\Ka) = 1-\varepsilon(\Ka),
\end{equation}
with the parameter
\begin{equation}\label{def-epsil}
\varepsilon(\Ka):=1-\sup_{(x_1,x_2)\,\in \XX^2}{\Vert(\delta_{x_1}-\delta_{x_2})\Ka\Vert_{\sf tv}}.
\end{equation}
Whenever
$\delta_{x_1}\Ka\ll\delta_{x_2}\Ka$ for any $(x_1,x_2)\in  \XX^2$  using the estimate (\ref{jhj}) as well as (\ref{ref-coupling-tv}) we readily find that
\begin{equation}\label{def-epsil-jh}
 \varepsilon(\Ka)\geq \jmath(\Ka)\geq \hbar(\Ka)\quad \mbox{\rm and}\quad
  \cchi(\Ka)\leq 1-\jmath(\Ka)\leq 1- \hbar(\Ka).
\end{equation}

%%% 
\subsubsection{Entropy contraction coefficient}

Given a convex function  $\Phi$ satisfying (\ref{Phi-ref}), the Dobrushin $\Phi$-contraction coefficient $\cchi_{\Phi}(\Ka)$ of  some Markov transition $\Ka(x,dy)$  from $\XX$ into $\YY$  is defined by 
\begin{equation}\label{defi-beta-Phi}
 \cchi_{\Phi}(\Ka):=\sup_{\nu_1,\nu_2\,\in \Ma_1(\XX)}\frac{D_{\Phi}(\nu_1\Ka,\nu_2\Ka)}{D_{\Phi}(\nu_1,\nu_2)},
\end{equation}
where the supremum is taken over all probability measures 
$\nu_1,\nu_2\in\Ma_1(S)$ such that
$D_{\Phi}(\nu_1,\nu_2)>0$. 
 Here, again, the terminology $\Phi$-Dobrushin coefficient comes from the fact that choosing $\Phi(u,v)=\Phi_0(u,v)=\frac{1}{2}|u-v|$ we recover the standard Dobrushin contraction coefficient, i.e., $\cchi_{\Phi_0}(\Ka)= \cchi(\Ka)$. 
For any convex function $\Phi$ satisfying (\ref{Phi-ref}) we recall that
\begin{equation}\label{H-contract}
\cchi_{\Phi}(\Ka)\leq  \cchi(\Ka).
\end{equation}
The estimate (\ref{H-contract}) ensures that the  Dobrushin coefficient is a universal contraction coefficient (see~\cite{cohen-0} for finite state spaces and~\cite{dm-03} for general measurable spaces). For further details on the Dobrushin contraction coefficient we refer the reader to the article~\cite{dob-56}, see also~\cite{dg-01,dm-03} and Chapter 4 in~\cite{dm-04}.

%%%%%%%%%%%%%%%%%%%%%%%%%%%%%
% Sinkhorn contraction 		%
% properties				%
%%%%%%%%%%%%%%%%%%%%%%%%%%%%%
\section{Sinkhorn contraction properties}\label{sec-sinkhorn-contract}
This section investigates the contraction and stability properties of Sinkhorn iterations and their associated Gibbs-loop semigroups. We first recall the strong-mixing regime, in which uniform positivity of the reference kernel yields Hilbert-metric contractions and exponential convergence. We then move beyond this setting by introducing Lyapunov drift and local minorization conditions tailored to the forward-backward Sinkhorn transitions. These conditions allow us to establish weighted-norm contraction inequalities, from which we derive exponential convergence rates in total variation, weighted norms, Wasserstein distances, and relative entropy. The results provide a unified framework for analyzing Sinkhorn dynamics on non-compact state spaces and for treating models that fall outside the classical strong-mixing assumptions.
\subsection{Strong mixing models}\label{strong-mixing-sec}
As underlined in \cite{fienberg,mosteller} in the context of finite state space models, Sinkhorn iterations of the form (\ref{s-2}) preserve the cross product Radon-Nikodym ratio in (\ref{it-Ka-def}) and, hence, for any $n\geq 0$ we have
\begin{equation}\label{cp-inv}
\hbar(\Sa_n)=\hbar(\Sa_0)=\hbar(\Qa).
\end{equation}
This property is crucial in the convergence analysis of Sinkhorn iterations using Hilbert-projective techniques, see for instance~\cite{borwein,chen,deligiannidis,franklin}. Recall that the Hilbert contraction coefficient of a Markov transition $\Qa$ can be explicitly computed in terms of  $\hbar(\Qa)$ (see Eq. \eqref{defi-tau-intro}), hence the exponential convergence of Sinkhorn iterations w.r.t. the Hilbert metric takes place if, and only if, $\hbar(\Qa)>0$. An elementary proof of this property is presented in Proposition~\ref{prop-h-bar-Kn} (in Appendix \ref{apHilbert}). 
From (\ref{jhj}), let us observe that
$$
\jmath(\Sa_n)\wedge\jmath(\Qa)\geq \hbar(\Sa_n)= \hbar(\Qa)\geq  \jmath(\Qa)^2.
$$
Thus, for any given convex function  $\Phi$ satisfying (\ref{Phi-ref}), we can use (\ref{defi-tau-intro}) and (\ref{H-contract}) to show that
$$
\hbar(\Qa)>0 \quad \text{if and only if} \quad \jmath(\Qa)>0,
$$
while 
\begin{equation}\label{def-cchi-S}
\jmath(\Qa)>0~~\text{implies}~~
\left\{
    \begin{array}{l}
    \cchi_H(\Sa_n)=\cchi_H(\Qa)<1,\\
    \\
    \text{$\Da_{1,1}(\Sa^{\circ})$ is satisfied with $\varepsilon(r)= \hbar(\Qa)>0$ for any $r\geq 1$, and}\\
    \\
    \sup_{n\geq 0}\cchi_{\Phi}(\Sa_n)\leq \cchi(\Sa):=\sup_{n\geq 0}\cchi(\Sa_n)\leq 1-\hbar(\Qa)<1.
    \end{array}\right.
\end{equation}
We can also obtain the sandwich-type contraction inequalities described below.
\begin{theo}\label{tx-st}
For any convex function $\Phi$ satisfying (\ref{Phi-ref}) we have the contraction estimates
\begin{equation}\label{ref-contract-even-i}
\begin{array}{rcccl}
 D_{\Phi}\left(\pi_{2(n+1)},\nu_V\right)&\leq& \cchi(\Sa)~ D_{\Phi}\left(\lambda_U ,\pi_{2n+1}\right)&\leq&  \cchi(\Sa)^2~D_{\Phi}\left(\pi_{2n},\nu_V \right),\\
 &&&&\\
  D_{\Phi}\left(\nu_V,\pi_{2(n+1)}\right)&\leq&\cchi(\Sa)~ D_{\Phi}\left(\pi_{2n+1},\lambda_U \right)&\leq &\cchi(\Sa)^2~D_{\Phi}\left(\nu_V ,\pi_{2n}\right).
 \end{array}
 \end{equation}
The above  inequalities remain valid when we replace $(\nu_V,\pi_{2n})$ by  $(\lambda_U,\pi_{2n+1})$. If we assume that $\jmath(\Qa)>0$, then $\cchi(\Sa)<1$ and, in addition, there exists some $c>0$ such that
 \begin{equation}\label{log-est-sm} 
\left\Vert\log{\frac{d \pi_{2n}}{d\nu_V}}\right\Vert\vee
\left\Vert\log{\frac{d \pi_{2n+1}}{d\lambda_U}}\right\Vert\leq c~ \cchi(\Sa)^{2n}
\quad \text{for any $n\geq 0$.}
\end{equation}
\end{theo}
\proof
The estimates (\ref{ref-contract-even-i}) are a direct consequence of (\ref{H-contract})  combined with the Markov transport formulae
 \begin{eqnarray*}
( \pi_{2(n+1)},\nu_V)&=&(\lambda_U \Sa_{2(n+1)},\pi_{2n+1}\Sa_{2(n+1)}), \quad\text{and}\\
(\lambda_U ,\pi_{2n+1})&=&(\pi_{2n}\Sa_{2n+1},\nu_V \Sa_{2n+1}),
 \end{eqnarray*}
which are themselves easily obtained from expressions (\ref{s-2}) and (\ref{fixed-points-gibbs}). The proof of (\ref{log-est-sm}) follows exactly the same arguments as the proof of the inequality (\ref{log-est}) in Theorem \ref{lm-tx-2}, hence it is skipped. \cqfd

The following result can be obtained from the series expansion given by Corollary~\ref{cor-1-i}.
\begin{cor}\label{cor-entropy-wb}
If we assume that $\jmath(\Qa)>0$, then $\cchi(\Sa)<1$ and, in addition, there exists some $c>0$ such that
$$
  \mbox{\rm Ent}(\PP~|~\Pa_{2n})  \leq   c~ \cchi(\Sa)^{2n}, \quad \text{for any $n\ge 0$}.
$$
\end{cor}

Several examples satisfying condition $ \jmath(\Qa)>0$ are discussed in Section 2.4 in~\cite{dg-01}. For instance, consider a parameter $\delta>0$ and drift function $x\mapsto b(x)$ with bounded oscillations from  $\XX=\RR^d$ into itself. The choice $\Qa(x,dy)~\propto~e^{-\delta\,\Vert y-b(x)\Vert}~dx$ yields
$$
\jmath(\Qa)\geq \exp{\left(-\delta\,\Delta(b)\right)},
\quad \text{with}\quad\Delta(b)=\sup_{(x_1,x_2)\in \XX^2}{\Vert b(x_1)-b(x_2)\Vert}.
$$
In continuous time settings, condition $ \jmath(\Qa)>0$ is also satisfied for Markov transitions of 
elliptic diffusions on compact manifolds $\XX$ --see for instance the pioneering work of Aronson~\cite{aronson}, Nash~\cite{nash} and Varopoulos~\cite{varopoulos} on Gaussian estimates for heat kernels on manifolds. Nevertheless, the above condition is not met for  Gaussian transitions.

\subsection{Gibbs-Loop Sinkhorn semigroups}\label{sec-sinkhorn-contract-reg}

For a positive function $\varphi\in \Ba(\XX)$ and $r>0$, we define its $r$-sublevel set by
\begin{equation}\label{def-sublevel}
C_{\varphi}(r,\XX):=\left\{x\in \XX~:~\varphi(x)\leq r\right\}.
\end{equation}
When there is no ambiguity, we write $C_{\varphi}(r)=\{\varphi\leq r\}$.

Lyapunov functions combined with local minorization conditions are classical tools in the analysis of time-inhomogeneous Markov chains; see, for instance, \cite{dpa} and the references therein. In the context of Gibbs-loop semigroups \eqref{gibbs-tv}, such conditions are formulated using Lyapunov-type functions $\varphi\in\Ba(\XX)$ and $\psi\in\Ba(\YY)$ satisfying $\varphi\wedge\psi\geq 1$, as described below.

\medskip
\noindent
\textbf{Condition} {\it $\Da_{\varphi,\psi}(\Sa^{\circ})$:}
\textit{
There exist $\epsilon\in(0,1)$, a constant $c>0$, and an integer $n_0\geq 1$ such that
$\pi_{2n_0}(\psi)\vee \pi_{2n_0+1}(\varphi)<\infty$ and, for all $n\geq n_0$,
\begin{equation}\label{lyap-eq-sink-o}
\Sa_{2n}^{\circ}(\psi)\leq \epsilon,\psi+c,
\qquad
\Sa_{2n+1}^{\circ}(\varphi)\leq \epsilon,\varphi+c.
\end{equation}
Moreover, there exist $r_0\geq 1$ and a function
$$
\varepsilon: [r_{0},\infty[ \rightarrow ]0,1]
$$
such that, for any $r\geq r_0$, $(x_1,x_2)\in C_{\varphi}(r)^2$, $(y_1,y_2)\in C_{\psi}(r)^2$, and $n\geq n_0$,
\begin{equation}\label{sink-locmin-o}
\Vert (\delta_{x_1}-\delta_{x_2})\Sa^{\circ}_{2n+1}\Vert{\sf tv}
\vee
\Vert (\delta_{y_1}-\delta_{y_2})\Sa^{\circ}_{2n}\Vert{\sf tv}
\leq 1-\varepsilon(r).
\end{equation}
}

Condition \eqref{sink-locmin-o} holds if and only if there exist probability measures
$\Xa\in\Ma_1(\XX)$ and $\Ya\in\Ma_1(\YY)$, possibly depending on $(n,r,(x_1,x_2))$ and $(n,r,(y_1,y_2))$, respectively, such that, for $i\in{1,2}$,
  \begin{equation}\label{sink-locmin-bis-o}
\Sa^{\circ}_{2n}(y_i,dy)\geq \varepsilon(r)~\Ya(dy)\quad \mbox{\rm and}\quad
\Sa_{2n+1}^{\circ}(x_i,dx)\geq \varepsilon(r)~\Xa(dx).
 \end{equation}
 The drift condition \eqref{lyap-eq-sink-o} is trivially satisfied for the constant functions
$\varphi=1=\psi$. In this case, $C_{\varphi}(r)=\XX$ and $C_{\psi}(r)=\YY$ for all $r\geq 1$, and
\eqref{sink-locmin-bis-o} reduces to a uniform minorization condition over the whole state spaces.

Let $\Ba_{\infty}(\XX)\subset\Ba(\XX)$ denote the subalgebra of locally bounded, uniformly positive functions $\varphi$ that grow at infinity, in the sense that $\sup_C\varphi<\infty$ for any compact set $C\subset\XX$ and that, for any
\begin{equation}
r\geq \varphi_\star:=\inf_\XX\varphi>0,
\label{eqPhi_star}
\end{equation}
the sublevel set $C_{\varphi}(r)$ is nonempty and compact.
We further define
 $$\Ba_0(\XX):=\{1/\varphi~:~\varphi\in \Ba_{\infty}(\XX)\}\subset \Ba_b(\XX)
 $$ 
 the algebra of bounded, positive functions that are locally lower bounded and vanish at infinity.
If $\varphi\in\Ba_{\infty}(\XX)$ and $\psi\in\Ba_{\infty}(\YY)$, then the sets
$C_{\varphi}(r)$ and $C_{\psi}(r)$ are compact, and \eqref{sink-locmin-bis-o} yields genuine local minorization conditions.

The assumption $\varphi\wedge\psi\geq 1$ is made solely for notational convenience.
Indeed, if $\varphi\in\Ba_{\infty}(\XX)$ and $\psi\in\Ba_{\infty}(\YY)$ satisfy
\eqref{lyap-eq-sink-o}, then the rescaled functions
$\varphi/\varphi_\star$ and $\psi/\psi_\star$ also satisfy the same drift inequalities.

The drift conditions \eqref{lyap-eq-sink-o} imply that, for all $n\geq n_0$,
 $$
 \Sa^{\circ}_{2n}\left(\Ba_{\psi}(\YY)\right)\subset\Ba_{\psi}(\YY)
 \quad\mbox{\rm and}\quad
  \Sa_{2n+1}^{\circ}\left(\Ba_{\varphi}(\XX)\right)\subset \Ba_{\varphi}(\XX).
 $$ 
 Moreover, whenever $\Da_{\varphi,\psi}(\Sa^{\circ})$ holds,
   \begin{equation}\label{uest-pi}
 \lambda_U(\varphi)\vee \nu_V(\psi)<\infty\quad \mbox{\rm and}\quad
 \sup_{n\geq n_0}(\pi_{2n}(\psi)\vee \pi_{2n+1}(\varphi))<\infty.\
\end{equation}
These bounds follow directly from the drift inequalities \eqref{lyap-eq-sink-o}; a detailed proof is given in Section~\ref{tech-proof-ap}.

\subsection{Weighted norm contractions}\label{main-theo-sec}
Combining the drift conditions \eqref{lyap-eq-sink-o} with the local minorization
conditions \eqref{sink-locmin-o} allows us to control the contraction coefficients of the Sinkhorn transitions.
More precisely, we obtain the following uniform contraction estimate for the forward–backward Markov kernels $\Sa^{\circ}_{2n}$ and $\Sa^{\circ}_{2n+1}$.

\newpage

\begin{theo} \label{theorem4.1}
Assume that Condition $\Da_{\varphi,\psi}(\Sa^{\circ})$ holds for some $n_0\geq 1$. Then there exist parameters $\kappa> 0$ and $\rho\in [0,1[$ depending  only on  $(\epsilon,c)$ in (\ref{lyap-eq-sink-o}) and on $\varepsilon(r)$ in (\ref{sink-locmin-o}) and  such that 
\begin{equation}\label{th-1-est}
\sup_{n\geq n_0}\left(
\cchi_{\psi_{\kappa}}(\Sa^{\circ}_{2n})\vee \cchi_{\varphi_{\kappa}}(\Sa^{\circ}_{2n+1})
\right)\leq \rho,
\quad \text{where} \quad 
(\varphi_{\kappa},\psi_{\kappa}):=(1+\kappa\varphi,1+\kappa\psi).
\end{equation}
\end{theo}

Theorem~\ref{theorem4.1} is an immediate consequence of Theorem~\ref{theo-contract-V} in Appendix~\ref{phi-psinorm-sec}.

The contraction estimate (\ref{th-1-est}) ensures the exponential decay of Sinkhorn semigroups. For instance,  for any $\mu_1,\mu_2\in \Ma_{1,\psi}(\YY)$ and $n_0\leq p\leq n$ we readily obtain
$$
\vertiii{ (\mu_1-\mu_2)\Sa^{\circ}_{2(p-1),2n}}_{\psi_{\kappa}}\leq \rho^{1+n-p}~
\vertiii{  \mu_1-\mu_2}_{\psi_{\kappa}}.
$$
Note that
$$
\kappa~\psi\leq \psi_{\kappa}\leq (1+\kappa)\psi
\quad\text{implies}\quad
\displaystyle\vertiii{  (\mu_1-\mu_2)\Sa^{\circ}_{2(p-1),2n}}_{\psi}\leq (1+1/\kappa)~\rho^{1+n-p}~
\vertiii{  \mu_1-\mu_2}_{\psi}.
$$
Choosing $p=n_0$ and $(\mu_1,\mu_2)=(\pi_{2(n_0-1)},\nu_V)$, the above estimate immediately yields, for any $n\geq n_0$, an upper bound on the $\psi$-norm. We summarize the above discussion with the following corollary.
\begin{cor}\label{cor-wnorm}
Under the assumptions of Theorem~\ref{theorem4.1} for any $n\geq n_0$ we have
\begin{eqnarray}
\vertiii{  \pi_{2n}-\nu_V}_{\psi}&\leq& (1+1/\kappa)~\rho^{1+n-n_0}~
\vertiii{  \pi_{2(n_0-1)}-\nu_V}_{\psi}\nonumber
\\
\vertiii{  \pi_{2n+1}-\lambda_U}_{\varphi}&\leq& (1+1/\kappa)~\rho^{1+n-n_0}~
\vertiii{  \pi_{2n_0-1}-\lambda_U}_{\varphi},\label{cor-wnorm-est}
\end{eqnarray}
with the parameters $\rho\in [0,1[$ and $\kappa>0$ as in (\ref{th-1-est}).
\end{cor}

Choosing  $(\mu_1,\mu_2)=(\delta_{y_1},\delta_{y_2})$ for some $(y_1,y_2)\in\YY^2$, we can also verify the merging property of Sinkhorn semigroups \cite{persi,saloff-zuniga,saloff-zuniga-2}
$$
\Vert \Sa^{\circ}_{2(n_0-1),2n}(y_1,\point)-\Sa^{\circ}_{2(n_0-1),2n}(y_2,\point)\Vert_{\sf tv}\leq \rho^{1+n-n_0}~\left(\psi_{\kappa}(y_1)+\psi_{\kappa}(y_2)\right).
$$
\begin{rmk}
The proximal sampler discussed in Remark~\ref{rmk-proximal}  has a Markov transition 
$\Sa^{\circ}_{1}$ with invariant measure $\lambda_U$. The weighted-norm contraction properties  of these models are simply based on the existence of a Lyapunov function $\varphi$ satisfying the r.h.s. condition in (\ref{lyap-eq-sink-o}) with $n=0$. As in (\ref{sink-locmin-o}) we also require a local minorization condition $\Vert (\delta_{x_1}-\delta_{x_2})\Sa^{\circ}_{1}\Vert_{\sf tv}
\leq 1-\varepsilon(r)$ on the sub-level sets $(x_1,x_2)\in  C_{\varphi}(r)^2$ of the Lyapunov function
 to ensure that $\cchi_{\varphi_{\kappa}}(\Sa^{\circ}_{1})<1$, for some parameter $\kappa\geq 0$ (see for instance Theorem~\ref{theo-contract-V}).
\end{rmk}

Condition $\Da_{\varphi,\psi}(\Sa^{\circ})$ is expressed in terms of the forward-backward transition kernels $\Sa^{\circ}_{2n}$ and $\Sa^{\circ}_{2n+1}$, which may not be easy to test -- except in the case of Gaussian and linear models. Several Lyapunov functions in this context are designed in Section \ref{sec-lin-gauss}. A stronger (but easier to test) condition can be provided in terms of the  transitions $\Sa_{2n}$ and $\Sa_{2n+1}$.

\medskip
\noindent\textbf{Condition} {\it $\Da_{\varphi,\psi}(\Sa)$ 
There exists  $\epsilon\in ]0,1[$, a constant $c>0$ and an integer $n_0\geq 0$ such that for any $n\geq n_0$,
 \begin{equation}\label{lyap-eq-sink}
 \Sa_{2n}(\psi)\leq \epsilon~\varphi+c,\qquad
  \Sa_{2n+1}(\varphi)\leq \epsilon~\psi+c.
 \end{equation}
Moreover, there exists some $r_{0}\geq 1$ and a function 
 $$
 \varepsilon: [r_{0},\infty[ \rightarrow ]0,1]
$$ 
such that for any $r\geq r_0$, $(x_1,x_2)\in  C_{\varphi}(r)^2$, $(y_1,y_2)\in  C_{\psi}(r)^2$ and
 $n\geq n_0$,
\begin{equation}\label{sink-locmin}
\Vert (\delta_{x_1}-\delta_{x_2})\Sa_{2n}\Vert_{\sf tv}\vee  \Vert (\delta_{y_1}-\delta_{y_2})\Sa_{2n+1}\Vert_{\sf tv}
\leq 1-\varepsilon(r).
\end{equation}
}

When Condition $\Da_{\varphi,\psi}(\Sa)$ is met, for any $n\geq n_0$ we readily see that  
 $$
 \Sa_{2n}\left(\Ba_{\psi}(\YY)\right)\subset\Ba_{\varphi}(\XX)
 \quad\mbox{\rm and}\quad
  \Sa_{2n+1}\left(\Ba_{\varphi}(\XX)\right)\subset \Ba_{\psi}(\YY).
 $$
 In addition, we have
 \begin{eqnarray}
 \pi_{2n+1}(\varphi)=\nu_V \Sa_{2n+1}(\varphi)&\leq& \epsilon~\nu_V(\psi)+c\label{ref-control-pi}
\\
\pi_{2n}(\psi)=\lambda_U \Sa_{2n}(\psi)&\leq& \epsilon~\lambda_U(\varphi)+c\quad \mbox{\rm  with  $(\epsilon,c)$ as in (\ref{lyap-eq-sink}).}\nonumber
\end{eqnarray}

 \begin{theo}\label{theo-1-intro}
 We have the following:
 \begin{itemize}
\item Condition $\Da_{\varphi,\psi}(\Sa^{\circ})$ holds as soon as condition $\Da_{\varphi,\psi}(\Sa)$ is satisfied and the estimates (\ref{cor-wnorm-est}) hold for some $n_0\geq 1$ and some parameters $\rho\in [0,1[$ and $\kappa>0$. 

\item When  $\Da_{\varphi,\psi}(\Sa)$ is satisfied there exist parameters $\kappa\geq 0$ and  $\rho\in [0,1[$ that only depend on the parameters  $(\epsilon,c)$ in (\ref{lyap-eq-sink}) and on the function $\varepsilon(r)$ in (\ref{sink-locmin}) and such that, for any $n\geq n_0$, 
 $$
\cchi_{\varphi_{\kappa},\psi_{\kappa}}(\Sa_{2n})\vee 
  \cchi_{\psi_{\kappa},\varphi_{\kappa}}(\Sa_{2n+1})
  <\rho,
\quad\text{where}\quad
  (\varphi_{\kappa},\psi_{\kappa}):=
  (1+\kappa\varphi,1+\kappa\psi).
$$
\end{itemize}
 \end{theo}

A proof of Theorem~\ref{theo-1-intro} is provided in Appendix \ref{phi-psinorm-sec} (on page~\pageref{theo-1-intro-proof}).

Let us emphasize that Condition~$\Da_{\varphi,\psi}(\Sa)$ is stronger than
Condition~$\Da_{\varphi,\psi}(\Sa^{\circ})$. In particular, the Gibbs-loop
transitions $\Sa^{\circ}_n$ may be stable even in situations where the
Sinkhorn transitions $\Sa_n$ are unstable. Intuitively, Gibbs-loop
transitions alternate between forward and backward steps that steer the
dynamics toward their invariant measures \eqref{fixed-points-gibbs},
whereas the Sinkhorn transitions in \eqref{s-2} asymptotically transport
one target marginal measure to the other.

As an illustration, Lemma~4.9 in~\cite{adm-24} (see also
Section~\ref{GS-sec}) shows that the Gibbs-loop transitions $\Sa^{\circ}_n$
of the Gaussian Sinkhorn algorithm are always stable. By contrast, in the
setting of regularized Gaussian entropic problems over a finite time
horizon $t$ (see Section~3.4 in~\cite{adm-24}), the Sinkhorn transitions
$\Sa_n$ may become unstable when the time horizon is not sufficiently
large; see Remarks~3.9 and~4.6 in~\cite{adm-24}.

%%%
%
%%%
\subsection{A Lyapunov stability theorem}\label{sec-lyap-i}

\subsubsection{Regularity conditions}
This section presents some sufficient conditions for $\Da_{\varphi,\psi}(\Sa)$, in terms of the potential functions $(U,V,W)$. 
These simple conditions enable us to obtain geometric rates for the convergence of the Sinkhorn iterations for a variety of models, including polynomial growth models (Section~\ref{sec-poly}) and heavy tailed models (Section~\ref{heavy-tail-sec}), as well as boundary states models including Beta and Weibull distributions (Section~\ref{bst-mod-sec} and Section~\ref{beta-w-sec}). 
Illustrations in the context of Linear-Gaussian models and finite  (non necessarily Gaussian) mixture marginals are also discussed respectively in Section~\ref{sec-lin-gauss} and Section~\ref{finite-mix-sec}.

Let us assume that $U,V,W$ are locally bounded and $e^{-W}\in\LL_{\infty}(\lambda\otimes\nu)$. In addition, $U,V$ have compact sub-level sets.
These conditions ensure there exists some $W_{\star}\in\RR$ such that
\begin{equation}\label{ref-Wmin}
W(x,y)\geq W_{\star}\quad \mbox{\rm for $(\lambda\otimes\nu)$-almost every $(x,y)\in(\XX\times\YY)$.}
\end{equation}
 In addition,
$(U,V)$ are bounded below by a real number so that
$e^{-U}$ and $e^{-V}$ are uniformly bounded. This excludes probability densities that tends to infinity at boundary states, such that some classes of Weibull or the Beta densities for some shape parameters.  Also note that for any $\delta>0$ we have
$$
(e^{\delta U},e^{\delta V})\in (\Ba_{\infty}(\XX)\times \Ba_{\infty}(\YY)).
$$

Consider the following conditions:

\medskip
\noindent\textbf{Condition} {\it $\Ha_{\delta}(U,V)$:
There is some $\delta\in ]0,1[$ such that the functions
$
(\Ua_{\delta},
\Va_{\delta})
$ defined in (\ref{UV-delta}) are locally bounded with compact sub-level sets and we have
\begin{equation}\label{huv}
\lambda_{(1-\delta)U}\left(e^{W^V}\right)\vee
~\nu_{(1-\delta)V}\left(e^{W_U}\right)
<\infty.
\end{equation}}

Note that \eqref{ref-Wmin} and \eqref{huv} together imply the inequalities
$$
\lambda_{(1-\delta)U}\left(1\right)\vee \nu_{(1-\delta)V}\left(1\right)\leq e^{-W_{\star}}~\left(\lambda_{(1-\delta)U}\left(e^{W^V}\right)\vee
~\nu_{(1-\delta)V}\left(e^{W_U}\right)\right)<\infty.
$$
This yields the exponential moment estimates
$$
\Ha_{\delta}(U,V)\Longrightarrow
\lambda_{(1-\delta)U}\left(1\right)=\lambda_U\left(\exp{(\delta U)}\right)<\infty\quad \mbox{\rm and}\quad \nu_{(1-\delta)V}(1)=\nu_V\left(\exp{(\delta V)}\right).
$$ 
In particular all absolute moments $\lambda_U(|U|^p)\vee \nu_V(|V|^p)<\infty$, with $p\geq 0$ are finite.
Given the reference measure $\Pa_0=(\lambda_U\times \Qa)$ also have
\begin{equation}\label{entropy-minP0}
\Ha_{\delta}(U,V)\Longrightarrow
\mbox{\rm Ent}((\lambda_U\otimes\nu_V)~|~\Pa_0)=(\lambda_U\otimes\nu_V)(W)-\nu_V(V)<\infty.
\end{equation}
To check this claim,  note that
\begin{eqnarray*}
(\lambda_U\otimes\nu_V)(W)&=&
\lambda_U(W^V)=W_{\star}+\lambda_U(W^V-W_{\star})\\
&\leq &W_{\star}+\lambda_U(e^{W^V-W_{\star}})=
W_{\star}+e^{-W_{\star}}~\lambda_{(1-\delta)U}(e^{-\Ua_{\delta}}).
\end{eqnarray*}
On the other hand, for any $r>0$ we have
\begin{eqnarray*}
\lambda_{(1-\delta)U}(e^{-\Ua_{\delta}})&\leq& e^{-r}~\lambda_{(1-\delta)U}(1)+\lambda_{(1-\delta)U}(e^{-\Ua_{\delta}}~1_{\Ua_{\delta}\leq r})\\
&\leq &\lambda_{(1-\delta)U}(1)~\left(e^{-r}+\exp{\left(\sup_{\Ua_{\delta}(x)\leq r}\vert \Ua_{\delta}(x)\vert\right)}\right)\Longrightarrow (\ref{entropy-minP0}).
\end{eqnarray*}

\medskip
\noindent
\textbf{Condition} {\it $\Ha^{\prime}_{\delta}(U,V)$:
There is some $\delta\in ]0,1/2[$ such that
$
\Ua_{\delta}$ and $
\Va_{\delta}
$ are locally bounded with compact sub-level sets and we have
\begin{equation}\label{huv-prim}
\lambda_{(1-2\delta)U}\left(
1\right)\vee \nu_{(1-2\delta)V}\left(
1\right)<\infty.
\end{equation}}

Let us remark that
\begin{equation}\label{huv-prim-equiv}
\Ha^{\prime}_{\delta}(U,V)\Longrightarrow \Ha_{\delta}(U,V).
\end{equation}
To verify this claim, note that the regularity condition on $(\Ua_{\delta},\Va_{\delta})$ ensures that 
$$
(e^{\Ua_{\delta}},e^{\Va_{\delta}})\in (\Ba_{\infty}(\XX)\times \Ba_{\infty}(\YY)).
$$
When $\Ua_{\delta}$ is locally bounded with compact sub-level sets
$C_{\Ua}(r):=\{\Ua_{\delta}\leq r\}$ for any $0<\delta<1/2$ we have
\begin{eqnarray*}
\lambda_{(1-\delta)U}\left(e^{W^V}\right)&=&\lambda_{(1-2\delta)U}\left(
e^{-\Ua_{\delta}}\right)\\
&\leq &\left(e^{-r}+e^{c_{\delta}(r)}\right)
\lambda_{(1-2\delta)U}\left(
1\right)
 \quad \mbox{\rm with}\quad
c_{\delta}(r):=\sup_{x\in C_{\Ua}(r)} \vert \Ua_{\delta}(x)\vert. 
\end{eqnarray*}
In this context, we readily see that
\begin{equation}\label{huv-2}
\lambda_{(1-2\delta)U}\left(
1\right)\vee \nu_{(1-2\delta)V}\left(
1\right)<\infty\Longrightarrow (\ref{huv}).
\end{equation}
This ends the proof of (\ref{huv-prim-equiv}).

\begin{rmk}\label{rmk-W-bb}
Whenever $\Vert W\Vert<\infty$ we have $\Vert W_U\Vert\vee \Vert W^V\Vert<\infty$.
In this case $\Ha_{\delta}(U,V)$ is satisfied as soon as there is some $\delta\in ]0,1[$ such that
\begin{equation}\label{huv-3}
\lambda_{(1-\delta)U}\left(
1\right)\vee \nu_{(1-\delta)V}\left(
1\right)<\infty.
\end{equation}
\end{rmk}

\subsubsection{A Lyapunov stability theorem}

\begin{theo}\label{lm-tx-2}
If Condition $\Ha_{\delta}(U,V)$ is satisfied for some $\delta\in ]0,1]$, then Condition $\Da_{\varphi,\psi}(\Sa)$ is also satisfied with $n_0=1$ for the Lyapunov functions 
\begin{equation}\label{def-Lyap-UV} 
(\varphi,\psi)=(e^{\delta (U-U_{\star})},e^{\delta (V-V_{\star})})
\quad\mbox{with}\quad
(U_\star,V_\star):=(\inf_\XX U,\inf_{\YY} V)\in(\RR\times\RR).
\end{equation}
 In addition, there exist some parameter $\rho\in [0,1[$ and some constant $c>0$ such that
\begin{equation}\label{log-est} 
\left\Vert\log{\frac{d \pi_{2n}}{d\nu_V}}\right\Vert_{\psi}\vee
\left\Vert\log{\frac{d \pi_{2n+1}}{d\lambda_U}}\right\Vert_{\varphi}\leq c~\rho^{2n}
\quad \text{for any $n\geq 1$.}
\end{equation}
\end{theo}
A proof of Theorem~\ref{lm-tx-2} is provided in Appendix \ref{lm-tx-2-proof}.

 The above theorem combined with  the estimates (\ref{ref-control-pi}) ensures that
 \begin{eqnarray*}
\sup_{n\geq n_0} \pi_{2n+1}(\varphi)&\leq& \epsilon~\nu_V(\psi)+c=\epsilon~ e^{-\delta V_{\star}}~\nu_{(1-\delta)V}(1)+c
\\
\sup_{n\geq n_0} \pi_{2n}(\psi)&\leq& \epsilon~\lambda_U(\varphi)+c=\epsilon~ e^{-\delta U_{\star}}~\lambda_{(1-\delta)U}(1)+c
\end{eqnarray*}
 with the parameters $(\epsilon,c)$ as in (\ref{lyap-eq-sink}).
 The rescaling of the potential functions (\ref{def-Lyap-UV}) is only made to have 
$\varphi\wedge\psi\geq 1$.
Combining Theorem~\ref{lm-tx-2} with Theorem~\ref{theo-1-intro} we obtain the following corollary. 
\begin{cor}\label{lcorUV}
Assume $\Ha_{\delta}(U,V)$ is satisfied for some $\delta\in ]0,1]$. In this situation, there exists  some parameters  $n_0\geq 0$,
$\rho\in [0,1[$ and $c_1,c_2\geq 0$ such that for any $n\geq n_0$ we have
\begin{eqnarray}
 \Vert \pi_{2n}-\nu_V\Vert_{\tiny tv}\leq c_1~ \vertiii{ \pi_{2n}-\nu_V}_{\psi}&\leq& c_2~\rho^{n-n_0}~
\vertiii{  \pi_{2n_0}-\nu_V}_{\psi}
\nonumber\\
\Vert \pi_{2n+1}-\lambda_U\Vert_{\tiny tv}\leq c_1~\vertiii{ \pi_{2n+1}-\lambda_U}_{\varphi}&\leq& c_2~\rho^{n-n_0}~
\vertiii{  \pi_{2n_0+1}-\lambda_U}_{\varphi},\label{cor-wnorm-est-Hd}
\end{eqnarray}
with the  functions 
$
(\varphi,\psi)=(e^{\delta U},e^{\delta V})
$.
\end{cor}
The r.h.s. term in (\ref{cor-wnorm-est-Hd}) can be estimate with (\ref{ref-control-pi}). Also note that
$$
(\ref{def-Lyap-UV})\Longrightarrow 
\lambda_U(\varphi)=e^{-\delta U_{\star}}~\lambda_{(1-\delta)U}(1)\quad\mbox{\rm and}\quad
\nu_V(\psi)=e^{-\delta V_{\star}}~\nu_{(1-\delta)V}(1)
$$
This yields the exponential decays.
\begin{cor}\label{cor-ref-i}
Assume $\Ha_{\delta}(U,V)$ is satisfied for some $\delta\in ]0,1]$. In this situation, there exists  some parameters  $n_0\geq 0$,
$\rho\in [0,1[$ and $c_1,c_2\geq 0$ such that for any $n\geq n_0$ we have
\begin{eqnarray}
 \Vert \pi_{2n}-\nu_V\Vert_{\tiny tv}\leq c_1~ \vertiii{ \pi_{2n}-\nu_V}_{\psi}&\leq& c_2~\rho^{n}~\left(1+\lambda_U(\varphi)+\nu_V(\psi)\right)
\nonumber\\
\Vert \pi_{2n+1}-\lambda_U\Vert_{\tiny tv}\leq c_1~\vertiii{ \pi_{2n+1}-\lambda_U}_{\varphi}&\leq& c_2~\rho^{n}~\left(1+\lambda_U(\varphi)+\nu_V(\psi)\right),\label{cor-wnorm-est-Hd2}
\end{eqnarray}
with the  functions $(\psi,\varphi)$ as in Corollary~\ref{lcorUV}.
\end{cor}

Assume that $(U,V)$ and $(W_{U},W^{V})$ are continuous functions and we have
$$
\delta U-W^{V}\geq \Ua_{\delta}^-
\quad\mbox{\rm and}\quad
\delta V-W_{U}\geq \Va_{\delta}^-
,$$
for some functions $(\Ua_{\delta}^-,\Va_{\delta}^-)$ with compact level sets. In this case, 
the functions $(\Ua_{\delta},\Va_{\delta})$ are locally bounded with compact level sets.
As underlined in~\cite{ajay-dp-25}, using (\ref{ref-dom-theta-0}) and Proposition~\ref{propLyapunovInequalities}, the continuity condition in the above assertion can be relaxed. More precisely, Condition $\Da_{\varphi,\psi}(\Sa)$, Theorem~\ref{lm-tx-2} as well as Corollaries~ \ref{lcorUV} and~\ref{cor-ref-i}
hold
 when $(\Ua_{\delta}^-,\Va_{\delta}^-)$ are locally bounded with compact sub-level sets.
Assume we have $$W(x,y)\leq a_W(x)+b_W(y)$$ for some locally bounded functions $(a_W,b_W)$ such that
$\lambda_U(a_W)\vee\nu_V(b_W)<\infty$.  Note that
$$
\Ua^-_{\delta}(x):=\delta U(x)-a_W(x)-\nu_V(b_V)\leq \Ua_{\delta}(x)\leq \delta U(x)-W_{\star}.
$$
with $W_{\star}$ as in (\ref{ref-Wmin}).
This shows that $(\Ua_{\delta},\Ua^-_{\delta})$ and similarly $(\Va_{\delta},\Va^-_{\delta})$ 
with $\Va^-_{\delta}:=\delta V-b_W$ are locally bounded.  Moreover $(\Ua^-_{\delta},\Va^-_{\delta})$ have compact 
compact sub-level sets if and only if $(\delta U-a_W,\delta V-b_W)$ have compact 
compact sub-level sets.

\subsection{Wasserstein exponential decays}\label{sec-wass}

Assume that  the Lyapunov functions $(\psi,\varphi)$  defined in (\ref{def-Lyap-UV})
satisfy the estimates
\begin{equation}\label{UV-exp2wa}
e^{\delta U(x_1)}+
e^{\delta U(x_2)}\geq \varrho_{\varphi}(x_1,x_2)\quad\mbox{\rm and}\quad
e^{\delta V(y_1)}+
e^{\delta V(y_2)}\geq \varrho_{\psi}(y_1,y_2)
\end{equation}
for some semi-metric $\varrho_{\varphi},\varrho_{\psi}$. In this case, by (\ref{gtow}) 
 there exists  $c_u,c_v>0$ such that
 $$
 \WW_{\varrho_\varphi}(\pi_{2n+1},\lambda_U) \leq c_u~\vertiii{\pi_{2n+1}-\lambda_U}_{\varphi}
 \quad
\mbox{\rm and}\quad
\WW_{\varrho_\psi}(\pi_{2n},\nu_V) \leq\vertiii{\pi_{2n}-\nu_V}_{\psi}
 $$
with the Wasserstein  semi-metrics $ \WW_{\varrho_\varphi}, \WW_{\varrho_\psi}$ defined in (\ref{def-WW-varrho}). This yields the following Wasserstein estimates.

\begin{cor}\label{cor-ww-txt}
Assume $\Ha_{\delta}(U,V)$ is satisfied for some $\delta\in ]0,1]$. In addition, the Lyapunov functions $(e^{\delta U},e^{\delta V})$ satisfy (\ref{UV-exp2wa}) for some semi-metric $(\varrho_{\varphi},\varrho_{\psi})$.
 In this situation,  there exists some  $n_0\geq 0$ and $c_u,c_v\geq 0$ such that for any $n\geq n_0$ we have
\begin{eqnarray}
 \WW_{\varrho_\psi}(\pi_{2n},\nu_V)&\leq &c_v~\rho^{n-n_0}~
\vertiii{ \pi_{2n_0}-\nu_V}_{\psi}\nonumber\\
 \WW_{\varrho_\varphi}(\pi_{2n+1},\lambda_U)&\leq &c_u~\rho^{n-n_0}~
\vertiii{ \pi_{2n_0+1}-\lambda_U}_{\varphi},\label{cor-wnorm-est-Hd-v0}
\end{eqnarray}
In the above display, $n_0\geq 1$ and $\rho\in [0,1[$ are as in (\ref{cor-wnorm-est-Hd}) and  $(\psi,\varphi)$ are the Lyapunov functions that depend on $\delta$ defined in (\ref{def-Lyap-UV}).
\end{cor}

Assume there exists some semi-metric $\varrho_U,\varrho_V$ and some constants $c_u,c_v\in\RR$ such that 
\begin{eqnarray}
U(x_1)+U(x_2)&\geq& c_u+ \varrho_{U}(x_1,x_2),\nonumber\\
V(y_1)+V(y_2)&\geq& c_v+ \varrho_{V}(y_1,y_2).\label{dUV-r}
\end{eqnarray}
In this case, by (\ref{eq-kanto-estimate}), for any $p\geq 1$ there exists some constant $c_p>0$ such that
$$
\WW_{p,\varrho_U}(\pi_{2n+1},\lambda_U)^p\leq c_p~\vertiii{\pi_{2n+1}-\lambda_U}_{\varphi}
\quad
\mbox{\rm and}\quad
\WW_{p,\varrho_V}(\pi_{2n},\nu_V)^p\leq c_p~\vertiii{\pi_{2n}-\nu_V}_{\psi},
$$
with the $p$-Wasserstein semi-metrics $\WW_{p,\varrho_U},\WW_{p,\varrho_V}$ defined in (\ref{def-WW-varrho}). This yields the following Wasserstein estimates.

\begin{cor}
Assume $\Ha_{\delta}(U,V)$ is satisfied for some $\delta\in ]0,1]$. In addition, the target potential functions $(U,V)$ satisfy (\ref{dUV-r}) for some semi-metric $(\varrho_U,\varrho_V)$.
 In this situation, for any $p\geq 1$ there exists some $c_p\geq 0$ such that for any $n\geq n_0$ we have
\begin{eqnarray}
\WW_{p,\varrho_V}(\pi_{2n},\nu_V)^p&\leq &c_p~\rho^{n-n_0}~
\vertiii{ \pi_{2n_0}-\nu_V}_{\psi},\nonumber\\
\WW_{p,\varrho_U}(\pi_{2n+1},\lambda_U)^p&\leq &c_p~\rho^{n-n_0}~
\vertiii{  \pi_{2n_0+1}-\lambda_U}_{\varphi},\label{cor-wnorm-est-Hd-v}
\end{eqnarray}
In the above display, $n_0\geq 1$ and $\rho\in [0,1[$ are as in (\ref{cor-wnorm-est-Hd}) and  $(\psi,\varphi)$ are the Lyapunov functions that depend on $\delta$ defined in (\ref{def-Lyap-UV}).
\end{cor}

\begin{rmk}
Using Proposition~\ref{prop-uest-pi-ratio-2} in Appendix \ref{appSchrodinger} and following exactly the same argument as in the proof of Theorem~\ref{lm-tx-2}, the estimate (\ref{log-est}) remains valid when the weaker condition $\Da_{\varphi,\psi}(\Sa)$ holds for some $n_0\geq 0$, provided that
$$
~\left(\Vert \frac{d\pi_{2n_0}}{d\nu_V}\Vert_{\psi}\vee \Vert\frac{d\nu_V}{d\pi_{2n_0}}\Vert_{\psi}\right)\vee
\left(\Vert \frac{d\pi_{2n_0+1}}{d\lambda_U}\Vert_{\varphi}\vee \Vert\frac{d\lambda_U}{d\pi_{2n_0+1}}\Vert_{\varphi}\right)<\infty.
$$
\end{rmk}

\subsection{Entropy exponential decays}\label{sec-entrr}
Let us discuss some fairly straightforward consequences of the estimate (\ref{log-est}). First note that (\ref{log-est}) ensures the convergence of the series 
\begin{eqnarray}
\lim_{n\rightarrow\infty}U_{2n}=\UU&:=&U_0+\sum_{p\geq 0}\log{\frac{d \pi_{2p+1}}{d\lambda_U}}\quad \mbox{\rm in $\Ba_{\varphi}(\XX)$,}
\nonumber\\
\lim_{n\rightarrow\infty}
V_{2n}=\VV&:=&V_0+\sum_{p\geq 0}
\log{\frac{ d\pi_{2p}}{d\nu_V}}\quad \mbox{\rm in $\Ba_{\psi}(\YY)$},\label{hyp-cv}
\end{eqnarray}
where $U_0=U$ and $V_0=0$. Applying Lebesgue's dominated convergence theorem, the integral equations stated in (\ref{prop-schp}) converge as $n\rightarrow\infty$ to a
 system of integral equations
$$
\UU=U+\log{\Qa(e^{-\VV})}\quad\mbox{\rm and}\quad
\VV=V+\log{\Ra(e^{-\UU})}.
$$
More precisely, using the estimates in (\ref{log-est}) it is straightforward to prove the following corollary.
\begin{cor}
Under the assumptions of Theorem~\ref{lm-tx-2} for any $n\geq 1$ we have
\begin{equation}\label{ac-series} 
\|\UU-U_{2n}\|_{\psi}\vee \|\VV-V_{2n}\|_{\varphi}\leq \frac{c}{1-\rho^2}~\rho^{2n}
\end{equation}
with the constant $c>0$ as in (\ref{log-est}). In addition, we have the relative entropy estimates
\begin{equation}\label{ac-ent-series} 
 \mbox{\rm Ent}\left(\lambda_U~|~\pi_{2n+1}\right)\leq c~ \lambda_U(\varphi)~
\rho^{2n}~
\quad
\mbox{and}\quad
 \mbox{\rm Ent}\left(\nu_V~|~\pi_{2n}\right)\leq c~ \nu_V(\psi)~
\rho^{2n}.
\end{equation}
\end{cor}
With this notation at hand, the bridge distribution   (\ref{def-entropy-pb-v2}) has the form
\begin{equation}\label{sinhorn-entropy-form-Sch-lim}
\lim_{n\to\infty} \Pa_{n}(d(x,y)) = \PP(d(x,y)):=e^{-W(x,y)}~\lambda_{\UU}(dx)~\nu_{\VV}(dy).
\end{equation}
Moreover, we have
\begin{eqnarray*}
  \mbox{\rm Ent}(\PP~|~\Pa_{2n})  &=&\lambda_{U}(U_{2n}-\UU)+\nu_V(V_{2n}-\VV)\\
  &=&\lambda_{U}\left(\sum_{q\geq  n}\log{\frac{d\lambda_U}{d \pi_{2q+1}}}\right)+\nu_V\left(\sum_{q\geq n}
\log\frac{d\nu_V}{ d\pi_{2q}}\right).
\end{eqnarray*}
The exponential estimate (\ref{log-est}) also justifies the interchange of summation and integration. The next corollary is now a direct consequence of the relative entropy estimates in (\ref{ac-ent-series}). 
\begin{cor}\label{cor-1-i}
Under the assumptions of Theorem~\ref{lm-tx-2} we obtain
\begin{eqnarray*}
  \mbox{\rm Ent}(\PP~|~\Pa_{2n}) 
  &=&\sum_{p\geq n} \mbox{\rm Ent}\left(\lambda_U~|~\pi_{2p+1}\right)+\sum_{p\geq n} \mbox{\rm Ent}\left(\nu_V~|~\pi_{2p}\right)
  \quad \text{for any $n \ge 1$}.
\end{eqnarray*}
In addition, there exists some $c>0$ such that, for any $n\geq 1$, we have the entropy exponential decay estimates
$$
  \mbox{\rm Ent}(\PP~|~\Pa_{2n})  \leq   c~ \rho^{2n}
~ (\lambda_U(\varphi)+\nu_V(\psi)).
$$
\end{cor}

%%%%%%%%%%%%%%%%%%%%%%%%%%%%%
% Some illustrations		%
%%%%%%%%%%%%%%%%%%%%%%%%%%%%%
\section{Some illustrations}\label{sec-illustrations}

In this section, we illustrate the scope and flexibility of the theoretical results developed above through a collection of representative examples. We show how the Lyapunov, drift-minorization, and weighted contraction conditions can be verified in a variety of settings, ranging from polynomial and heavy-tailed models to boundary-driven, semi-compact, and Gaussian frameworks. These examples highlight how the abstract assumptions translate into explicit growth, integrability, and regularity conditions on the target potentials and transition kernels.

\subsection{Polynomial growth models}\label{sec-poly}
Consider some finite dimensional normed vector space $\XX=\YY$ equipped with measure $\lambda=\nu$ such that for any $p>0$ and $\tau>0$ we have
\begin{equation}\label{elambda}
\int e^{-\tau\Vert x\Vert^{p}}~\lambda(dx)<\infty.
\end{equation}
 In this context, we have $\lambda_U(dx)=e^{-U(x)}~\lambda(dx)$ and
$\nu_V(dx)=\lambda_V(dx)=e^{-V(x)}~\lambda(dx)$. 

Assume that the functions $U,V$ and $W_{U},W^{V}$ are continuous.
In addition, there exists some parameters $c_w,c_u,c_v\in\RR$ as well as $\varsigma_{u},\varsigma_{v},\tau_u,\tau_v,p_u,p_v> 0$ and $q_u,q_v\geq 0$ such that
 \begin{eqnarray}
W(x,y)&\leq& c_w+\varsigma_{u}~\Vert x\Vert^{p_u}+
\varsigma_{v}~\Vert y\Vert^{p_v},\nonumber\\
  U(x)&\geq& c_u+\tau_u~\Vert x\Vert^{p_u+q_u}\quad \mbox{\rm and}\quad
V(x)\geq c_v+ \tau_v~ \Vert x\Vert^{p_v+q_v}.\label{pq-u}
\end{eqnarray}

\begin{prop}\label{elambda-poly-prop}
Condition $\Ha^{\prime}_{\delta}(U,V)$ is satisfied for some $0<\delta<1/2$ whenever, for any parameter $\iota\in \{u,v\}$, one of the following assertions is satisfied
\begin{equation}\label{elambda-poly}
\left(q_\iota>0\right)\quad \mbox{or}\quad \left(
q_\iota=0\quad \mbox{and}\quad \tau_\iota/2>\varsigma_\iota\right).
\end{equation}
\end{prop}

\proof
The exponential moment condition (\ref{elambda}) ensures that
(\ref{huv-2}) holds for any $0<\delta<1/2$ and the probability measures $\lambda_U$ and $\lambda_V$ have finite $p_u$ and $p_v$ absolute moments
$$
m^{p_u}_u:=\int \Vert x\Vert^{p_u}\lambda_U(dx)<\infty \quad \mbox{\rm and}\quad
m^{p_v}_v:=\int \Vert x\Vert^{p_u}\lambda_V(dx)<\infty.
$$
Also note that
 \begin{eqnarray*}
\Ua_{\delta}(x)=\delta U(x)-W^V(x)&\geq& (\delta c_u-c_w-\varsigma_{v}~m^{p_v})+\delta \tau_u~\Vert x\Vert^{p_u+q_u}-\varsigma_{u}~\Vert x\Vert^{p_u}.
\end{eqnarray*}
This shows that $\Ua_{\delta}$ has compact sub-level sets with $0<\delta<1/2$ as soon as
$$
q_u>0\quad \mbox{\rm or}\quad \left(
q_u=0\quad \mbox{\rm and}\quad\delta \tau_u>\varsigma_u\right).
$$
In the same vein, condition (\ref{elambda-poly}) with $a=v$ ensures
$\Va_{\delta}$ has compact sub-level sets.
\cqfd

For instance, for bi-Laplace transitions on $\XX=\RR=\YY$ we have
$$
\begin{array}{l}
\displaystyle\Qa(x,dy)=\frac{\varsigma}{2}~e^{-\varsigma |x-y|}~ dy\quad \mbox{\rm with $\varsigma>0$}\\
\\
\Longrightarrow
W(x,y)=c_w+\varsigma |x-y|\leq c_w+
\varsigma |x|+\varsigma |y|\quad \mbox{\rm with $c_w:=\log{(2/\varsigma)}$ }\\
\\
\Longrightarrow |W(x,z_1)-W(x,z_2)|\vee |W(z_1,y)-W(z_2,y)|\leq 
\varsigma~ |z_1-z_2|.
\end{array}$$

For linear Gaussian models, the r.h.s. condition in (\ref{elambda-poly}) is satisfied as soon as the covariance of the  linear Gaussian transition $\Qa$ is defined by (\ref{def-Qa}) is sufficiently large. In the same context of linear and Gaussian models, we show in Proposition~\ref{prop-gauss-model}  that this condition is in fact equivalent to condition  $\Ha_{\delta}(U,V)$.

By (\ref{pq-u}) there exists some constant $\overline{\tau}_u>0$ such that
\begin{eqnarray*}
U(x_1)+U(x_2)&\geq&  \overline{c}_u+
\overline{\tau}_u~\left(\Vert x_1\Vert+
\Vert x_2\Vert\right)^{p_u+q_u}\quad
\mbox{\rm with}\quad \overline{c}_u:=2c_u.
\end{eqnarray*}
This yields the estimate
$$
U(x_1)+U(x_2)\geq \overline{c}_u+\varrho_{U}(x_1,x_2)\quad
\mbox{\rm with}\quad
\varrho_{U}(x_1,x_2):=\overline{\tau}_u~\Vert x_1-x_2\Vert^{p_u+q_u}.
$$
Similarly, there exist some constants $ \overline{c}_u\in\RR$ and $\overline{\tau}_v>0$ such that
$$
V(x_1)+V(x_2)\geq   \overline{c}_v+
\varrho_{V}(x_1,x_2)\quad
\mbox{\rm with}\quad
\varrho_{V}(x_1,x_2):=\overline{\tau}_v~\Vert x_1-x_2\Vert^{p_v+q_v}.
$$
We conclude that $(U,V)$ satisfy (\ref{dUV-r}) with the semi-metric $(\varrho_U,\varrho_V)$ so that the 
Wasserstein estimates (\ref{cor-wnorm-est-Hd-v}) hold.

\subsection{Heavy tailed transitions}\label{heavy-tail-sec}
Consider the normed vector space $\XX=\YY$ and the reference measure $\lambda=\nu$  discussed in Section~\ref{sec-poly}. Assume some parameter  $c_w\in \RR$, $\tau_{p,w},\tau_{q,w}>0$  and $p_w\wedge q_w>1$ such that
$$
W(x,y)\leq 
c_w+\log\left(1+\tau_{p,w}~\Vert x\Vert^{p_w}+
\tau_{q,w}~\Vert y\Vert^{q_w}\right).
$$
For instance, for Cauchy transitions on $\XX=\RR=\YY$ we have
$$
\begin{array}{l}
\displaystyle\Qa(x,dy)=\frac{1}{\pi a}~\frac{1}{1+((x-y)/a)^2}~dy\quad \mbox{\rm with}\quad a>0\\
\\
\displaystyle\Longrightarrow W(x,y)=\log{(a\pi)}+\log{\left(1+((x-y)/a)^2\right)}\quad \mbox{\rm and}\quad\vert\partial_xW(x,y)\vert\vee \vert\partial_yW(x,y)\vert\leq 1/a\\
\\
\displaystyle\Longrightarrow c_w=\log{(a\pi)},\qquad \tau_{p,w}=2/a^2=\tau_{q,w}\quad \mbox{\rm and}\quad
p_w=q_w=2.
\end{array}$$

\subsubsection{Exponential targets} 

Consider some continuous target potential functions $U,V$ such that
$$
  U(x)\geq c_u+\tau_u~\Vert x\Vert^{p_u}\quad \mbox{\rm and}\quad
    V(x)\geq c_v+\tau_v~\Vert x\Vert^{q_v}\quad \mbox{\rm with $c_u,c_v\in\RR$ and     $p_u,q_v>0$}.
$$
Also assume that  $(W_{U},W^{V})$ are continuous functions.
Following the proof of Proposition~\ref{elambda-poly-prop}, the exponential moment condition (\ref{elambda}) ensures that
(\ref{huv-2}) holds for any $0<\delta<1/2$ and the probability measures $\lambda_U$ and $\nu_V=\lambda_V$ have finite $p_w$ and $q_w$ absolute moments.
By Jensen's inequality, we also have
 \begin{eqnarray*}
\Ua_{\delta}(x)&\geq&
\delta c_u+\delta\tau_u~\Vert x\Vert^{p_u}-c_w-\log\left(1+\tau_{p,w}~\Vert x\Vert^{p_w}+
\tau_{q,w}~m_v^{q_w}\right)\stackrel{\Vert x\Vert\rightarrow\infty}{\longrightarrow}\infty
\end{eqnarray*}
In the same vein we check that
\begin{equation}\label{infty-UV-delta}
\lim_{\Vert x\Vert\rightarrow\infty}(\delta U(x)-W^V(x))=\infty=\lim_{\Vert x\Vert\rightarrow\infty}(\delta V(x)-W_U(x)).
\end{equation}
We conclude that $\Ha^{\prime}_{\delta}(U,V)$ is met for any $0<\delta<1/2$.
Also note that $(U,V)$ satisfy (\ref{dUV-r}) with the semi-metrics 
$$
\varrho_{U}(x_1,x_2):=\overline{\tau}_u~\Vert x_1-x_2\Vert^{p_u}
\quad
\mbox{\rm and}\quad
\varrho_{V}(x_1,x_2):=\overline{\tau}_v~\Vert x_1-x_2\Vert^{q_v}
$$
for some parameters $\overline{\tau}_u,\overline{\tau}_v>0$ so that the 
Wasserstein estimates (\ref{cor-wnorm-est-Hd-v}) hold.

\subsubsection{Heavy tailed targets} 

Consider continuous target potential functions $U,V$ such that
$$
  U(x)\geq c_u+\log{\left(1+\tau_u~\Vert x\Vert^{p_u}\right)}\quad \mbox{\rm and}\quad
    V(x)\geq c_v+\log{\left(1+\tau_v~\Vert x\Vert^{q_v}\right)}
$$
with $c_u,c_v\in\RR$  and some parameters 
$$
p_u>1+2p_w\quad \mbox{\rm and}\quad
q_v>1+2q_w.
$$
Condition (\ref{huv-2}) is met as soon as
$$
(1-2\delta)p_u>1\quad \mbox{\rm and}\quad
(1-2\delta)q_v>1
$$
On the other hand, the probability measures $\lambda_U$ and $\lambda_V$ have finite $p_w$ and $q_w$ absolute moments as soon as
$$
p_u-p_w>1\quad \mbox{\rm and}\quad
q_v-q_w>1.
$$
Finally, (\ref{infty-UV-delta}) holds when
$$
\delta p_u-p_w>0\quad \mbox{\rm and}\quad
\delta q_v-q_w>0.
$$
This shows that $\Ha_{\delta}(U,V)$ for any $0<\delta<1/2$ such that
$$
\frac{p_w}{p_u}<\delta<\frac{1}{2}~\left(1-\frac{1}{p_u}\right)
\quad \mbox{\rm and}\quad
\frac{q_w}{q_v}<\delta<\frac{1}{2}~\left(1-\frac{1}{q_v}\right).
$$

Also note that the exponential functions $(e^{\delta U},e^{\delta V})$ satisfy (\ref{UV-exp2wa}) with the semi-metrics 
$$
\varrho_{\varphi}(x_1,x_2):=\overline{\tau}_u~\Vert x_1-x_2\Vert^{p_u}
\quad
\mbox{\rm and}\quad
\varrho_{\psi}(x_1,x_2):=\overline{\tau}_v~\Vert x_1-x_2\Vert^{q_v}
$$
for some parameters $\overline{\tau}_u,\overline{\tau}_v>0$ so that the 
Wasserstein estimates (\ref{cor-wnorm-est-Hd-v0}) hold.

\subsection{Some boundary states}\label{bst-mod-sec}

Consider a bounded Lipschitz domain $\XX=\YY\subset\RR^d$  equipped with the uniform measure $\lambda=\nu$.
 Denote by $\varrho_{\partial}(x,\partial\XX)$ the distance from $x\in\XX$ to the boundary $\partial\XX$. We assume that the target distribution $\lambda_U$ has a continuous potential function $U$ and satisfies the estimates
$$
\lambda_U(dx)\leq c_{q,u}~\varrho_{\partial}(x,\partial\XX)^{q_u}~\lambda(dx)
$$
for some parameters $c_{p,u},c_{q,u}>0$ and $q_u>0$. Similarly, 
$\nu_V=\lambda_V$ has a continuous potential function $V$ and satisfies the above estimates for some $c_{p,v},c_{q,v}>0$ and $q_v>0$. Finally, we assume that the integrated functions $(W^V,W_U)$ are continuous.

When $W(x,y)$ is bounded, Sinkhorn transitions are strongly mixing, in the sense that  $\cchi(\Sa)<1$ with the uniform contraction parameter $\cchi(\Sa)$ defined in (\ref{def-cchi-S}). The exponential decays of $\Phi$-entropies associated with these strong mixing models are discussed in Section~\ref{strong-mixing-sec}. 
Note that for any $x\in\XX$ we have
$$
U(x)\geq -\log{c_{q,u}}+q_u~\log{\frac{1}{\varrho_{\partial}(x,\partial\XX)}}
\Longleftrightarrow
e^{\delta U(x)}\geq \frac{1}{c_{q,u}^{\delta}}~\frac{1}{\varrho_{\partial}(x,\partial\XX)^{\delta q_u}}.
$$
By Remark~\ref{rmk-W-bb} the above models satisfy
$\Ha_{\delta}(U,V)$ for any $\delta\in ]0,1[$. In addition,
 following the proof of (\ref{ref-poly-lyap}) we check that the exponential function $e^{\delta U}$ satisfies (\ref{UV-exp2wa}) with the semi-metric 
 $$
  \varrho_{\varphi}(x_1,x_2)=\overline{c}_u~\left(\frac{1}{\varrho_{\partial}(x_1,\partial\XX)}-\frac{1}{\varrho_{\partial}(x_2,\partial\XX)}\right)^{\delta q_u}~1_{x_1\not=x_2}
 \quad \mbox{\rm
  for some parameter $\overline{c}_u>0$.}
  $$
  Similarly, the exponential function $e^{\delta V}$ satisfies (\ref{UV-exp2wa}) with the  metric $\varrho_{\psi}$ defined as  $ \varrho_{\varphi}$
 by replacing $(c_u,q_u)$ by $(c_v,q_v)$, so that the 
Wasserstein estimates (\ref{cor-wnorm-est-Hd-v0}) hold.

When $W(x,y)=W(y,x)$ is symmetric and non necessarily bounded, assume we have an estimate of the form
$$
W(x,y)\leq c_w+q_w~\log{\frac{1}{\varrho_{\partial}(x,\partial\XX)}}
$$
for some $c_w\in\RR$ and some parameters $0\leq q_w<(q_u\wedge q_v)/2$. In this case
condition $\Ha_{\delta}(U,V)$ holds for any $\frac{q_w}{q_u\wedge q_v}<\delta<1/2$.

\subsection{Beta marginal models}\label{beta-w-sec}

 Next we illustrate the boundary state models discussed in Section~\ref{bst-mod-sec}
with the open unit interval $
\XX=\YY=]0,1[
$  with boundary $\partial \XX=\{0,1\}$. The distance to the boundary is given explicitly by the function
$$
x\in \XX\mapsto \varrho_{\partial}(x,\partial\XX)=x~1_{]0,1/2[}+(1-x)~1_{[1/2,1[}(x)\quad\mbox{\rm and}\quad
\lambda(dx)=1_{]0,1[}(x)dx.
$$
Note that
$$
\frac{1}{2}~\varrho_{\partial}(x,\partial\XX)\leq
x(1-x)\leq \varrho_{\partial}(x,\partial\XX)
$$
and consider the reference Markov transition $\Qa$ defined by
the semi-circle distribution (\ref{semicircle-ex}).
Let $(\lambda_U,\lambda_V)$ be some Beta marginal distributions given by
$$
\lambda_U(dx)=\frac{1}{c_u}~x^{a_u}(1-x)^{b_u}~
\lambda(dx)\quad\mbox{\rm and}\quad
\lambda_V(dx)=\frac{1}{c_v}~x^{a_v}(1-x)^{b_v}~
\lambda(dx)
$$
for some normalizing constants $c_{u},c_v>0$ and some parameters
$a_u,a_v,b_u,b_v>1$. In this context, we have
$$
 \frac{1}{2^{a_u\vee b_u}c_u}~\varrho_{\partial}(x,\partial\XX)^{a_u\vee b_u}~\lambda(dx)\leq
\lambda_U(dx)\leq \frac{1}{c_u}~\varrho_{\partial}(x,\partial\XX)^{a_u\wedge b_u}~\lambda(dx).
$$
Also observe that
\begin{eqnarray*}
1\leq \frac{1}{1-(x-y)^2}&=&\frac{1}{(1-(x-y))(1+(x-y))}\\
&=&\frac{1}{y+(1-x)}~\frac{1}{x+(1-y)}\leq \frac{1}{x(1-x)}\wedge
\frac{1}{y(1-y)},
\end{eqnarray*}
which yields the estimates
\begin{eqnarray*}
W(x,y)&=&\log{q(x)}+\frac{1}{2}~\log{\left(\frac{1}{(1-(x-y))(1+(x-y))}\right)}\\
&\leq&  \frac{1}{2}~\log\left(\frac{1}{x(1-x)}\right)\quad \mbox{\rm with}\quad \frac{1}{3}\leq q(x):=\int_0^1 \sqrt{1-(x-y)^2}dy\leq 1.
\end{eqnarray*}
Note that for any $(x,y)\in ]0,1[^2$ we obtain
$$
\vert\partial_x\sqrt{1-(x-y)^2}\vert=\frac{|y-x|}{\sqrt{((1-x)+y)((1-y)+x)}} \leq g(y):=\frac{1}{\sqrt{y(1-y)}}\Longrightarrow \lambda(g)<\infty
$$
as well as
$$
\vert \partial_x\log{(1-(x-y)^2)}\vert\leq g(x)^2\wedge g(y)^2.
$$
In this example, the functions $(U,V)$ and $(W_{U},W^{V})$ are differentiable functions.
Condition (\ref{huv-prim}) is clearly met for any $0<\delta<1/2$ and we have
$$
\delta U(x)-W^{V}(x)\geq \delta\log{(c_u)}+(\delta a_u-1/2)\log{\frac{1}{x}}+(\delta b_u-1/2)\log{\frac{1}{1-x}}.
$$
This shows that $\Ha^{\prime}_{\delta}(U,V)$ is met for any 
$
((a_u\wedge b_u)\wedge (a_v\wedge b_v))^{-1}<2\delta<1 
$. We also have
 $$
U(x_1)+U(x_2)\geq \overline{c}_u:+ \varrho_U(x_1,x_2)
 $$
 with  $\overline{c}_u:=2\log{(c_u)}$ and the metric
 $$
 \varrho_U(x_1,x_2):=a_u~\left\vert \log{\frac{1}{x_1}}-\log{\frac{1}{x_2}}\right\vert+b_u
\left\vert\log{\frac{1}{1-x_1}}-\log{\frac{1}{1-x_2}}\right\vert.
 $$
 
 We conclude that $U$ satisfies  (\ref{dUV-r}).
 Similarly $V$ satisfies  (\ref{dUV-r}) with a metric $\varrho_V$ defined as $\varrho_U$ by replacing $(a_u,b_u,c_u)$ by  $(a_v,b_v,c_v)$. 
 
Recall that $uv\geq (u+v)/2$ and $u^pv^q\geq (uv)^{p\wedge q}$ for any $u,v\geq 1$ and any parameters $p,q>0$ we see that
 \begin{eqnarray*}
 e^{\delta U(x)}&\geq& \frac{c_u^{\delta}}{2}~\left(
 \left(\frac{1}{x}\right)^{\delta a_u}+\left(\frac{1}{1-x}\right)^{\delta b_u}\right),\\
  e^{\delta U(x)}&\geq& c_u^{\delta}~\left(\frac{1}{x(1-x)}\right)^{\delta (a_u\wedge b_u)}\geq
\overline{c}_u^{\delta}~\frac{1}{\varrho_{\partial}(x,\partial\XX)^{\delta (a_u\wedge b_u)}}\quad\mbox{\rm with}\quad \overline{c}_u:=  2^{(a_u\wedge b_u)}c_u.
 \end{eqnarray*}
 Following the proof of (\ref{ref-poly-lyap}) we check that the exponential function $e^{\delta U}$ satisfies (\ref{UV-exp2wa}) with 
 $$
  \varrho_{\varphi}(x_1,x_2)=\overline{c}_u~\left(\frac{1}{\varrho_{\partial}(x_1,\partial\XX)}-\frac{1}{\varrho_{\partial}(x_2,\partial\XX)}\right)^{\delta (a_u\wedge b_u)}
 $$
 as well as with the semi-metric
 $$
 \varrho_{\varphi}(x_1,x_2)=\overline{c}_u~\left(\left\vert \frac{1}{x_1}-\frac{1}{x_2}\right\vert^{\delta a_u}+\left\vert \frac{1}{1-x_1}-\frac{1}{1-x_2}\right\vert^{\delta b_u}\right)
 $$
 for some parameter $\overline{c}_u>0$. Similarly, the exponential function $e^{\delta V}$ satisfies (\ref{UV-exp2wa}) with the complete metric $\varrho_{\psi}$ defined as  $ \varrho_{\varphi}$
 by replacing $(a_u,b_u,c_u)$ by $(a_v,b_v,c_v)$, so that the 
Wasserstein estimates (\ref{cor-wnorm-est-Hd-v0}) hold.

\subsection{Weibull marginal models}\label{w-sec}

We equip the open half-line $\XX=\YY=]0,\infty[$ with the Lebesgue measure $\lambda(dx)=\nu(dx)=1_{]0,\infty[}(x)dx$ and
consider the Weibull marginal measure
$$
\displaystyle\lambda_U(dx)=\frac{a_u+1}{b_u}~\left(\frac{x}{b_u}\right)^{a_u}
e^{-(x/b_u)^{a_u+1}}~\lambda(dx)
$$
with potential
$$
U(x)=c_u+(x/b_u)^{a_u+1}+a_u~\log{\frac{1}{x}},\quad
\mbox{\rm where}\quad c_u:=\log{\frac{b_u}{a_u+1}}+a_u\log{b_u}
$$
and $a_u,b_u>0$. Note that (\ref{huv-2}) holds for any $0<\delta<1/2$ and the probability measures $\lambda_U$ and $\lambda_V$ have finite $p_u$ and $p_v$ absolute moments
$m^{p_u}_u$ and $m^{p_v}_v$ for any $p_u,p_v$.
We define similarly $\lambda_V$ with some parameters $a_v,b_v>0$.

Consider a continuous function $W(x,y)$ satisfying (\ref{pq-u}) for some 
parameters $c_w\in\RR$ as well as $p_u,p_v> 0$. In this context, we have
$$
W^{V}(x)\leq  c_w+\varsigma_{u}~x^{p_u}+
\varsigma_{v}~m_v^{p_v}\quad \mbox{\rm with}\quad
m^{p_v}_v:=\int_0^{\infty}~x^{p_u}\lambda_V(dx)<\infty.
$$
This implies that
$$
\delta U(x)-W^{V}(x)\geq
\overline{c}_u+\delta a_u~\log{\frac{1}{x}}+\frac{\delta}{b_u^{a_u+1}}~x^{a_u+1}
-\varsigma_{u}~x^{p_u}
\quad\mbox{\rm
with}\quad
\overline{c}_u:=\delta c_u-c_w-
\varsigma_{v}~m_v^{p_v}
$$
Assume that $(W_{U},W^{V})$ are both continuous functions. In this situation,
  $\Ua_{\delta}$ has compact sub-level sets with $0<\delta<1/2$ as soon as
$$
a_u+1>p_u\quad \mbox{\rm or}\quad\left(
p_u=a_u+1\quad\mbox{\rm and}\quad \varsigma_{u}~b_u^{a_u+1}<\delta<1/2\right).
$$
The r.h.s. condition requires that $\varsigma_{u}<b_u^{-(a_u+1)}/2$. The analysis of these Weibull models now follows the same line of arguments as the ones we used for the Beta models. %, thus it is left to the reader.

\subsection{Exponential marginal models}\label{expo-sec}

We equip the half-line $\XX=\YY=[0,\infty[$ with the Lebesgue measure $\lambda(dx)=1_{[0,\infty[}(x)dx$ and
consider the exponential marginal measures
$$
\displaystyle\lambda_U(dx)=\tau_u~e^{-\tau_u x}~\lambda(dx)\quad \mbox{\rm and}\quad 
\lambda_V(dx)=\tau_v~e^{-\tau_v x}~\lambda(dx)\quad\mbox{\rm with}\quad\tau_u\wedge \tau_v>0\\
$$
for which we have the potentials
$$ 
U(x)=c_u+\tau_u x\quad
\mbox{\rm and}\quad\quad V(x)=c_v+\tau_v x\quad
\mbox{\rm with}\quad (c_u,c_v):=(-\log{\tau_u},-\log{\tau_u}).
$$
Note that the functions $x\in [0,\infty[\mapsto \tau_{\iota} x\in [0,\infty[$
with $\iota\in\{u,v\}$  have compact sub-level sets.
Let $\Qa(x,dy)$ be the exponential transition defined for any $x\in [0,\infty[$ by
$$
\displaystyle\Qa(x,dy)=\frac{\varsigma}{1+(1-e^{-\lambda x})}~e^{-\varsigma |x-y|}~ \lambda(dy)\quad \mbox{\rm with $0<\varsigma<(\tau_u\wedge\tau_v)/2$}\\
$$
which yields the cost function
$$
W(x,y)\leq c_w+\varsigma |x-y|\leq c_w+
\varsigma ~x+\varsigma~ y,\quad \mbox{\rm with $c_w:=\log{(2/\varsigma)}$ }.
$$
In this context, condition $\Ha^{\prime}_{\delta}(U,V)$ holds for any 
$
({\varsigma}/{\tau_u})\vee
({\varsigma}/{\tau_v}) <\delta<{1}/{2}
$.
\subsection{Semi-compact models}\label{sec-semi-discrete}

  In many practical situations of interest, the target state spaces are different.  
  Next we illustrate this situation with a semi-discrete model and a semi-compact models.

 \subsubsection*{Semi-discrete model.} 

Consider a finite or countable subset $\XX\subset\RR^d$
equipped with the counting measure $\lambda$ and the discrete topology.  
In this context, the potential function $U$ is a function
$$
 U~:~x\in \XX\mapsto U(x) \in[0,\infty[\quad \mbox{\rm such that}\quad \sum_{x\in\XX}e^{-U(x)}=1.
$$
Let $\nu_V(dy)=e^{-V(y)}~dy$ be a Boltzmann-Gibbs probability 
measure on $\RR^d$ and assume that $(U,V,W)$ are continuous functions satisfying the polynomial 
growth conditions (\ref{pq-u}). In this context, the sub-level sets of $U$ are finite, thus compact in $\XX$.

Following the proof of
Proposition~\ref{elambda-poly-prop} these conditions ensure 
that  $\Ha^{\prime}_{\delta}(U,V)$ holds for some $0<\delta<1/2$ as soon as condition (\ref{elambda-poly}) is satisfied. The analysis of heavy tailed transitions follows word-for-word the same lines of arguments as in
Section~\ref{heavy-tail-sec}, this it is skipped.

 \subsubsection*{Semi-compact model.} 

Consider the uniform measure $\lambda_U$ on a compact space $\XX\subset \YY=\RR^d$ associated with a null potential $U=0$, which yields $\lambda_U=\lambda$.
 Let $\nu_V$ be a centered Gaussian distribution with unit variance  on $\YY$ and assume that the transition potential is chosen so that
$$
W(x,y)\leq w(x)+a~\Vert y\Vert^{p}\Longrightarrow
W_U(y)\leq \lambda(w)+b~\Vert y\Vert^{p}
$$
for some continuous function $w(x)$ and some parameters $a\leq 2$ and $p\leq 1$.
These models satisfy  $\Ha^{\prime}_{\delta}(U,V)$ as soon as
\begin{equation}\label{re-pdl}
(p=2~~\mbox{\rm and}~~b<\delta/2)\quad \mbox{\rm or}\quad
(p<2)
\end{equation}
We check this claim using the fact that
$$
(\ref{re-pdl})\Longrightarrow
\delta V(y)-W_U(y)\geq c+\frac{\delta}{2}~\Vert y\Vert^2-b~\Vert y\Vert^{p}
\stackrel{\Vert y\Vert\rightarrow\infty}{\longrightarrow} +\infty
\quad \mbox{\rm for some $c\in\RR$. }
$$

\subsection{Linear-Gaussian models}\label{sec-lin-gauss}
Consider the state spaces $\XX=\YY=\RR^d$, for some integer $d\geq 1$ and the Lebesgue measure $\lambda(dz)=\nu(dz)=dz$. Hereafter, points $x$ in the Euclidean space $\XX=\RR^d$ are represented by $d$-dimensional column vectors (or, equivalently, by $d \times 1$ matrices). 
 
Let $\Ga l_d$ denote the general linear group of $(d\times d)$-invertible matrices, and $\Sa^+_d\subset \Ga l_d$ the subset 
 of positive definite matrices. We sometimes use the L\" owner partial ordering notation $v_1\geq v_2$ to mean that a symmetric matrix $v_1-v_2$ is positive semi-definite (equivalently, $v_2 - v_1$ is negative semi-definite), and $v_1>v_2$ when $v_1-v_2$ is positive definite (equivalently, $v_2 - v_1$ is negative definite). Given $v\in\Sa_d^+$ we denote by $v^{1/2}$  the principal (unique) symmetric square root.

Denote by $\nu_{m,\sigma}$ the Gaussian distribution on $\RR^d$ with mean $m\in\RR^d$ and covariance matrix $\sigma \in \Sa^+_d$. In addition, let $g_{\sigma}$ denote the probability density function of the distribution $\nu_{0,\sigma}$, with covariance matrix $\sigma \in \Sa^+_d$.  
 Consider the target probability measures $(\lambda_U,\nu_V)=(\nu_{m,\sigma},\nu_{\overline{m},\overline{\sigma}})$ for some given parameters
$$
(m,\overline{m})\in(\RR^d\times\RR^{d})\quad \text{and}\quad 
 (\sigma,\overline{\sigma})\in (\Sa^+_d\times \Sa^+_{d}).
$$
Equivalently, the log-density functions $(U,V)$ are defined by 
\begin{equation}\label{def-U-V}
U(x)=-\log{g_{\sigma}(x-m)}
\quad \mbox{\rm and}\quad V(y):=-\log{g_{\overline{\sigma}}(y-\overline{m})}.
\end{equation}
The linear Gaussian transition $\Qa$ is defined by (\ref{def-Qa}), with the potential
\begin{equation}\label{def-W}
W(x,y)=-\log{g_{\tau}(y-(\alpha+\beta x))}\quad
\end{equation}
for some given $(\alpha,\beta,\tau)\in  \left(\RR^{d}\times\Ga l_d\times \Sa^+_{d}\right)$.
In this context, we have $\jmath(\Qa)=0$ and Hilbert projective techniques cannot be applied. However, a refined exponential stability analysis of the Sinkhorn algorithm for Gaussian models has been developed in~\cite{adm-24}. For instance, following Remark 5.4 in~\cite{adm-24}, there exists some $n_0\geq 1$ such that for any $n\geq n_0$ we have the estimate
\begin{equation}\label{def-rho-t}
\Vert\pi_{2n}-\nu_{\overline{m},\overline{\sigma}}\Vert_{\sf tv}\leq c~\rho^{n}\quad \mbox{\rm with}\quad\rho=(1+\lambda_{\sf min}(\Delta_{\varpi}))^{-1}<1.
\end{equation}
The matrix $\Delta_{\varpi}$ in expression \eqref{def-rho-t} above is defined as
\begin{equation}\label{def-fix-ricc-1}
\Delta_{\varpi}:=\frac{\varpi}{2}+\left(\varpi+\left(\frac{\varpi}{2}\right)^2\right)^{1/2}
\geq \varpi+(I+\varpi^{-1})^{-1}
\end{equation}
with
$$ 
 \varpi:=(\gamma\,\gamma^{\prime})^{-1}\quad\mbox{\rm and}\quad  \gamma:=\overline{\sigma}^{1/2}~\tau^{-1}\beta~\sigma^{1/2}.
$$
The drift-minorization techniques developed in this paper also enable us to obtain exponential convergence on weighted Banach spaces.  
\begin{prop}\label{prop1-lyap-gauss}
Condition $\Da_{\varphi,\psi}(\Sa^{\circ})$ is satisfied with the Lyapunov functions
\begin{equation}\label{first-set-glyap}
\varphi(x):=1+(x-m)^{\prime}\sigma^{-1}(x-m)\quad\mbox{\rm and}\quad
\psi(y):=1+(y-\overline{m})^{\prime}\overline{\sigma}^{-1}(y-\overline{m}).
\end{equation}
\end{prop}
A proof is provided in the Appendix \ref{GS-sec} (on page~\pageref{lyap-sec-Gauss}).

The choice of the Lyapunov functions is not unique. For instance, Condition $\Da_{\varphi,\psi}(\Sa^{\circ})$ is also met for any $\delta>0$ with the exponential Lyapunov functions
\begin{eqnarray}
\varphi(x)&:=&\exp{\left(\delta\,((x-m)^{\prime}\sigma^{-1}(x-m))^{1/2}\right)} \quad \text{and}\nonumber\\
\psi(y)&:=&\exp{\left(\delta\,((y-\overline{m})^{\prime}\,\overline{\sigma}^{-1}(y-\overline{m}))^{1/2}\right)}.
\label{second-set-glyap}
\end{eqnarray}
There also exists some $\delta_0>0$ such that $\Da_{\varphi,\psi}(\Sa^{\circ})$ is met for any $0<\delta<\delta_0$  with the exponential Lyapunov functions
\begin{eqnarray}
\varphi(x)&:=&\exp{\left(\delta\,(x-m)^{\prime}\sigma^{-1}(x-m)\right)} \quad \text{and}\nonumber\\
\psi(y)&:=&\exp{\left(\delta\,(y-\overline{m})^{\prime}\overline{\sigma}^{-1}(y-\overline{m})\right)}.\label{ref-gauss-lyap}
\end{eqnarray}

A detailed proof of the above assertions is provided in Appendix \ref{GS-sec} (on page~\pageref{lyap-sec-Gauss}).

We end this section with an equivalent formulation of condition $\Ha_{\delta}(U,V)$ in terms of the parameters of the Gaussian potentials (\ref{def-U-V}) and (\ref{def-W}).

\begin{prop}\label{prop-gauss-model}
Condition  $\Ha_{\delta}(U,V)$ is satisfied for the Lyapunov functions in (\ref{ref-gauss-lyap}) if, and only if,
\begin{equation}\label{cc-prop}
\beta^{\prime}\tau^{-1}\beta<(\delta\wedge (1-\delta))~\sigma^{-1}\quad \text{and}\quad
\tau^{-1}<(\delta\wedge (1-\delta))~\overline{\sigma}^{-1}.
\end{equation}
\end{prop}
A proof is provided in Appendix \ref{tech-proof-ap} (on page~\pageref{prop-gauss-model-proof}).

%%%
%
%%%
\subsection{Finite mixture models}\label{finite-mix-sec}

 Consider some  (fully supported) discrete probability measures  $\xi_0$ and $\xi_1$ on some some finite states $E_0$ and $E_1$ and some collection of probability measures
 $(\lambda_{U_i})_{i\in E_0}$ and  $(\nu_{V_j})_{j\in E_1}$
 for some continuous potential functions $U_i:x\in\XX\mapsto U(x)\in \RR$ and $V_j:y\in\XX\mapsto V(y)\in \RR$ indexed by $i\in E_0$ and $j\in E_1$.
The mixture distributions $\lambda_U$ and $\nu_V$ associated with these potential functions take the form
\begin{equation}\label{mix-intro}
U=-\log{\sum_{i\in E_0} \xi_0(i)~e^{-U_i}}\quad\mbox{\rm and}\quad
V=-\log{\sum_{j\in E_1} \xi_1(j)~e^{-V_j}}.
\end{equation}

\begin{prop}\label{prop-mixture-intro}
Condition $\Ha_{\delta}(U,V)$ is satisfied for some $\delta\in ]0,1[$ if conditions $\Ha_{\delta}(U_i,V_j)$ are met for the same $\delta\in ]0,1[$ and almost every $(i,j)\in (E_0\times E_1)$.
%we have
%$$
%\delta U_i-W^{V_j}\in \Ba_{\infty}(\XX)\quad\mbox{and}\quad
%\delta V_i-W_{U_j}\in \Ba_{\infty}(\YY)
%$$
%In addition, the functions $W_{U_j},W^{V_j}$ are locally bounded and we have
%$$
%\lambda_{(1-\delta)U_i}\left(e^{W^{V_j}}\right)\vee
%~\nu_{(1-\delta)V_j}\left(e^{W_{U_i}}\right)
%<\infty$$

\end{prop}
A proof is provided in Appendix \ref{tech-proof-ap} (on page~\pageref{prop-mixture-intro-proof}).

Let us now consider the mixture model (\ref{mix-intro}) with  $\XX=\YY=\RR^d$, the potential $W(x,y)$ defined in Eq. (\ref{def-W}) and the log-densities
\begin{equation}\label{gauss-mix}
-U_i(x)=\log{g_{\sigma_i}(x-m_i)}
\quad \mbox{\rm and}\quad -V_j(y):=\log{g_{\overline{\sigma}_j}(y-\overline{m}_j)},
\end{equation}
where $(m_i,\overline{m}_j)_{(i,j)\in (E_0\times E_0)}\in(\RR^d\times\RR^{d})^{E_0\times E_1}$ and 
 $(\sigma_i,\overline{\sigma}_j)_{(i,j)\in (E_0\times E_0)}\in(\Sa^+_d\times \Sa^+_{d})^{E_0\times E_1}$ are given parameters. Combining Proposition~\ref{prop-mixture-intro} with Proposition~\ref{prop-gauss-model} we readily obtain the following result.

\begin{prop}
The Gaussian mixture model defined in (\ref{mix-intro}) with the Gaussian log-densities of (\ref{gauss-mix}) defined on $\XX=\YY=\RR^d$ satisfies Condition $\Ha_{\delta}(U,V)$ for some $\delta\in ]0,1[$ provided that
\begin{equation}\label{condition-multivariate-mixture} 
\beta^{\prime}\tau^{-1}\beta<(\delta\wedge (1-\delta))~\sigma_i^{-1}\quad \mbox{and}\quad
\tau^{-1}<(\delta\wedge (1-\delta))~\overline{\sigma}_j^{-1}
\end{equation}
for almost every  $(i,j)\in (E_0\times E_1)$.
\end{prop}

Let $(X_i)_{1\leq i\leq N}$ be $N$ independent and identically distributed samples from some (arbitrary) data distribution 
$\eta(dx)$  on $\XX=\RR^d$. The Gaussian-kernel estimation model defined by
 \begin{equation}\label{gaussian-mix-est}
\lambda_{U}(dx):=\frac{1}{N}\sum_{1\leq i\leq N}g_{\sigma}(x-X_i)~dx
 \end{equation}
is a well known Gaussian mixture approximation of the distribution $\eta(dx)$ (as the number of samples $N\rightarrow\infty$). The covariance matrix $\sigma$ is a bandwidth parameter that controls the amount and orientation of the smoothing. Note that (\ref{gaussian-mix-est}) coincides with (\ref{mix-intro}) when $(m_i,\sigma_i)=(X_i,\sigma)$ with the counting measure $\xi_0(i)=1/N$ on the index set $E_0=\{1,\ldots,N\}$. For a more detailed discussion on these non parametric statistical estimation techniques we refer to the pioneering articles~\cite{parzen,rosenblatt} as well as the books~\cite{chacon,simonov,wand}.

The Gaussian mixture distribution (\ref{gaussian-mix-est}) can be used to estimate complex and unknown data distributions  arising in generative models, see for instance~\cite{doucet-bortoli} and references therein.

\begin{rmk}\label{rmk-er}
Schr\"odinger bridges and the Sinkhorn algorithm are often used in practice to solve high dimensional optimal transport problems. In our context, choosing $\tau=tI$ and $(\alpha,\beta)=(0,I)$ in (\ref{def-W}), the Schr\"odinger bridge in the limit $t\rightarrow 0$ of vanishing regularization solves the Monge-Kantorovich problem and converge towards the optimal $2$-Wasserstein coupling. In this case, $\Sa_0=Q$ coincides with the heat kernel on a time horizon $t$ (see also Remark 2.3 in~\cite{adm-24}), namely,
$$
\Qa(x,dy)=(2\pi t)^{-d/2}~\exp{\left(-\frac{1}{2t}~\Vert y-x\Vert^2\right)}~dy.
$$
with the Euclidian norm $\Vert z\Vert^2=z^{\prime}z$.
In this context, the l.h.s. condition in (\ref{condition-multivariate-mixture}) applied to (\ref{gaussian-mix-est}) becomes
\begin{equation}\label{condition-multivariate-mixture-ex1} 
\frac{1}{t}~I<(\delta\wedge (1-\delta))~\sigma^{-1}.
\end{equation}
%in the limit $t\rightarrow 0$ of vanishing regularization the Schr\"odinger bridge solves the Monge-Kantorovich problem and converge towards the optimal $2$-Wasserstein coupling.
%Note that in this case $\Sa_0=Q$ coincides with the heat kernel on a time horizon $t$ (see also Remark 2.3 in~\cite{adm-24}), namely
%$$
%\Qa(x,dy)=(2\pi t)^{-d/2}~\exp{\left(-\frac{1}{2t}~\Vert y-x\Vert^2\right)}~dy.
%$$
%In this context, the l.h.s. condition in (\ref{condition-multivariate-mixture}) applied to (\ref{gaussian-mix-est}) becomes
%\begin{equation}\label{condition-multivariate-mixture-ex1} 
%\frac{1}{t}~I<(\delta\wedge (1-\delta))~\sigma^{-1}.
%\end{equation}
\end{rmk}

In many cases, the Markov transition $\Qa(x,dy)$ of interest is generated by a stochastic diffusion between time $0$ and a time horizon $t$ --see, e.g., \cite{leonard,stromme2023sampling,kato2024large}. In our context, if we choose $(\alpha,\beta)=(0,e^{tA})$ and $\tau=\tau_t:=\int_0^t e^{sA}e^{sA^{\prime}}ds$ for some Hurwitz matrix $A$ and some $t>0$, $\Qa(x,dy)$ coincides with the transition semigroup of the Ornstein-Uhlenbeck diffusion (see Remark 2.3 in~\cite{adm-24}). In this case, the l.h.s. condition in (\ref{condition-multivariate-mixture}) applied to (\ref{gaussian-mix-est}) is met for any $t\geq t_0$ provided that
\begin{equation}\label{condition-multivariate-mixture-ex2} 
e^{tA^{\prime}}\tau_{t_0}^{-1}e^{tA}<(\delta\wedge (1-\delta))~\sigma^{-1}.
 \end{equation}
Both transitions gradually transform the data distribution $\lambda_U$ into a simpler distribution $\lambda_U\Qa$ by injecting noise into the data. 
Also note that  (\ref{condition-multivariate-mixture-ex1}) and (\ref{condition-multivariate-mixture-ex2}) are met when the time horizon is chosen sufficiently large.

The choice of the prior distribution $\nu_V$ is rather flexible. For instance, we can choose a Gaussian distribution, say $\nu_{0,\overline{\sigma}}$, satisfying the r.h.s. condition in (\ref{cc-prop}). 
Starting from such Gaussian random states, the backward Sinkhorn transitions (indexed by odd indices) eventually generate samples back into the initial distribution by slowly removing the noise distortion. Conditions (\ref{condition-multivariate-mixture-ex1}) and (\ref{condition-multivariate-mixture-ex2}) reflect the fact that the forward noise injection as well as the reverse denoising procedure in diffusion-based generative models often require very large time horizons.

%%%%%%%%%%%%%%%%%%%%%%%%%%%%%%%%%%%%%%%%%%%%%%%%%%
% Acknowledgements
%%%%%%%%%%%%%%%%%%%%%%%%%%%%%%%%%%%%%%%%%%%%%%%%%%
\section*{Acknowledgements}

We would like to thank the anonymous reviewers for their excellent suggestions for improving the paper. Their detailed comments greatly improved the presentation of the article.

This work is supported by the Innovation and Talent Base for Digital Technology and Finance (B21038) and ‘the Fundamental Research Funds for the Central Universities’, Zhongnan University of Economics and Law (2722023EJ002).

J. Miguez acknowledges the financial support of the Office of Naval Research (award no. N00014-22-1-2647), MICIU/AEI/10.13039/501100011033/FEDER UE (grants no. PID2021-125159NB-I00 TYCHE and PID2024-158181NB-I00 NISA), and Community of Madrid, under grant IDEA-CM (TEC-
2024/COM-89).

%%%%%%%%%%%%%%%%%%%%%%%%%%%%%%%%%%%%%%%%%%%%%%
%% Appendices
%%%%%%%%%%%%%%%%%%%%%%%%%%%%%%%%%%%%%%%%%%%%%%

\appendix

\section*{Appendix}
This appendix compiles technical material and supplementary results that bolster the paper's central arguments.
Appendix~\ref{GS-sec} focuses on Gaussian Sinkhorn algorithms. This material is taken from~\cite{adm-24}. We review some properties of the Riccati matrix differential equations and the random maps associated with Sinkhorn transitions.  We also provide equivalent formulations of the drift-minorization conditions for the Lyapunov functions presented  in Proposition~\ref{prop1-lyap-gauss}.
Appendix~\ref{apHilbert} is devoted to basic properties of Hilbert's projective metric.
Appendix~\ref{phi-psinorm-sec} recalls several analytical tools related to weighted Banach spaces, operator norms, and contraction coefficients that are used throughout the stability analysis. For a more thorough discussion on weighted norms and contraction coefficients, we refer the reader to the articles~\cite{dg-01,dm-03,horton-dp,delgerber2025} as well as the book~\cite{penev}.

Appendix~\ref{appSchrodinger} reviews Sinkhorn potentials. The first part, Appendix~\ref{appSchrodinger-com}, presents some commutation formulae.
Building on this framework,  several quantitative estimates are provided in Appendix~\ref{appSchrodinger-quant}. Section~\ref{appSchrodinger-lyap} is dedicated to 
Lyapunov inequalities  and minorization conditions. The Proof of Theorem~\ref{lm-tx-2} is provided in Appendix~\ref{lm-tx-2-proof}
and Appendix~\ref{tech-proof-ap} contains further technical arguments.

%%%%%%%%%%%%%%%%%%%%%%%%%%%%%
% The Gaussian Sinkhorn 	%
% algorithm					%
%%%%%%%%%%%%%%%%%%%%%%%%%%%%%
\section{The Gaussian Sinkhorn algorithm}\label{GS-sec}
  
The Frobenius matrix norm of a given matrix $v\in \RR^{d_1\times d_2}$ is defined by
$
\left\Vert v\right\Vert_{F}^2=\tr(v^{\prime}v)
$, 
with the trace operator $\tr(\cdot)$ and  $v^{\prime}$ the transpose
of the matrix $v$. Note that, for $x\in\RR^{d \times 1}=\RR^d$, the Frobenius norm $\Vert x\Vert_F=\sqrt{x^{\prime}x}$ coincides with the Euclidean norm. The spectral norm is defined by $\Vert v\Vert_2=\sqrt{\lambda_{\sf max}(v^{\prime}v)}$. 
 
%%%
%
%%%
\subsection*{Riccati difference equations}
 
Let us associate with some  $\varpi\in\Sa^+_d$  the increasing map $\mbox{\rm Ricc}_{\varpi}$ from $\Sa^0_d$ into  {$\Sa^+_d$} defined by
\begin{equation}\label{ricc-maps-def}
\begin{array}{lccl}
\mbox{\rm Ricc}_{\varpi}: &\Sa^0_d &\mapsto &\Sa^+_d\\
&v &\leadsto &\mbox{\rm Ricc}_{\varpi}(v):=(I+(\varpi+v)^{-1})^{-1}\\ 
\end{array}.
\end{equation}
If we define the matrices
$$
 \varpi^{-1}:= \gamma\,\gamma^{\prime}
 \quad\mbox{\rm and}\quad
 \overline{\varpi}^{-1}:=\gamma^{\prime}\gamma\in\Sa^+_{d},
\quad\mbox{\rm with}\quad  \gamma:=\overline{\sigma}^{1/2}~\tau^{-1}\beta~\sigma^{1/2},
$$
and $\upsilon_{0}:=\overline{\sigma}^{-1/2}\tau~\overline{\sigma}^{-1/2}$, then it can be shown (see Theorem 4.3 in~\cite{adm-24}) that 
the recursions
\begin{equation}\label{k-ricc-0}
\upsilon_{2n+1}^{-1}=
I+  \gamma^{\prime}~\upsilon_{2n}~ \gamma
\quad \mbox{and}\quad
\upsilon_{2(n+1)}^{-1}=I+ \gamma~\upsilon_{2n+1}~  \gamma^{\prime}, 
\end{equation}
satisfy the Riccati difference equations
\begin{equation}\label{k-ricc}
\upsilon_{2(n+1)}=\mbox{\rm Ricc}_{\varpi}\left(\upsilon_{2n}\right)\quad
\mbox{and}\quad
\upsilon_{2n+1}=\mbox{\rm Ricc}_{ \overline{\varpi}}\left(\upsilon_{2n-1}\right), \quad \text{for $n\ge 0$}.
\end{equation} 
Let us also define
$$
\tau_{2n}:=
\overline{\sigma}^{1/2}~\upsilon_{2n}~\overline{\sigma}^{1/2}
\quad\mbox{and}\quad
\tau_{2n+1}:=
\sigma^{1/2}~\upsilon_{2n+1}~\sigma^{1/2},
$$
as well as the gain matrices
$$
\beta_{2n}:=\tau_{2n}~\tau^{-1}\beta\quad \mbox{and}\quad
\beta_{2n+1}:=\tau_{2n+1}~\beta^{\prime}\tau^{-1}.
$$
The recursions for $\tau_n$ and $\beta_n$ are initialized with $(\beta_0,\tau_0)=(\beta,\tau)$, where $(\beta,\tau)$ are the parameters of the reference Gaussian measure. 

%%%
%
%%%
\subsection*{Sinkhorn transitions}
 
We also denote by $G$ a $d$-dimensional centred Gaussian random variable (r.v.) with unit covariance. The Sinkhorn transitions in (\ref{s-2}) can be expressed as
$$
\Sa_{2n}(x,dy)=\PP\left(\Za_{2n}(x)\in dy\right)\quad\mbox{\rm and}\quad
\Sa_{2n+1}(y,dx)=\PP\left(\Za_{2n+1}(y)\in dx\right)
$$
with random maps defined, for any $n\geq 0$, by
\begin{eqnarray*}
\Za_{2n}(x)&=&m_{2n}+\beta_{2n}(x-m)+\tau_{2n}^{1/2}~G \quad \text{and}\\
\Za_{2n+1}(y)
&=&m_{2n+1}+\beta_{2n+1}(y-\overline{m})+\tau_{2n+1}^{1/2}~G.
\end{eqnarray*}
The parameters $m_n$ in the latter expressions are defined, for any $n\geq 0$, by the recursion
\begin{eqnarray*}
m_{2n+1}&=&m+
\beta_{2n+1} (\overline{m}-m_{2n}),\\
m_{2(n+1)}&=&\overline{m}+\beta_{2(n+1)} ~(m-m_{2n+1})
\end{eqnarray*}
with initial condition
$
m_0=\alpha_0+\beta_0 m$, where, again, $(\alpha_0,\beta_0)=(\alpha,\beta)$ are parameters of the reference Gaussian measure. Note that
$
\pi_n=\nu_{m_{n},\sigma_n}
$
with
$$\sigma_{2n}:=
\beta_{2n}~\sigma~\beta_{2n}^{\prime}+\tau_{2n}
\quad \mbox{and}\quad\sigma_{2n+1}:=\beta_{2n+1}~\overline{\sigma}~\beta_{2n+1}^{\prime}+\tau_{2n+1}.
$$

%%%
%
%%%
\subsection*{Gibbs-loop process}
The transitions of  the Gibbs-loop process introduced in Section~\ref{GL-sec}
 can also be rewritten as
$$
 \Sa^{\circ}_{2n}(y_1,dy_2):=\PP\left(\Za^{\circ}_{2n}(y_1)\in dy_2\right)
 \quad\mbox{\rm and}\quad 
  \Sa^{\circ}_{2n+1}(x_1,dx_2):=\PP\left(\Za^{\circ}_{2n+1}(x_1)\in dx_2\right)
$$
with the random maps
\begin{eqnarray*}
\Za^{\circ}_{2n}(y)&:=&\overline{m}+\beta^{\circ}_{2n}(y-\overline{m})+\left(\tau^{\circ}_{2n}\right)^{1/2}~G\\
\Za^{\circ}_{2n+1}(x)&:=&
m+\beta_{2n+1}^{\circ}(x-m)+\left(\tau_{2n+1}^{\circ}\right)^{1/2}~G
\end{eqnarray*}
defined in terms of the parameters
$$
\begin{array}{rclcrcl}
\beta_{2n+1}^{\circ}&:=&\beta_{2n+1}~\beta_{2n},&&
\tau_{2n+1}^{\circ}&:=&\tau_{2n+1}+\beta_{2n+1}~\tau_{2n}~
\beta_{2n+1}^{\prime},\\
&&&&&&\\
\beta^{\circ}_{2n}&:=&\beta_{2n}~\beta_{2n-1},&&
\text{and}~~\tau^{\circ}_{2n}&:=&\tau_{2n}+
\beta_{2n}~\tau_{2n-1}~\beta_{2n}^{\prime}.
\end{array}$$

By Lemma 4.9 in~\cite{adm-24}, we have
$$
0\leq \overline{\sigma}^{-1/2}~\beta^{\circ}_{2(n+1)}~\overline{\sigma}^{1/2}=I-\upsilon_{2(n+1)}\leq (I+\varpi)^{-1}
$$
as well as
$$
0\leq \sigma^{-1/2}~\beta_{2n+1}^{\circ}~\sigma^{1/2}
=I-\upsilon_{2n+1}\leq (I+\overline{\varpi})^{-1}.
$$
Moreover, Corollary 4.4 in~\cite{adm-24} yields the uniform estimates 
$$
\overline{\sigma}^{1/2}(I+\varpi^{-1})^{-1}\overline{\sigma}^{1/2}\leq \tau_{2n}\leq \overline{\sigma}\quad\mbox{and}\quad
\sigma^{1/2}(I+\overline{\varpi}^{-1})^{-1}\sigma^{1/2}\leq \tau_{2n+1}\leq \sigma, ~~\forall n\ge 1.
$$

%%%
%
%%%
\subsection*{Lyapunov functions}\label{lyap-sec-Gauss}
From the fixed point equations (\ref{fixed-points-gibbs}) we have
\begin{equation}
\overline{\sigma}=\beta^{\circ}_{2n}\overline{\sigma}\left(\beta^{\circ}_{2n}\right)^{\prime}+\tau^{\circ}_{2n}\geq \tau^{\circ}_{2n}
\quad \mbox{and}\quad
\sigma=\beta^{\circ}_{2n+1}\sigma\left(\beta^{\circ}_{2n+1}\right)^{\prime}+\tau^{\circ}_{2n+1}\geq \tau^{\circ}_{2n+1},
\label{eq2stars}
\end{equation}
which, in turn, imply that
\begin{equation}
\overline{\sigma}_{\varpi}:=\overline{\sigma}^{1/2}(I+\varpi^{-1})^{-1}\overline{\sigma}^{1/2}\leq \tau_{2n}^{\circ}\leq \overline{\sigma}
\quad\mbox{and}\quad\sigma_{\overline{\varpi}}:=
\sigma^{1/2}(I+\overline{\varpi}^{-1})^{-1}\sigma^{1/2}\leq \tau_{2n+1}^{\circ}\leq \sigma.
\label{eq1star}
\end{equation}

Also note that if $\overline{\sigma}^{-1}\leq (\tau_{2n}^{\circ})^{-1}$ then $(\tau_{2n}^{\circ})^{1/2}\overline{\sigma}^{-1}(\tau_{2n}^{\circ})^{1/2}\leq I$. Moreover, we also readily see that
$$
\begin{array}{l}
\EE\left[(\Za^{\circ}_{2n}(y)-\overline{m})^{\prime}\overline{\sigma}^{-1}(\Za^{\circ}_{2n}(y)-\overline{m})\right]\\
\\
=(y-\overline{m})^{\prime}\left(\beta^{\circ}_{2n}\right)^{\prime}\overline{\sigma}^{-1}
\beta^{\circ}_{2n}(y-\overline{m})
+\tr\left(\overline{\sigma}^{-1/2}\tau_{2n}^{\circ}\overline{\sigma}^{-1/2}\right)\\
\\
=(y-\overline{m})^{\prime}\overline{\sigma}^{-1/2}\left[\left(\overline{\sigma}^{-1/2}\beta^{\circ}_{2n}\overline{\sigma}^{1/2}\right)^{\prime}
\left(\overline{\sigma}^{-1/2}\beta^{\circ}_{2n}\overline{\sigma}^{1/2}\right)\right]\overline{\sigma}^{-1/2}(y-\overline{m})
+\tr\left(\overline{\sigma}^{-1/2}\tau_{2n}^{\circ}\overline{\sigma}^{-1/2}\right),
\end{array}$$
which implies 
\begin{eqnarray}
\EE\left[(\Za^{\circ}_{2n}(y)-\overline{m})^{\prime}\overline{\sigma}^{-1}(\Za^{\circ}_{2n}(y)-\overline{m})\right]
&=& (y-\overline{m})^{\prime}\overline{\sigma}^{-1/2}\left[I-\upsilon_{2n}\right]^2\overline{\sigma}^{-1/2}(y-\overline{m}) \nonumber \\
&& +\tr\left(\overline{\sigma}^{-1/2}\tau_{2n}^{\circ}\overline{\sigma}^{-1/2}\right) \nonumber \\
&\leq&  (1+\lambda_{\sf min}(\varpi))^{-2}~(y-\overline{m})^{\prime}\overline{\sigma}^{-1}(y-\overline{m}) +d.
\label{eqEstimate}
\end{eqnarray}
We are now in position to prove Proposition~\ref{prop1-lyap-gauss}.

%%%
%
%%%
\subsection*{Proof of Proposition~\ref{prop1-lyap-gauss}}

Consider the function
$$
\psi(y):=1+(y-\overline{m})^{\prime}\overline{\sigma}^{-1}(y-\overline{m}).
$$
The estimate in \eqref{eqEstimate} ensures that
$$
\Sa_{2n}^\circ(\psi)(y)=\EE\left[\psi(\Za^{\circ}_{2n}(y))\right] 
\leq
(1+\lambda_{\sf min}(\varpi))^{-2}~\psi(y)+\left(
    \left(1- (1+\lambda_{\sf min}(\varpi))^{-2}\right)+d
\right).
$$

Next, let us observe that
$$
\begin{array}{l}
 \displaystyle \Sa^{\circ}_{2n}(y_1,dy_2) =\frac{1}{\sqrt{\mbox{\rm det}(2\pi \tau_{2n}^{\circ})}}\\
 \\
 \displaystyle\times~\exp{\left(-\frac{1}{2}~\left((y_2-\overline{m})-\beta^{\circ}_{2n}(y_1-\overline{m}))\right)^{\prime}\left(\tau^{\circ}_{2n}\right)^{-1}\left((y_2-\overline{m})-\beta^{\circ}_{2n}(y_1-\overline{m}))\right)\right)}~dy_2
\end{array}$$
and
$$
\begin{array}{l}
\displaystyle\left((y_2-\overline{m})-\beta^{\circ}_{2n}(y_1-\overline{m}))\right)^{\prime}\left(\tau^{\circ}_{2n}\right)^{-1}\left((y_2-\overline{m})-\beta^{\circ}_{2n}(y_1-\overline{m}))\right)\\
\\
\displaystyle \leq 2
\left(y_2-\overline{m}\right)^{\prime}\left(\tau^{\circ}_{2n}\right)^{-1}\left(y_2-\overline{m}\right)+2\left(\beta^{\circ}_{2n}(y_1-\overline{m}))\right)^{\prime}\left(\tau^{\circ}_{2n}\right)^{-1}\left(\beta^{\circ}_{2n}(y_1-\overline{m}))\right)\\
\\
\displaystyle \leq 2
\left(y_2-\overline{m}\right)^{\prime}\overline{\sigma}_{\varpi}^{-1}\left((y_2-\overline{m})\right)+2~(1+\lambda_{\sf min}(\varpi^{-1}))~\left(y_1-\overline{m}\right)^{\prime}(\beta^{\circ}_{2n})^{\prime}\overline{\sigma}^{-1}\beta^{\circ}_{2n}\left(y_1-\overline{m}\right),
\end{array}
$$
where the second inequality follows from \eqref{eq1star}. Moreover, arguing again as above we arrive at
\begin{equation}
\begin{array}{l}
\displaystyle\left((y_2-\overline{m})-\beta^{\circ}_{2n}(y_1-\overline{m}))\right)^{\prime}\left(\tau^{\circ}_{2n}\right)^{-1}\left((y_2-\overline{m})-\beta^{\circ}_{2n}(y_1-\overline{m}))\right)\\
\\
\displaystyle \leq 2
\left(y_2-\overline{m}\right)^{\prime}\overline{\sigma}_{\varpi}^{-1}\left(y_2-\overline{m}\right)+2~(1+\lambda_{\sf min}(\varpi^{-1}))(1+\lambda_{\sf min}(\varpi))^{-2}~\left(y_1-\overline{m}\right)^{\prime}\overline{\sigma}^{-1}\left(y_1-\overline{m}\right).
\end{array}
\label{eq3ang}
\end{equation}
The inequality \eqref{eq3ang} above together with \eqref{eq2stars} readily yields
$$
\begin{array}{l}
 \displaystyle \Sa^{\circ}_{2n}(y_1,dy_2) \geq \frac{1}{\sqrt{\mbox{\rm det}(2\pi \overline{\sigma})}}~\\
 \\
 \displaystyle\times \exp{\left(-(1+\lambda_{\sf min}(\varpi^{-1}))(1+\lambda_{\sf min}(\varpi))^{-2}~\left(y_1-\overline{m}\right)^{\prime}\overline{\sigma}^{-1}\left(y_1-\overline{m}\right)\right)}\\
 \\
 \displaystyle\times~\exp{\left(-\left(y_2-\overline{m}\right)^{\prime}\overline{\sigma}_{\varpi}^{-1}\left(y_2-\overline{m}\right)\right)} ~dy_2.
\end{array}$$
Thus for any $y_1$ such that
$$
\psi(y_1)=1+(y_1-\overline{m})^{\prime}\overline{\sigma}^{-1}(y_1-\overline{m})\leq r
$$
we have
$$
\begin{array}{l}
 \displaystyle \Sa^{\circ}_{2n}(y_1,dy_2) \geq \varepsilon(r)~\nu_{\overline{m},\overline{\sigma}_{\varpi}/2}(dy_2),
\end{array}$$
where
$$
\varepsilon(r)=\frac{1}{\sqrt{\mbox{\rm det}(2\overline{\sigma}\overline{\sigma}_{\varpi}^{-1})}}~\\
 \\
 \displaystyle\times \exp{\left(-(1+\lambda_{\sf min}(\varpi^{-1}))(1+\lambda_{\sf min}(\varpi))^{-2}~(r-1)\right)}.
$$
This shows that the function $\psi$, and similarly the function $\varphi$  defined in (\ref{first-set-glyap}) satisfy 
the drift-minorization conditions stated in (\ref{lyap-eq-sink-o}) and (\ref{sink-locmin-o}). \cqfd

%%%
%
%%%
\subsection*{Lyapunov functions in \eqref{second-set-glyap} and \eqref{ref-gauss-lyap}}

Using the triangle inequality for Mahalanobis distances we check that
$$
\begin{array}{l}
\left[(\Za^{\circ}_{2n}(y)-\overline{m})^{\prime}\overline{\sigma}^{-1}(\Za^{\circ}_{2n}(y)-\overline{m})\right]^{1/2}\\
\\
\leq 
\left[(\beta^{\circ}_{2n}(y-\overline{m}))^{\prime}\overline{\sigma}^{-1}(\beta^{\circ}_{2n}(y-\overline{m}))\right]^{1/2}+
\left[(\left(\tau^{\circ}_{2n}\right)^{1/2}~G)^{\prime}\overline{\sigma}^{-1}(\left(\tau^{\circ}_{2n}\right)^{1/2}~G)\right]^{1/2}\\
\\
\leq  (1+\lambda_{\sf min}(\varpi))^{-1}~((y-\overline{m})^{\prime}\overline{\sigma}^{-1}(y-\overline{m}))^{1/2}+
\left[G^{\prime}G\right]^{1/2},
\end{array}$$
which yields the almost sure estimate
$$
\left[(\Za^{\circ}_{2n}(y)-\overline{m})^{\prime}\overline{\sigma}^{-1}(\Za^{\circ}_{2n}(y)-\overline{m})\right]^{1/2}
\leq 
(1+\lambda_{\sf min}(\varpi))^{-1}~((y-\overline{m})^{\prime}\overline{\sigma}^{-1}(y-\overline{m}))^{1/2}+
\left[G^{\prime}G\right]^{1/2}.
$$ 
This implies  that the function $\psi$, and similarly the function $\varphi$  defined in (\ref{second-set-glyap}) satisfy 
the drift-minorization conditions stated in (\ref{lyap-eq-sink-o}) and (\ref{sink-locmin-o}).

The analysis of the Lyapunov functions (\ref{ref-gauss-lyap}) is based on the following technical lemma.
\begin{lem}\label{lem-tex-2}
For any matrices $\sigma>\tau$ and any $z\in\RR^d$ and $0\leq \delta\leq 1$ we have
$$
\EE\left[\exp{\left(\frac{\delta}{2}(z+\tau^{1/2}G)^{\prime}\sigma^{-1}(z+\tau^{1/2}G)\right)}\right]
=\displaystyle\sqrt{\frac{\mbox{\rm det}(\sigma)}{\mbox{\rm det}(\sigma-\delta\tau)}}~
\exp{\left(\frac{\delta}{2}z^{\prime}(\sigma-\delta\tau)^{-1}z\right)}.
$$
\end{lem} 
\proof
For any $g\in\RR^d$ we have
$$
\frac{\delta}{2}~(z+g)^{\prime}\sigma^{-1}(z+g)-\frac{1}{2}~g^{\prime}\tau^{-1}g
=\frac{\delta}{2}~z^{\prime}\sigma^{-1}z+\frac{\delta^2}{2}~z^{\prime}\sigma^{-1}\cchi\sigma^{-1}z-\frac{1}{2}~(g-\delta\cchi\sigma^{-1}z)^{\prime}\cchi^{-1}
(g-\delta\cchi\sigma^{-1}z),
$$
with
$$
0<\cchi:=(\tau^{-1}-\delta \sigma^{-1})^{-1}
=-\frac{1}{\delta}~\sigma+\frac{1}{\delta}\sigma\tau^{-1}(\tau^{-1}-\delta \sigma^{-1})^{-1}.
$$
We note that
\begin{eqnarray*}
\sigma^{-1}+\delta \sigma^{-1}\cchi\sigma^{-1}&=&
\sigma^{-1}+\delta \sigma^{-1}\left(-\frac{1}{\delta}~\sigma+\frac{1}{\delta}\sigma\tau^{-1}(\tau^{-1}-\delta \sigma^{-1})^{-1}\right)\sigma^{-1}\\
&=&\tau^{-1}(\tau^{-1}-\delta \sigma^{-1})^{-1}\sigma^{-1}=(\sigma-\delta\tau)^{-1}
\end{eqnarray*}
and finally observe that
$$
\frac{1}{\mbox{\rm det}(\tau \cchi^{-1})}
=\frac{\mbox{\rm det}(\sigma)}{\mbox{\rm det}(\sigma-\delta\tau)}.
$$
\cqfd

Using Lemma~\ref{lem-tex-2} we check that
$$
\begin{array}{l}
\displaystyle\EE\left[\exp{\left(\frac{\delta}{2}\left[(\Za^{\circ}_{2n}(y)-\overline{m})^{\prime}\overline{\sigma}^{-1}(\Za^{\circ}_{2n}(y)-\overline{m})\right]\right)}\right]\\
\\
=\displaystyle\EE\left[\exp{\left(\frac{\delta}{2}\left[\left(\beta^{\circ}_{2n}(y-\overline{m})+\left(\tau^{\circ}_{2n}\right)^{1/2}~G\right)^{\prime}\overline{\sigma}^{-1}\left(\beta^{\circ}_{2n}(y-\overline{m})+\left(\tau^{\circ}_{2n}\right)^{1/2}~G\right)\right]\right)}\right]\\
\\
=\displaystyle
\frac{1}{\sqrt{\mbox{\rm det}(I-\delta~ \overline{\sigma}^{-1/2}\tau^{\circ}_{2n}\overline{\sigma}^{-1/2})}}~\exp{\left(\frac{\delta}{2}~(y-\overline{m})^{\prime}~(\beta^{\circ}_{2n})^{\prime}\left(\overline{\sigma}-\delta~\tau^{\circ}_{2n}\right)^{-1}\beta^{\circ}_{2n}~(y-\overline{m})\right)}.
\end{array}
$$
On the other hand, we have
$$
 \overline{\sigma}^{-1/2}\tau^{\circ}_{2n}\overline{\sigma}^{-1/2}\leq I
$$
and arguing as above we conclude that
\begin{eqnarray}
(\beta^{\circ}_{2n})^{\prime}\left(\overline{\sigma}-\delta~\tau^{\circ}_{2n}\right)^{-1}\beta^{\circ}_{2n}
&=&
(\overline{\sigma}^{-1/2}\beta^{\circ}_{2n})^{\prime}\left(I-\delta~\overline{\sigma}^{-1/2}\tau^{\circ}_{2n}\overline{\sigma}^{-1/2}\right)^{-1}\overline{\sigma}^{-1/2}\beta^{\circ}_{2n}\\
&\leq& \frac{1}{1-\delta}~(\beta^{\circ}_{2n})^{\prime}\overline{\sigma}^{-1}\beta^{\circ}_{2n} \leq   \frac{1}{1-\delta} (1+\lambda_{\sf min}(\varpi))^{-2}~\overline{\sigma}^{-1}.
\end{eqnarray}
If we now consider the function
$$\psi(y):=\exp{\left(\frac{\delta}{2}(y-\overline{m})^{\prime}\overline{\sigma}^{-1}(y-\overline{m})\right)}$$
then we have
\begin{eqnarray}
\frac{1}{\psi(y)}~\EE\left(
\psi\left(\Za^{\circ}_{2n}(y)\right)\right)%\\
%\\
&\leq& %\displaystyle
\frac{1}{\sqrt{\mbox{\rm det}(I-\delta~I)}}~\times \nonumber\\
&& \times~\exp{\left(-\frac{\delta}{2}~ \left(1-\frac{1}{1-\delta} (1+\lambda_{\sf min}(\varpi))^{-2}\right)~(y-\overline{m})^{\prime}~\overline{\sigma}^{-1}~(y-\overline{m})\right)}.
%\end{array}
\nonumber
\end{eqnarray}
Choosing
$$
0<\delta<1-\frac{1}{(1+\lambda_{\sf min}(\varpi))^2}
$$
we check that
$$
\frac{1}{(1-\delta)(1+\lambda_{\sf min}(\varpi))^2}<1\quad\text{and}\quad
\lim_{|y|\to\infty} \frac{1}{\psi(y)}~\Sa^{\circ}_{2n}(
\psi)(y) = 0.
$$
This implies  that the function $\psi$, and similarly the function $\varphi$, defined in (\ref{ref-gauss-lyap}) satisfy 
the drift-minorization conditions stated in (\ref{lyap-eq-sink-o}) and (\ref{sink-locmin-o}).

%%%%%%%%%%%%%%%%%%%%%%%%%%%%%
% Hilbert metric 			%
% contractions				%
%%%%%%%%%%%%%%%%%%%%%%%%%%%%%
\section{Hilbert metric contractions} \label{apHilbert}

We recall the Hilbert contraction coefficient formula (\cite{bauer,birkhoff,birkhoff-2,hopf}, see also~\cite{atar, cohen, legland})
\begin{equation}\label{defi-tau}
\cchi_H(\Ka)=\tanh{\left(\frac{1}{4}\log\sup_{(x_1,x_2)\in \XX^2}H(\delta_{x_1}\Ka,\delta_{x_2}\Ka)\right)}.
\end{equation}
%The proof of formula (\ref{defi-tau}) can be found in~\cite{bauer,birkhoff,birkhoff-2,hopf}, see also~\cite{atar, cohen, legland}. 
Also note that
$$
H(\delta_{x_1}\Ka,\delta_{x_2}\Ka)=\sup_{(y_1,y_2)\in  \YY^2}
\left(\imath_{x_1,x_2}(\Ka)(y_1)~
\imath_{x_2,x_1}(\Ka)(y_2)\right)
$$
with the Radon-Nikodym function $\imath_{x_1,x_2}(\Ka)(y)$ defined in (\ref{i-def}). 
We also have
$$
\sup_{(x_1,x_2)\in S^2}H(\delta_{x_1}\Ka,\delta_{x_2}\Ka)=1/\hbar(\Ka)
$$
with the parameter $\hbar(\Ka)$ defined in (\ref{it-Ka-def}).
Using the fact that
$$
\tanh\left(\frac{1}{4}~\log{u}\right)=\frac{u^{1/4}-u^{-1/4}}{u^{1/4}-u^{-1/4}}=\frac{u^{1/2}-1}{u^{1/2}+1}=\frac{1-u^{-1/2}}{1+u^{-1/2}}
$$
we conclude that
\begin{equation}\label{prop-tau}
\cchi_H(\Ka)=\frac{1-\sqrt{\hbar(\Ka)}}{1+\sqrt{\hbar(\Ka)}}.
\end{equation}

\begin{prop}\label{prop-h-bar-Kn}
Consider the collection of Sinkhorn transitions  $\Sa_n$ defined in
 (\ref{def-Pa-n}). For any $n\geq 0$,
 $$
\hbar(\Sa_n)=\hbar(\Qa)\quad \mbox{and}\quad
\cchi_H(\Sa_n)=\frac{1-\sqrt{\hbar(\Qa)}}{1+\sqrt{\hbar(\Qa)}}.
 $$ 
\end{prop}
\proof
For any $n\geq 0$, we have 
\begin{eqnarray}
\imath_{y_1,y_2}(\Sa_{2n+1})(x_1)~\imath_{y_2,y_1}(\Sa_{2n+1})(x_2)
&=&
\frac{(d\delta_{x_1}\Sa_{2n}/d\pi_{2n})(y_1)}{(d\delta_{x_1}\Sa_{2n}/d\pi_{2n})(y_2)}\frac{(d\delta_{x_2}\Sa_{2n}/d\pi_{2n})(y_2)}{(d\delta_{x_2}\Sa_{2n}/d\pi_{2n})(y_1)} \nonumber\\
&=&\imath_{x_2,x_1}(\Sa_{2n})(y_2)~\imath_{x_1,x_2}(\Sa_{2n})(y_1) \nonumber
\end{eqnarray}
which, together with \eqref{it-Ka-def}, yields
$$
\hbar(\Sa_{2n+1})=
\hbar(\Sa_{2n}).
$$
In the same vein, for any $n\geq 1$ we have
\begin{eqnarray*}
\imath_{x_1,x_2}(\Sa_{2n})(y_1)~\imath_{x_2,x_1}(\Sa_{2n})(y_2) 
&=&
\frac{(d\delta_{y_1}\Sa_{2n-1}/d\pi_{2n-1})(x_1)}{(d\delta_{y_1}\Sa_{2n-1}/d\pi_{2n-1})(x_2)}~\frac{(d\delta_{y_2}\Sa_{2n-1}/d\pi_{2n-1})(x_2)}{(d\delta_{y_2}\Sa_{2n-1}/d\pi_{2n-1})(x_1)}\\
&=&
\imath_{y_2,y_1}(\Sa_{2n-1})(x_2)~\imath_{y_1,y_2}(\Sa_{2n-1})(x_1),
\end{eqnarray*}
hence
$$
\hbar(\Sa_{2n})=
\hbar(\Sa_{2n-1}).
$$
Recalling that $\Sa_0=\Qa$, this ends the proof of the proposition.
\cqfd

%%%%%%%%%%%%%%%%%%%%%%%%%%%%%
% Weighted norm	 			%
% contractions				%
%%%%%%%%%%%%%%%%%%%%%%%%%%%%%
\section{Weighted norm contractions}\label{phi-psinorm-sec}

We associate with some $\varphi\in \Ba(\XX)$ and $\psi\in \Ba(\YY)$ such that $\varphi(x)\wedge\psi(y)\geq 1$ the following condition.

{\it $\La_{\varphi,\psi}(\Ka)$
There exists some parameters $\epsilon\in ]0,1[$ and $c<\infty$ such that
 \begin{equation}\label{lyap-eq}
 \Ka(\psi)\leq \epsilon~\varphi+c.
 \end{equation}
In addition, there exist some constant $r_{0}\geq 1$ and some function 
$$
\begin{array}{ccccc}
\varepsilon &:&[r_0,\infty[ &\rightarrow &]0,1]\\
&&r &\mapsto &\varepsilon(r)\\
\end{array}
$$
such that, for any $r\geq r_0$ and $(x_1,x_2)\in  C_{\varphi}(r)^2$,
 \begin{equation}\label{minloc}
 \chi_{\varphi}(r,\Ka):=\sup_{(x_1,x_2)\in  C_{\varphi}(r)^2}\Vert (\delta_{x_1}-\delta_{x_2})\Ka\Vert_{\sf tv}<1-\varepsilon(r).
 \end{equation}
}
 
By (\ref{ref-coupling-tv}), condition (\ref{minloc}) is equivalent to the existence of some $\Ya\in \Ma_1(\YY)$ (that may depends on $(r,x_1,x_2)$) satisfying for $i\in\{1,2\}$ the local minorization condition 
 \begin{equation}\label{minloc-2}
 \Ka(x_i,dy)\geq \varepsilon(r)~\Ya(dy).
 \end{equation}
For instance, assume that $\delta_{x_1}\Ka\ll\delta_{x_2}\Ka$ for any $x_1,x_2\in \XX$ and set
 \begin{equation}\label{i-min}
\jmath_{\varphi,\psi}(r,\Ka):=\inf_{(x_1,x_2)\in  C_{\varphi}(r)}~
\inf_{y\in  C_{\psi}(r)}\imath_{x_1,x_2}(\Ka)(y).
 \end{equation}
In this notation, we see that if $\jmath_{\varphi,\psi}(r,\Ka)\geq \varepsilon(r)$ then \eqref{minloc-2} holds.

Note that the strength of conditions  (\ref{lyap-eq}), (\ref{minloc})  depends on the strength of the functions $(\varphi,\psi)$. When the functions are bounded, the geometric drift condition  (\ref{lyap-eq})  is  trivially met, but in this case (\ref{minloc-2})  is a uniform contraction condition on the space $\XX$.  
  
When $\varphi\in\Ba_{\infty}(\XX)$ the $r$-sub-level sets $C_{\varphi}(r)$ are compact and condition (\ref{minloc-2}) is met as soon as we have
\begin{equation}\label{ref-P-min}
\Ka(x,dy)\geq q(x,y)~\nu(dy)
\end{equation}
for some  Radon positive measure $\nu$ on $\YY$ and some density function $q$, satisfying the local minorization condition  
 \begin{equation}\label{min-p}
\inf_{C_{\varphi}(r)\times C_{\psi}(r)}q>0\quad \mbox{\rm and}\quad 0<\nu(C_{\psi}(r))<\infty
  \end{equation}
for any $r\geq r_0$. For locally compact Polish spaces condition $0<\nu(C_{\psi}(r))<\infty$ is satisfied if $\psi$ has compact sub-levels sets $C_{\psi}(r)$ with non empty interior and $\nu$ is a Radon measure of full support. Also note that the l.h.s. minorization condition (\ref{min-p}) is satisfied when $(x,y)\in (\XX\times\YY)^{\circ}\mapsto q(x,y)$ is a continuous positive function on the interior $(\XX\times\YY)^{\circ}:=(\XX^{\circ}\times\YY^{\circ})$ of   $(\XX\times\YY)$.
 Also note that
 \begin{equation}\label{ref-dom-theta-0} 
\text{if} \quad \Ka(\psi)/\varphi\leq \Theta\in \Ba_0(\XX) 
\quad \text{then (\ref{lyap-eq}) holds.}
\end{equation}
 Indeed, if the inequality above holds then for any $0<\epsilon< \Vert\Theta\Vert$, the set $A_{\epsilon}:=\{\Theta\geq \epsilon\}=C_{1/\Theta}(1/\epsilon)$ is a non empty compact subset and we have
\begin{equation}\label{ref-0-CC}
 \Ka(\psi)\leq \epsilon~ 1_{S-A_{\epsilon}} \varphi+
 1_{A_{\epsilon}}(\Theta~\varphi)\leq  \epsilon~ \varphi+c_{\epsilon}\quad \mbox{\rm with}\quad c_{\epsilon}:=\Vert\Theta\Vert
 \sup_{A_{\epsilon}}\varphi.
\end{equation}

\begin{theo}\label{theo-contract-V}
Assume Condition $\La_{\varphi,\psi}(\Ka)$ is satisfied for some uniformly positive functions $(\varphi,\psi)\in (\Ba(\XX)\times\Ba(\YY))$. Then, there exist constants $\kappa\geq 0$ and $\rho\in [0,1[$ that only depend on the parameters  $(\epsilon,c)$ in (\ref{lyap-eq}) and on the function $\varepsilon(r)$ in (\ref{minloc}) and such that
 \begin{equation}\label{minlyap-K}
  \cchi_{{\varphi_{\kappa},\psi_{\kappa}}}(\Ka)\leq \rho\quad \mbox{with}\quad (\varphi_{\kappa},\psi_{\kappa}):=
  (1+\kappa\varphi,1+\kappa\psi).
 \end{equation}
\end{theo}

Theorem \ref{theo-contract-V} is a consequence  of Lemma 2.3 in~\cite{horton-dp}, see also~\cite{delgerber2025}, as well as  Theorem 8.2.21 in~\cite{penev} and Lemma 2.1 in~\cite{dpa}. Several examples satisfying the Lyapunov condition (\ref{lyap-eq}) and the local minorization conditions (\ref{minlyap-K}) are discussed in~\cite{dpa,penev}. 

%%%
%
%%%
\subsection*{Proof of Theorem~\ref{theo-1-intro}}\label{theo-1-intro-proof}

We check that if Condition $\Da_{\varphi,\psi}(\Sa)$ is satisfied for some integer $n_0\geq 0$, then (\ref{lyap-eq-sink-o}) and (\ref{sink-locmin-o}) also hold for any $n>n_0$.

Assume that Condition $\Da_{\varphi,\psi}(\Sa)$ holds for some $n_0\geq 0$. From (\ref{lyap-eq-sink}), for any $n>n_0$ we have the estimates
 \begin{eqnarray}
 \Sa^{\circ}_{2n}(\psi)&=&
  \Sa_{2n-1} (\Sa_{2n}(\psi))\leq \epsilon~
    \Sa_{2n-1}(\varphi)+c\leq \epsilon^2~\psi+(1+\epsilon)~c \quad \text{and}\nonumber\\
 \Sa^{\circ}_{2n+1}(\varphi)    &=&  \Sa_{2n} (\Sa_{2n+1}(\varphi))\leq 
 \epsilon~
    \Sa_{2n}(\psi)+c\leq \epsilon^2~\varphi+(1+\epsilon)~c. \label{fint}
 \end{eqnarray}
 Also note that for any $n>n_0$ the inequality \eqref{sink-locmin} in Condition $\Da_{\varphi,\psi}(\Sa)$ implies that
 \begin{equation}
  \Vert ((\delta_{y_1}-\delta_{y_2}) \Sa_{2n-1}) \Sa_{2n}\Vert_{\sf tv}\vee  \Vert ((\delta_{x_1}-\delta_{x_2}) \Sa_{2n})\Sa_{2n+1}\Vert_{\sf tv}\leq  1-\varepsilon(r).
 \label{eqContract0}
 \end{equation}
 The estimates in \eqref{fint} imply that \eqref{lyap-eq-sink-o} holds, while the inequality \eqref{eqContract0} implies that \eqref{sink-locmin-o} is satisfied.

Moreover, combining the fixed point equations in \eqref{fixed-points-gibbs} with (\ref{fint}) we arrive at
\begin{equation}
 \sup_{n\geq n_0}\left(\Vert \Sa_{2n}(\psi)\Vert_{\varphi}\vee \Vert \Sa_{2n+1}(\varphi)\Vert_{\psi}\right)<\infty\quad\mbox{\rm and}\quad (\lambda_U(\varphi)\vee \nu_V(\psi))\leq \frac{c}{1-\epsilon},
 \label{eqFinite}
\end{equation}
and (\ref{lyap-eq-sink}) together with \eqref{eqFinite} yields
$$
 \sup_{n\geq n_0}(\pi_{2n}(\psi)\vee \pi_{2n+1}(\varphi))\leq c_1~(\lambda_U(\varphi)\vee \nu_V(\psi))
$$
for some $c_1>0$. Therefore, Condition  $\Da_{\varphi,\psi}(\Sa^{\circ})$ is satisfied.

Finally, we observe that \eqref{lyap-eq-sink} ensures that inequality \eqref{lyap-eq} is satisfied for the Markov transitions $\Sa_{2n}$ and $\Sa_{2n+1}$, while \eqref{sink-locmin} yields
$$
\chi_{\varphi}(r,\Sa_{2n})\vee   \chi_{\psi}(r,\Sa_{2n+1})  \leq 1-\varepsilon(r)
$$
and, hence, the inequality \eqref{minloc} also holds for $\Sa_{2n}$ and $\Sa_{2n+1}$. As a consequence, Theorem \ref{theo-contract-V} can be applied to show that there are parameters $\kappa\geq 0$ and $\rho\in [0,1[$ such that
$$
\cchi_{{\varphi_{\kappa},\psi_{\kappa}}}(\Sa_{2n}) \vee \cchi_{{\psi_{\kappa},\varphi_{\kappa}}}(\Sa_{2n+1})\leq \rho.
$$
\cqfd

%%%%%%%%%%%%%%%%%%%%%%%%%%%%%
% Schrödinger potential		%
% functions					%
%%%%%%%%%%%%%%%%%%%%%%%%%%%%%
\section{Sinkhorn potential functions} \label{appSchrodinger}

This appendix collects a number of technical results concerning the Sinkhorn potential functions associated with the alternating Sinkhorn dynamics. We first establish a set of commutation formulae linking the successive marginals and Radon–Nikodym derivatives along the iteration. These identities are then used to derive quantitative bounds on the Sinkhorn potentials and on the corresponding density ratios, both in uniform and weighted norms. Under suitable integrability and Lyapunov-type assumptions, we further obtain refined estimates leading to drift and local minorization inequalities for the Sinkhorn transition kernels. Together, these results provide the analytical backbone needed to verify the contraction and stability properties used in the proofs of the main theorems.

\subsection{Commutation formulae}\label{appSchrodinger-com}

For any $n\geq 0$ we have the commutation formulae
\begin{equation}\label{comm} 
\Sa_{2n}\left(\frac{d\nu_V}{d\pi_{2n}}\right)= \frac{d \pi_{2n+1}}{d\lambda_U}\quad \mbox{and}\quad
 \Sa_{2(n+1)}\left(\frac{d\pi_{2n}}{d\nu_V}\right)= \frac{d\lambda_U}{d \pi_{2n+1}}.
\end{equation}
In the same vein, for any $n\geq 1$ we have
\begin{equation}\label{comm-2} 
\Sa_{2n-1}\left(\frac{d\lambda_U}{d \pi_{2n-1}}\right)=\frac{d\pi_{2n}}{d\nu_V}
 \quad \mbox{and}\quad
   \Sa_{2n+1}\left(\frac{d \pi_{2n-1}}{d\lambda_U}\right)=\frac{d\nu_V}{d\pi_{2n}}.
\end{equation}
A proof for (\ref{comm}) and (\ref{comm-2}) is provided in Appendix \ref{tech-proof-ap} (on page~\pageref{form-series-proof}).

Using the commutation properties (\ref{comm}) and (\ref{comm-2}) for any $n\geq p$ we readily check that
\begin{eqnarray}
\frac{d\pi_{2n}}{d\nu_V}&=&\Sa^{\circ}_{2n}\left(\frac{d\pi_{2(n-1)}}{d\nu_V}\right)=\Sa^{\circ}_{2n,2p}\left(\frac{d\pi_{2p}}{d\nu_V}\right) \quad \text{and}
\nonumber \\
  \frac{d \pi_{2n+1}}{d\lambda_U}&=& \Sa_{2n+1}^{\circ}\left(\frac{d \pi_{2n-1}}{d\lambda_U}\right)=\Sa^{\circ}_{2n+1,2p-1}\left(\frac{d \pi_{2p-1}}{d\lambda_U}\right),\label{gibbs-ratio}
\end{eqnarray}
with the backward semigroups
$$
\Sa^{\circ}_{2n,2p}:=\Sa^{\circ}_{2n}\ldots \Sa^{\circ}_{2(p+1)}
\quad\mbox{\rm and}\quad
\Sa^{\circ}_{2n+1,2p-1}:=\Sa^{\circ}_{2n+1}\ldots \Sa^{\circ}_{2p+1}.
$$
By (\ref{comm}) and (\ref{comm-2}) we also obtain
\begin{eqnarray}
\frac{d\lambda_U}{d \pi_{2n+1}}
&=& \left(\Sa_{2(n+1)}\,\Sa^{\circ}_{2n,2p}\,\Sa_{2p-1}\right)\left(\frac{d\lambda_U}{d \pi_{2p-1}}\right) \quad \text{and}\nonumber \\
\frac{d\nu_V}{d\pi_{2n}}&=& \left(  \Sa_{2n+1}\,\Sa^{\circ}_{2n-1,2p-1}\,\Sa_{2(p-1)}\right)\left(\frac{d\nu_V}{d\pi_{2(p-1)}}\right).
\label{gibbs-ratio-bb}
\end{eqnarray}

\subsection{Some quantitative estimates}\label{appSchrodinger-quant}

\begin{prop}\label{prop-bound-jQ}
When $ \jmath(\Qa)>0$, for any $n\geq 0$ we have the uniform estimates
\begin{equation}\label{bound-jQ}
\Vert\frac{d\pi_{2n}}{d\nu_V}\Vert\vee \Vert\frac{d\nu_V}{d\pi_{2n}}\Vert\leq \hbar(\Qa)^{-2}\quad
\mbox{and}\quad
\Vert \frac{d\pi_{2n+1}}{d\lambda_U}\Vert\vee \Vert \frac{d\lambda_U}{d\pi_{2n+1}}\Vert\leq\hbar(\Qa)^{-2}.
\end{equation}
\end{prop}
\proof Recalling that $\jmath(\Sa_n)\wedge\jmath(\Qa)\geq \hbar(\Sa_n)= \hbar(\Qa)$, for any given $x\in\XX$ we have 
$$
\hbar(\Qa)~\Sa_{2n}(x,dy_2)\leq \Sa_{2n}^{\circ}(y_1,dy_2)=
(\Sa_{2n-1}\Sa_{2n})(y_1,dy_2)\leq \frac{1}{\hbar(\Qa)}~\Sa_{2n}(x,dy_2).
$$
This yields, for any $(y_1,\overline{y}_1)\in\YY^2$, the estimate
$$
\hbar(\Qa)^2~
\Sa_{2n}^{\circ}(\overline{y}_1,dy_2)\leq \Sa_{2n}^{\circ}(y_1,dy_2)\leq \frac{1}{\hbar(\Qa)^2}~
\Sa_{2n}^{\circ}(\overline{y}_1,dy_2)
$$
and integrating w.r.t. $\nu_V(d\overline{y}_1)$ and $\pi_{2(n-1)}(dy_1)$ we find that
$$
\hbar(\Qa)^2~
\nu_V(dy_2)\leq\pi_{2n}(dy_2)\leq \frac{1}{\hbar(\Qa)^2}~
\nu_V(dy_2).
$$
Similarly, integrating w.r.t. $\nu_V(dy_1)$ and $\pi_{2(n-1)}(d\overline{y}_1)$ we also obtain
$$
\hbar(\Qa)^2~\pi_{2n}(dy_2)
\leq\nu_V(dy_2)\leq \frac{1}{\hbar(\Qa)^2}~
\pi_{2n}(dy_2).
$$
From the inequalities above we conclude that
$$
\frac{d\pi_{2n}}{d\nu_V}\vee \frac{d\nu_V}{d\pi_{2n}}\leq \hbar(\Qa)^{-2}
\quad\text{and}\quad
 \frac{d\pi_{2n+1}}{d\lambda_U}\vee \frac{d\lambda_U}{d\pi_{2n+1}}\leq\hbar(\Qa)^{-2}. 
$$
\cqfd

More generally these Radon Nikodym functions are unbounded and we need to carry a more refined analysis.

\begin{prop}\label{prop-uest-pi-ratio}
If Condition $\Da_{\varphi,\psi}(\Sa^{\circ})$ holds for some $n_0\geq 1$ then there exists $c_{\epsilon}>0$ such that 
 \begin{equation}\label{uest-pi-ratio}
\sup_{n\geq n_0}\Vert \frac{d\pi_{2n}}{d\nu_V}\Vert_{\psi}\leq c_{\epsilon}~\Vert \frac{d\pi_{2(n_0-1)}}{d\nu_V}\Vert_{\psi}\quad \mbox{\rm and}\quad
\sup_{n\geq n_0}\Vert \frac{d\pi_{2n+1}}{d\lambda_U}\Vert_{\varphi}\leq c_{\epsilon}~\Vert \frac{d\pi_{2n_0-1}}{d\lambda_U}\Vert_{\varphi}.
\end{equation}
\end{prop}

\proof Using (\ref{gibbs-ratio}) and (\ref{lyap-eq-sink-o}), we see that, for any $n\geq n_0$,
\begin{eqnarray*}
\frac{d\pi_{2n}}{d\nu_V}&\leq&\Vert \frac{d\pi_{2(n_0-1)}}{d\nu_V}\Vert_{\psi}~\left(\Sa^{\circ}_{2n}\ldots \Sa^{\circ}_{2n_0}\right)\left(\psi\right)\\
&\leq& \Vert \frac{d\pi_{2(n_0-1)}}{d\nu_V}\Vert_{\psi}\left(\epsilon^{n-n_0}\psi+\frac{c}{1-\epsilon}\right)\leq c_{\epsilon}~
\Vert \frac{d\pi_{2(n_0-1)}}{d\nu_V}\Vert_{\psi}~ \psi,
\quad \text{with}~c_{\epsilon}:=\left(1+\frac{c}{1-\epsilon}\right).
\end{eqnarray*}
This implies that 
$$
\sup_{n\geq n_0}\Vert \frac{d\pi_{2n}}{d\nu_V}\Vert_{\psi}\leq c_{\epsilon}~\Vert \frac{d\pi_{2(n_0-1)}}{d\nu_V}\Vert_{\psi}\quad \text{and, similarly,}\quad
\sup_{n\geq n_0}\Vert \frac{d\pi_{2n+1}}{d\lambda_U}\Vert_{\varphi}\leq c_{\epsilon}~\Vert \frac{d\pi_{2n_0-1}}{d\lambda_U}\Vert_{\varphi}.
$$
\cqfd

\begin{prop}\label{prop-uest-pi-ratio-2}
In addition, when Condition $\Da_{\varphi,\psi}(\Sa)$ holds for some $n_0\geq 0$,  there exists $c_{\epsilon}>0$ such that 
$$
\sup_{n\geq n_0}\left(\Vert \frac{d\pi_{2n+1}}{d\lambda_U}\Vert_{\varphi}\vee \Vert\frac{d\lambda_U}{d\pi_{2n+1}}\Vert_{\varphi}\right)\leq c_{\epsilon}~
\left(\Vert \frac{d\pi_{2n_0+1}}{d\lambda_U}\Vert_{\varphi}\vee \Vert\frac{d\lambda_U}{d\pi_{2n_0+1}}\Vert_{\varphi}\right),
$$
as well as
$$\sup_{n\geq n_0}\left(\Vert \frac{d\pi_{2n}}{d\nu_V}\Vert_{\psi}\vee \Vert\frac{d\nu_V}{d\pi_{2n}}\Vert_{\psi}\right)\leq c_{\epsilon}~\left(\Vert \frac{d\pi_{2n_0}}{d\nu_V}\Vert_{\psi}\vee \Vert\frac{d\nu_V}{d\pi_{2n_0}}\Vert_{\psi}\right).
$$
\end{prop}

\proof If Condition $\Da_{\varphi,\psi}(\Sa)$ is satisfied for some integer $n_0\geq 0$ then, by Theorem~\ref{theo-1-intro},  the drift condition (\ref{lyap-eq-sink-o})  also holds for any $n>n_0$, which implies that
$$
\sup_{n\geq n_0}\Vert \frac{d\pi_{2n}}{d\nu_V}\Vert_{\psi}\leq c_{\epsilon}~\Vert \frac{d\pi_{2n_0}}{d\nu_V}\Vert_{\psi}\quad \text{and, similarly,}\quad
\sup_{n\geq  n_0}\Vert \frac{d\pi_{2n+1}}{d\lambda_U}\Vert_{\varphi}\leq c_{\epsilon}~\Vert \frac{d\pi_{2n_0+1}}{d\lambda_U}\Vert_{\varphi}.
$$
Arguing as in the proof of Proposition~\ref{prop-uest-pi-ratio}, now combining the drift conditions (\ref{lyap-eq-sink}) with (\ref{gibbs-ratio-bb}), we find that
$$
\sup_{n\geq n_0}\Vert\frac{d\nu_V}{d\pi_{2n}}\Vert_{\psi}\leq c_{\epsilon}~
\Vert\frac{d\nu_V}{d\pi_{2n_0}}\Vert_{\psi}
\quad \mbox{\rm and}\quad
\sup_{n\geq n_0}\Vert\frac{d\lambda_U}{d\pi_{2n+1}}\Vert_{\varphi}\leq c_{\epsilon}~
\Vert\frac{d\lambda_U}{d\pi_{2n_0+1}}\Vert_{\varphi}.
$$
\cqfd

Under the stronger Condition $\Ha_{\delta}(U,V)$, more refined estimates can be obtained. With that end, consider the functions
\begin{equation}\label{Wa-U-V}
\Wa_{V,U}(x):=\log{\Qa_V(e^{W_{U}})(x)}\quad\mbox{\rm and}\quad \Wa^{U,V}(y):=\log{\Ra_U(e^{W^{V}})(y)},
\end{equation}
with the integrated cost functions $(W_{U},W^{V})$ defined in (\ref{int-costs}) and the integral operators
$$
\Qa_{V}(x,dy)=e^{-W(x,y)}~\nu_V(dy)\quad \mbox{\rm and}\quad
\Ra_{U}(y,dx)=e^{-W(x,y)}~\lambda_U(dx).
$$
The next technical lemma is pivotal.
\begin{lem}\label{lem-ineq}
For any $n\geq 0$ we have
\begin{equation}\label{Q-U-V}
\begin{array}{rcl}
\log{\Qa(\exp{(-V_{2n})})}&\geq& -W^{V}+\nu_V(V) \quad \text{and}\\
&&\\
-\nu_V(V)+\Wa^{U,V}\geq \log{\Ra(\exp{(-U_{2n})})}&\geq& -W_{U},
 \end{array}
\end{equation}
as well as
$$
U_{2n}\geq U-W^V+\nu_V(V).
$$
In addition, for any $n\geq 1$ we have
\begin{equation}\label{Q-U-V-bis}
\Wa_{V,U}\geq \log{\Qa(\exp{(-V_{2n})})}\geq -W^{V}+\nu_V(V)
\end{equation}
as well as
\begin{equation}\label{UU-VV}
\begin{array}{rcccl}
U-W^V+\nu_V(V)&\leq& U_{2n}&\leq& U+\Wa_{V,U} \quad \text{and}\\
&&&&\\
V-W_U&\leq& V_{2n}&\leq& V+\Wa^{U,V}-\nu_V(V).
\end{array}
\end{equation}
\end{lem}

\proof Formula (\ref{form-series}) yields the following monotone properties
\begin{equation}\label{prop-series} 
\begin{array}{rcccccl}
\nu_V(V_{2(n+1)})&\leq&
\nu_V(V_{2n})&=&\nu_V(V_{2(n+1)})+\mbox{\rm Ent}(\nu_V~|~\pi_{2n})&\leq& \nu_V(V_0)=0,\\
&&&&&&\\
\lambda_U(U_{2(n+1)})&\leq &\lambda_U(U_{2n})&=&\lambda_U(U_{2(n+1)})+\mbox{\rm Ent}(\lambda_U~|~\pi_{2n+1})&\leq& \lambda_U(U_0).
\end{array}
\end{equation}
Applying Jensen's inequality, for any $n\geq 0$ we have
\begin{eqnarray*}
\log{\Qa(\exp{(-V_{2n})})(x)}&=&\log{\int~\nu_V(dy)~\exp{(-W(x,y)+V(y)-V_{2n}(y))}}\\
&\geq&
-W^V(x)+\nu_V(V)-\nu_V(V_{2n})\geq -W^{V}(x)+\nu_V(V).
\end{eqnarray*}
By (\ref{prop-series}) we also have $\lambda_U(U)\geq \lambda_U(U_{2n})$ and, therefore, for any $n\geq 0$ we obtain
 \begin{eqnarray*}
\log{\Ra(\exp{(-U_{2n})})(y)}&=&\log{\int~\lambda_U(dx)~\exp{(-W(x,y)+U(x)-U_{2n}(x))}}\\&\geq&
-W_{U}(y)+\lambda_U(U)-\lambda_U(U_{2n})\geq -W_{U}(y)
\end{eqnarray*}
using Jensen's inequality again. This ends the proof of \eqref{Q-U-V}.

The integral equations in (\ref{prop-schp}) imply that for any $n\geq 0$ we have
$$
U_{2n}\geq U-W^{V}+\nu_V(V)
\quad\mbox{\rm and}\quad
V_{2(n+1)} \geq V-W_{U}.
$$
Using again (\ref{prop-schp}) we see that
\begin{eqnarray*}
U_{2(n+1)}-U=\log{\Qa(e^{-V_{2(n+1)}})}&\leq &
\log{\Qa_V(e^{W_{U}})}=\Wa_{V,U} \quad \text{and}\\
V_{2(n+1)}-V=\log{\Ra(e^{-U_{2n}})}&\leq &-\nu_V(V)+
\log{\Ra_U(e^{W^{V}})}=-\nu_V(V)+\Wa^{U,V}.
\end{eqnarray*}
\cqfd

\begin{prop}\label{prop-W-uvb}
If Condition $\Ha_{\delta}(U,V)$ is satisfied for some $\delta\in ]0,1]$, then
\begin{equation}\label{W-uvb}
\Vert e^{\Wa^{U,V}}\Vert_{\LL_{\infty}(\nu)}\vee
\Vert e^{\Wa_{V,U}}\Vert_{\LL_{\infty}(\lambda)}<\infty.
\end{equation}
In addition, the exists some constant $c>0$ such that, for any $n\geq 1$,
\begin{equation}\label{W-uvb-2}
\frac{d\pi_{2n}}{d\nu_V}\vee \frac{d\nu_V}{d\pi_{2n}}\leq c~e^{W_U}\quad \mbox{and}\quad
\frac{d\pi_{2n+1}}{d\lambda_U}\vee\frac{d\lambda_U}{d\pi_{2n+1}}\leq c~e^{W^V}.
\end{equation}
\end{prop}
\proof

Since $U$ is locally bounded with compact sub-level sets for any $\delta\geq 0$ and 
any sufficiently large  $r>0$ we have
$$
e^{-\delta U}=e^{-\delta U}~1_{\{U>r\}}+
e^{-\delta U}~1_{\{U\leq r\}}\leq \varsigma^U_{\delta}(r)\quad \mbox{\rm and}\quad
\Vert e^{-\delta U}\Vert_{\LL_{\infty}(\lambda)}^{-1}\geq 1/\varsigma^U_{\delta}(r)
$$
with the parameters
$$
 \varsigma^U_{\delta}(r):=e^{-\delta r}+e^{\delta c_U(r)}
 \quad \mbox{\rm with}\quad
c_U(r):=\sup_{x\in C_U(r)} \vert U(x)\vert 
\quad \mbox{\rm and}\quad
C_U(r)=\{x~:~U(x)\leq r\}
$$
This yields
\begin{equation}
\lambda_{(1-\delta)U}\left(e^{W^V}\right)=\lambda_{U}\left(e^{W^V}~e^{\delta U}\right)\geq 
\Vert e^{-\delta U}\Vert_{\LL_{\infty}(\lambda)}^{-1}~\lambda_{U}\left(e^{W^V}\right).
\nonumber%\label{eq65.2}
\end{equation}
In the same vein, we obtain
\begin{equation}
\nu_{(1-\delta)V}\left(e^{W_U}\right)=\nu_{V}\left(e^{W_U}~e^{\delta V}\right)\geq 
\Vert e^{-\delta V}\Vert_{\LL_{\infty}(\nu)}^{-1}~\nu_{V}\left(e^{W_U}\right)
\nonumber%\label{eq65.3}
\end{equation}
and conclude that
\begin{eqnarray*}
\Vert e^{\Wa_{V,U}}\Vert_{\LL_{\infty}(\lambda)}&=&
\Vert \Qa_{V}(e^{W_{U}})\Vert_{\LL_{\infty}(\lambda)}\\
&\leq& \Vert e^{-W}\Vert_{\LL_{\infty}(\lambda\otimes\nu)}~\nu_{V}\left(e^{W_U}\right)\\
&\leq &
\Vert e^{-\delta V}\Vert_{\LL_{\infty}(\nu)}~\Vert e^{-W}\Vert_{\LL_{\infty}(\lambda\otimes\nu)}~\nu_{(1-\delta)V}\left(e^{W_U}\right)
\end{eqnarray*}
as well as
$$
\Vert e^{\Wa^{U,V}}\Vert_{\LL_{\infty}(\nu)}=
\Vert\Ra_{U}(e^{W^V})\Vert\leq \Vert e^{-\delta U}\Vert_{\LL_{\infty}(\lambda)}~
\Vert e^{-W}\Vert_{\LL_{\infty}(\lambda\otimes\nu)}~\lambda_{(1-\delta)U}\left(e^{W^V}\right),
$$
where $\lambda_{(1-\delta)U}\left(e^{W^V}\right) \vee \nu_{(1-\delta)V}\left(e^{W_U}\right) < \infty$ by Condition $\Ha_\delta(U,V)$. This ends the proof of (\ref{W-uvb}). 

%Recalling that $U_{2n+1}=U_{2n}$ and $V_{2(n+1)}=V_{2n+1}$ and 
Using (\ref{UU-VV}) we see that, for any $n\geq 1$,
$$
-\Wa_{V,U}-W^{V}+\nu_V(V)\leq
U_{2(n+1)}-U_{2n}\leq \Wa_{V,U}+W^{V}-\nu_V(V)
$$
as well as
$$
-\Wa^{U,V}-W_U+\nu_V(V)
\leq V_{2(n+1)}-V_{2n} \leq 
\Wa^{U,V}+W_U-\nu_V(V).
$$
Then, recalling that $U_{2n+1}=U_{2n}$ and $V_{2(n+1)}=V_{2n+1}$ and using Eqs. \eqref{eqProofConmut1} and \eqref{eqProofConmut2} in the proof of the commutation formulae (\ref{comm}), we readily obtain the estimates
$$
e^{\nu_V(V)}\Vert e^{\Wa^{U,V}}\Vert_{\LL_{\infty}(\lambda)}^{-1}~e^{-W_U}\leq
\frac{d\pi_{2n}}{d\nu_V}\leq e^{-\nu_V(V)}\Vert e^{\Wa^{U,V}}\Vert_{\LL_{\infty}(\lambda)}~e^{W_U}
$$
as well as
$$
e^{\nu_V(V)}\Vert e^{\Wa_{V,U}}\Vert_{\LL_{\infty}(\lambda)}^{-1}~e^{-W^V}\leq
\frac{d\pi_{2n+1}}{d\lambda_U}\leq e^{-\nu_V(V)}\Vert e^{\Wa_{V,U}}\Vert_{\LL_{\infty}(\lambda)}~e^{W^V}
$$
\cqfd

\begin{rmk}
The functions $(\Wa_{V,U},\Wa^{U,V})$ defined in (\ref{Wa-U-V}) can be computed explicitly for the the Gaussian model defined in (\ref{def-U-V}) and (\ref{def-W}). 

We set $\sigma_{\beta}=\beta\sigma\beta^{\prime}$ and $\beta_{\tau}:=\tau^{-1}\beta$. With this notation at hand, after some elementary computations (cf. for instance Exemple 6.6 in~\cite{adm-24}) we have
$$
\Wa_{V,U}(x)=\frac{1}{2}~\tr(\tau^{-1} \sigma_{\beta})+
\frac{1}{2}~(\beta_{\tau}(x-m))^{\prime}
(\overline{\sigma}-\tau)~(\beta_{\tau}~(x-m))+(\overline{m}-m_0)^{\prime}\beta_{\tau}(x-m)
$$
and
\begin{equation}\label{c-mu-2}
\Wa^{U,V}(y)=\frac{1}{2}~\tr(\tau^{-1}\overline{\sigma})+
\frac{1}{2}
(y-\overline{m})^{\prime}\left(\beta_{\tau}\,\sigma\,\beta_{\tau}^{\prime}-\tau^{-1}\right)
(y-\overline{m})-(y-\overline{m})^{\prime}\tau^{-1}(\overline{m}-m_0).
\end{equation}
Note that one readily obtain the relationships
$$
(\ref{W-uvb})
\Longleftrightarrow
\left(\overline{\sigma}<\tau\quad\text{and}\quad
\beta\,\sigma\,\beta^{\prime}<\tau\right)
\Longleftrightarrow
\left(\tau^{-1}<\overline{\sigma}^{-1}\quad \text{and}\quad \beta^{\prime}\tau^{-1}\beta<\sigma^{-1}
\right)
\Longleftarrow (\ref{cc-prop}).
$$
\end{rmk}

Using (\ref{W-uvb-2}), Condition $\Ha_{\delta}(U,V)$ ensures that for any $n\geq 1$ we have
$$
e^{-\delta V}~\left(\frac{d\pi_{2n}}{d\nu_V}\vee\frac{d\nu_V}{d\pi_{2n}}\right)\leq c~e^{-\Va_{\delta}}
\quad\mbox{\rm and}\quad
e^{-\delta U}~~\left(\frac{d\pi_{2n+1}}{d\lambda_U}\vee\frac{d\lambda_U}{d\pi_{2n+1}}\right)\leq
c~e^{-\Ua_{\delta}}
$$
for some constant $c<\infty$. These estimates immediately yield the following result.
\begin{prop}\label{prop-norm-unif}
Under the assumptions of Proposition~\ref{prop-W-uvb},
for any $n\geq 1$ we have
$$
\frac{d\pi_{2n}}{d\nu_V}\vee\frac{d\nu_V}{d\pi_{2n}}\in\Ba_{\psi}(\YY)\quad
\mbox{\rm and}\quad
\frac{d\pi_{2n+1}}{d\lambda_U}\vee\frac{d\lambda_U}{d\pi_{2n+1}}
\in\Ba_{\varphi}(\XX).
$$
In addition, we have the uniform estimates
$$
\sup_{n\geq 1}~\left(\Vert \frac{d\pi_{2n}}{d\nu_V}\Vert_{\psi}\vee \Vert \frac{d\nu_V}{d\pi_{2n}}\Vert_{\psi}\right)<\infty
\quad\mbox{and}\quad
\sup_{n\geq 1}~\left(\Vert \frac{d\pi_{2n+1}}{d\lambda_U}\Vert_{\varphi}\vee \Vert \frac{d\lambda_U}{d\pi_{2n+1}}\Vert_{\varphi}\right)<\infty.
$$
\end{prop}

%%%
%
%%%
\subsection{Lyapunov inequalities}\label{appSchrodinger-lyap}
The result below is a fairly direct consequence of the estimates presented in Lemma \ref{lem-ineq}.

\begin{lem}\label{prop-KQR}
For any $n\geq 1$ we have the estimates
$$
e^{\nu_V(V)}~e^{-(\Wa_{V,U}(x)+\Wa^{U,V}(y))}\leq 
\frac{d\delta_x\Sa_{2n}}{d\delta_x\Qa_V}(y)\leq e^{-\nu_V(V)}~e^{W^V(x)+W_{U}(y)}
$$
and
$$
e^{\nu_V(V)}~e^{-(\Wa_{V,U}(x)+\Wa^{U,V}(y))}
\leq \frac{d\delta_y\Sa_{2n+1}}{d\delta_y\Ra_{U}}(x)\leq e^{-\nu_V(V)}~e^{W^V(x)+W_U(y)}.
$$
\end{lem}

\proof Using the Sinkhorn transition operators as expressed in \eqref{sinhorn-transitions-form-Sch} together with (\ref{Q-U-V-bis}) and (\ref{UU-VV}), we find that, for any $n\geq 1$,
\begin{eqnarray*}
\Sa_{2n}(x,dy)&=&\frac{\Qa_V(x,dy)e^{(V-V_{2n})(y)}}{\Qa(e^{-V_{2n}})(x)}
\geq e^{\nu_V(V)} e^{-(\Wa_{V,U}(x)+\Wa^{V,U}(y))}
\Qa_V(x,dy)~
\end{eqnarray*}
and, in the same vein,
\begin{eqnarray*}
\Sa_{2n+1}(y,dx)&=&\frac{\Ra_U(y,dx)~e^{(U-U_{2n})(x)}}{\Ra(e^{-U_{2n}})(y)}\geq e^{\nu_V(V)}~e^{-(\Wa_{V,U}(x)+\Wa^{U,V}(y))}~\Ra_U(y,dx).
\end{eqnarray*}

Similarly, one can resort to (\ref{Q-U-V-bis}) and (\ref{UU-VV}) again to obtain the lower bounds
\begin{eqnarray*}
\Sa_{2n}(x,dy)&=&\frac{\Qa_V(x,dy)e^{V-V_{2n}(y)}}{\Qa(e^{-V_{2n}})(x)}\leq  e^{-\nu_V(V)}~e^{W^{V}(x)}~\Qa_V(x,dy)~e^{W_{U}(y)}
\end{eqnarray*}
and, in the same way,
\begin{eqnarray*}
\Sa_{2n+1}(y,dx)&=&\frac{\Ra_U(y,dx)~e^{(U-U_{2n})(x)}}{\Ra(e^{-U_{2n}})(y)}\leq 
e^{-\nu_V(V)}~e^{W_U(y)}~\Ra_U(y,dx)~e^{W^V(x)}.
\end{eqnarray*}
\cqfd

With these estimates at hand, we can now state and prove the following Lyapunov inequalities.

\begin{prop} \label{propLyapunovInequalities}
Assume that Condition $\Ha_{\delta}(U,V)$ is satisfied for some $\delta\in ]0,1]$. 
Then, the functions $(\varphi,\psi)=(e^{\delta U},e^{\delta V})$ satisfy the Lyapunov inequalities 
\begin{equation}\label{Lyap-S}
\Sa_{2n}\left(\psi\right)/\varphi \leq 
c_{\delta}~e^{-\Ua_{\delta}}~\in\Ba_0(\XX)
\quad \mbox{and}\quad
\Sa_{2n-1}\left(\varphi\right)/\psi \leq 
c_{\delta}~e^{-\Va_{\delta}}~
\in\Ba_0(\YY)
\end{equation}
for every $n\geq 1$ and some finite constant $c_{\delta}<\infty$.
\end{prop}

\proof
By Lemma~\ref{prop-KQR} we have
\begin{eqnarray*}
\Sa_{2n}(x,dy)~\psi(y)/\varphi(x) &\leq &
e^{-\nu_V(V)}~e^{-\Ua_{\delta}(x)}~\Qa_{(1-\delta)V}(x,dy)~e^{W_{U}(y)}, 
\quad \text{and}
\\
\Sa_{2n+1}(y,dx)~\varphi(x)/\psi(y) &\leq &
e^{-\nu_V(V)}~e^{-\Va_{\delta}(y)}~\Ra_{(1-\delta)U}(y,dx)~e^{W^V(x)},
\end{eqnarray*}
which (via \eqref{Wa-U-V}) imply the inequalities
\begin{eqnarray*}
\Sa_{2n}\left(\psi\right)/\varphi&\leq  &
e^{-\nu_V(V)}~e^{-\Ua_{\delta}}~e^{\Wa_{(1-\delta)V,U}}, \quad \text{and}
\\
\Sa_{2n+1}\left(\varphi\right)/\psi &\leq &
e^{-\nu_V(V)}~e^{-\Va_{\delta}}~e^{\Wa^{(1-\delta)U,V}}.
\end{eqnarray*}
On the other hand, we have the estimates
$$
\Vert e^{\Wa_{(1-\delta)V,U}}\Vert_{\LL_{\infty}(\lambda)}=
\Vert\Qa_{(1-\delta)V}(e^{W_{U}})\Vert_{\LL_{\infty}(\lambda)}\leq \Vert e^{-W}\Vert_{\LL_{\infty}(\lambda\otimes\nu)}~\nu_{(1-\delta)V}\left(e^{W_U}\right)
$$
and
$$
\Vert e^{\Wa^{(1-\delta)U,V}}\Vert_{\LL_{\infty}(\nu)}=
\Vert\Ra_{(1-\delta)U}(e^{W^V})\Vert_{\LL_{\infty}(\nu)}\leq 
\Vert e^{-W}\Vert_{\LL_{\infty}(\lambda\otimes\nu)}~\lambda_{(1-\delta)U}\left(e^{W^V}\right),
$$
where $\nu_{(1-\delta)V}\left(e^{W_U}\right) \vee \lambda_{(1-\delta)U}\left(e^{W^V}\right) < \infty$ when Condition $\Ha_\delta(U,V)$ holds. This yields (\ref{Lyap-S})
with the constant
$$
c_{\delta}:=e^{-\nu_V(V)}~\Vert e^{-W}\Vert_{\LL_{\infty}(\lambda\otimes\nu)}~\left(\nu_{(1-\delta)V}\left(e^{W_U}\right)\vee
\lambda_{(1-\delta)U}\left(e^{W^V}\right)\right).
$$
By (\ref{Q-U-V}) we have
\begin{eqnarray*}
\Sa_{1}(y,dx)&=&\frac{\Ra(y,dx)~e^{-U(x)}}{\Ra(e^{-U})(y)}\leq e^{W_U(y)} ~\Ra(y,dx)~e^{-U(x)}
\end{eqnarray*}
This implies that
$$
e^{-\delta V(y)}
\Sa_{1}\left(e^{\delta U}\right)(y) \leq 
e^{-(\delta V(y)-W_U(y))} ~\Vert e^{-W}\Vert~\lambda_{(1-\delta)U}\left(1\right)
$$
Equivalently, some finite constant $c_{\delta}<\infty$ we have
$$
\Sa_{1}\left(\varphi\right)/\psi \leq 
c_{\delta}~e^{-\Va_{\delta}}~
$$
This ends the proof of the proposition.
\cqfd

\subsection{Proof of Theorem~\ref{lm-tx-2}}\label{lm-tx-2-proof}

By (\ref{ref-dom-theta-0}), the Lyapunov inequalities (\ref{Lyap-S}) in Proposition \ref{propLyapunovInequalities} ensure that the drift conditions stated in (\ref{lyap-eq-sink}) hold with $n_0=1$ and $(\varphi,\psi)=(e^{\delta (U-U_{\star})},e^{\delta (V-V_{\star})})$.
Recall from the proof of Proposition~\ref{prop-W-uvb} that $e^{-\delta U}$ and $e^{-\delta V}$ are uniformly bounded for any $\delta\geq 0$.
Since $\XX$ and $\YY$ can be exhausted respectively by the compact sub-levels sets
$C_{U}(r)$ and $C_{V}(r)$ of $U$ and $V$, there exists some $r_1>1$ such that $\lambda(C_{U}(r))\wedge \nu(C_V(r))>0$  for $r=r_1$ and thus for any $r\geq r_1$. This ensures there exists some $r_0>0$ such that for any $r\geq r_0$ we have
$$
\lambda(C_{\varphi}(r))\wedge \nu(C_{\psi}(r))>0.
$$

Following the proof of Proposition~\ref{prop-W-uvb}, for $\lambda$-almost any $x\in C_{\varphi}(r)$  we have
\begin{eqnarray*}
e^{\Wa_{V,U}(x)}=\Qa_V(e^{W_{U}})(x)&=&\int~e^{W_U(y)-W(x,y)}~\nu_V(dy)\\
&\leq& \Vert e^{-W}\Vert_{\LL_{\infty}(\lambda\otimes\nu)}~\nu_V(e^{W_U})\\
&\leq&
\Vert e^{-W}\Vert_{\LL_{\infty}(\lambda\otimes\nu)}~\Vert e^{-\delta V}\Vert_{\LL_{\infty}(\nu)}~ \nu_{(1-\delta)V}\left(e^{W_U}\right)
\end{eqnarray*}
and we conclude that
$$
w_{V,U}(r):=\sup_{C_{\varphi}(r)}\Wa_{V,U}<\infty.
$$
Similarly, for $\nu$-almost any $y\in C_{\psi}(r)$  we have
\begin{eqnarray*}
e^{\Wa^{U,V}(y)}=\Ra_U(e^{W^{V}})(y)&=&\int~e^{-W(x,y)+W^{V}(x)}~\lambda_U(dx)\\
&\leq &\Vert e^{-W}\Vert_{\LL_{\infty}(\lambda\otimes\nu)}~\lambda_U(e^{W^{V}})\\
&\leq& 
\Vert e^{-W}\Vert_{\LL_{\infty}(\lambda\otimes\nu)}~\Vert e^{-\delta U}\Vert_{\LL_{\infty}(\lambda)}~\lambda_{(1-\delta)U}\left(e^{W^{V}}\right)
\end{eqnarray*}
which yields
$$
w^{U,V}(r):=\sup_{C_{\psi}(r)}\Wa^{U,V}<\infty.
$$
On the other hand, using Lemma~\ref{prop-KQR} for any  $x\in C_{\varphi}(r)$ and $r\geq r_0$ we arrive at
\begin{eqnarray*}
\Sa_{2n}(x,dy)&\geq& 
e^{\nu_V(V)}~e^{-(W(x,y)+\Wa_{V,U}(x)+\Wa^{U,V}(y))}~\nu_V(dy)~1_{C_{\psi}(r)}(y)\\
&\geq & e^{-w(r)}~\nu_V\left(C_{\psi}(r)\right)~\Ya_r(dy)\quad\mbox{\rm with}\quad
\Ya_r(dy):=\frac{\nu_V(dy)~1_{C_{\psi}(r)}(y)}{\nu_V\left(C_{\psi}(r)\right)}
\end{eqnarray*}
and the parameter
$$
w(r):=\sup_{C_{\varphi}(r)\times C_{\psi}(r)}W+
w_{V,U}(r)+w^{U,V}(r)-\nu_V(V).
$$
Note that for any $r\geq r_0$ we have
\begin{eqnarray}
\nu_V\left(C_{\psi}(r)\right)&=&e^{-V_{\star}}~\nu\left(e^{-(V-V_{\star})}~1_{(V-V_{\star}) \leq \delta^{-1}\log{r}}\right) \nonumber\\
&\geq& w_0(r):=e^{-V_{\star}}~e^{-(\delta^{-1}\log{r})}~\nu\left(C_{\psi}(r)\right)>0 \nonumber
\end{eqnarray}
and this yields, for any  $x\in C_{\varphi}(r)$ and $r\geq r_0$, the estimate
$$
\Sa_{2n}(x,dy)\geq \varepsilon_0(r)
~\Ya_r(dy)
\quad\mbox{\rm with}\quad
\varepsilon_0(r):=w_0(r)~e^{-w(r)}.
$$
In the same vein, for any $y\in C_{\psi}(r)$ we have
\begin{eqnarray*}
\Sa_{2n+1}(y,dx)&\geq &e^{\nu_V(V)}~e^{-(W(x,y)+\Wa_{V,U}(x)+\Wa^{U,V}(y))}~\lambda_U(dx)~1_{C_{\varphi}(r)}(x)
\end{eqnarray*}
and, arguing as above, we arrive at
$$
\lambda_U\left(C_{\varphi}(r)\right)\geq
w_1(r):=e^{-U_{\star}}~e^{-(\delta^{-1}\log{r})}~\lambda\left(C_{\varphi}(r)\right)>0.
$$
Thus, for any $y\in C_{\psi}(r)$ we have
$$
\Sa_{2n+1}(y,dx)\geq \varepsilon_1(r)~\Xa_r(dx)
$$
where
$$
\varepsilon_1(r):=w_1(r)~e^{-w(r)}
\quad\mbox{\rm and}\quad
\Xa_r(dx):=\frac{\lambda_U(dx)~1_{C_{\varphi}(r)}(x)}{\lambda_U\left(C_{\varphi}(r)\right)}.
$$
The above estimates ensure that the local minorization condition (\ref{sink-locmin}) is met with $n_0=1$ and the function
$$
\varepsilon(r):= \varepsilon_0(r)\wedge \varepsilon_1(r)>0.
$$
Hence, we conclude that Condition $\Da_{\varphi,\psi}(\Sa)$ is satisfied  with $n_0=1$.

Using (\ref{fixed-points-gibbs}) and (\ref{gibbs-ratio}), for any $n\geq p\geq 1$ we verify that 
$$
\frac{d\pi_{2n}}{d\nu_V}(y)-1=\left[(\delta_y-\nu_V)\left(\Sa^{\circ}_{2n}\ldots \Sa^{\circ}_{2(p+1)}\right)\right]\left(\frac{d\pi_{2p}}{d\nu_V}\right),
$$
which implies the inequality
$$
\left\vert \frac{d\pi_{2n}}{d\nu_V}(y)-1\right\vert\leq \Vert \frac{d\pi_{2p}}{d\nu_V}\Vert_{\psi} ~\vertiii{(\delta_y-\nu_V)\left(\Sa^{\circ}_{2n}\ldots \Sa^{\circ}_{2(p+1)}\right)}_{\psi}.
$$
On the other hand, combining Proposition~\ref{prop-norm-unif} with  the norm equivalence (\ref{norm-equivalence}) and Theorem~\ref{theo-1-intro}, we see that there exists some $\kappa>0 $ and some constants $c_{\kappa,i}>0$, $i=1,2$, such that
\begin{eqnarray*}
~\vertiii{(\delta_y-\nu_V)\left(\Sa^{\circ}_{2n}\ldots \Sa^{\circ}_{2(p+1)}\right)}_{\psi}&\leq &c_{\kappa,1}
~\vertiii{(\delta_y-\nu_V)\left(\Sa^{\circ}_{2n}\ldots \Sa^{\circ}_{2(p+1)}\right)}_{\psi_{\kappa}}\\
&\leq & c_{\kappa,1}~\rho^{2(n-p)}~\vert \delta_y-\nu_V\vert_{\psi_{\kappa}}\leq  c_{\kappa,2}~ \rho^{2(n-p)}~(\psi(y)+\nu_V(\psi)).
\end{eqnarray*}
Therefore, for any $n\geq 1$, we obtain
\begin{equation}\label{eqPart-1}
\left\|
    \frac{
        d\pi_{2n}
    }{
        d\nu_V
    }-1
\right\|_{\psi} \leq
c_{\kappa,3}~ \rho^{2n}~
\quad\text{and, similarly,}\quad
\left\|\frac{d\nu_V}{d\pi_{2n}}-1\right\|_{\psi}\leq
c_{\kappa,4}~\rho^{2n},
\end{equation}
for some constants $c_{\kappa,i}$, $i=3,4$. By the same argument, using (\ref{fixed-points-gibbs}) we have, for any $n\geq p\geq 1$,
$$
  \frac{d \pi_{2n+1}}{d\lambda_U}(x)-1=  \left[(\delta_x-\lambda_U)(\Sa^{\circ}_{2n+1}\ldots \Sa^{\circ}_{2p+1})\right]\left(\frac{d \pi_{2p-1}}{d\lambda_U}\right)
$$
and arrive at
\begin{equation}\label{eqPart-2}
\left\|\frac{d \pi_{2n+1}}{d\lambda_U}-1\right\|_{\varphi}\vee \left\|\frac{d\lambda_U}{d \pi_{2n+1}}-1\right\|_{\varphi}\leq c_{\kappa,5}~\rho^{2n}~
\end{equation}
for any $n\geq 1$ and some constant $c_{\kappa,5}$.

Finally, note that for any function $f(x)>0$ we have
\begin{eqnarray*}
\left\vert\log{f(x)}\right\vert&=& \log{\left(1+(f(x)-1)\right)}~1_{ f(x)\geq 1}+
\log{\left(1+\frac{1}{ f(x)}-1\right)}~1_{ {1}/{ f(x)}\geq 1}\nonumber\\
&\leq&\left(f(x)-1\right)~1_{ f(x)\geq 1}+\left(\frac{1}{ f(x)}-1\right)~1_{ {1}/{ f(x)}\geq 1},
\end{eqnarray*}
where the last assertion comes from the fact that $ \log{(1+u)}\leq u$,  for any $1+u>0$.
This implies that
$$
\Vert \log{f}\Vert_{\varphi}=\Vert \log{1/f}\Vert_{\varphi}\leq \Vert f-1\Vert_{\varphi}+
\Vert 1/f-1\Vert_{\varphi}.
$$
The inequality above combined with \eqref{eqPart-1} and \eqref{eqPart-2} yields the estimates
$$
\left\Vert\log{\frac{d \pi_{2n}}{d\nu_V}}\right\Vert_{\psi}\vee
\left\Vert\log{\frac{d \pi_{2n+1}}{d\lambda_U}}\right\Vert_{\varphi}\leq c_{6,\kappa}~ \rho^{2n},
$$
for some finite constant $c_{\kappa,6}$. This completes the proof of (\ref{log-est}).
\cqfd

\section{Some technical proofs}\label{tech-proof-ap}

\subsection*{Proof of  (\ref{comm}) and (\ref{comm-2})}\label{form-series-proof} 
We first note, using \eqref{s-2} and \eqref{sinhorn-entropy-form-Sch}, that 
\begin{equation}
\frac{ d\Pa_{2n}}{ d\Pa_{2n+1}}(x,y)=\frac{d\pi_{2n}}{d\nu_V}(y)=\exp{\left(V_{2n+1}(y)-V_{2n}(y)\right)}
\label{eqProofConmut1}
\end{equation}
and, in the same vein, 
\begin{equation}
\frac{ d\Pa_{2n+1}}{ d\Pa_{2(n+1)}}(x,y)=\frac{d\pi_{2n+1}}{d\lambda_U}(x)=\exp{\left(U_{2(n+1)}(x)-U_{2n+1}(x)\right)}.
\label{eqProofConmut2}
\end{equation}
If we additionally observe that
\begin{eqnarray*}
 \Pa_{2(n+1)}(d(x,y))&=&  \frac{d\lambda_U}{d \pi_{2n+1}}(x)~\Pa_{2n}(d(x,y))~\frac{d\nu_V}{d\pi_{2n}}(y) \quad \text{and}\\
  \Pa_{2n+1}(d(x,y))&=&  \frac{d\lambda_U}{d \pi_{2n-1}}(x)~\Pa_{2n-1}(d(x,y))~\frac{d\nu_V}{d\pi_{2n}}(y)
\end{eqnarray*}
then, using (\ref{def-Pa-n}), we obtain the formulae
\begin{eqnarray*}
 \Sa_{2(n+1)}(x,dy)&=&  \frac{d\lambda_U}{d \pi_{2n+1}}(x)~\Sa_{2n}(x,dy)~\frac{d\nu_V}{d\pi_{2n}}(y) \quad \text{and}\\
  \Sa_{2n+1}(y,dx)&=& ~\frac{d\nu_V}{d\pi_{2n}}(y)~ \Sa_{2n-1}(y,dx)~\frac{d\lambda_U}{d \pi_{2n-1}}(x),
\end{eqnarray*}
and the commutation properties now follow from elementary manipulations.\cqfd

\subsection*{Proof of (\ref{uest-pi})}\label{uest-pi-proof}

Assume Condition $\Da_{\varphi,\psi}(\Sa^{\circ})$ is satisfied. Combining (\ref{lyap-eq-sink-o}) with (\ref{fixed-points-gibbs}) we arrive at
$$
\nu_V(\psi)= \nu_V(\Sa_{2n}^{\circ}(\psi))\leq \epsilon~\nu_V(\psi)+c
\quad \text{which implies} \quad
\nu_V(\psi)\leq\frac{c}{1-\epsilon}.
$$
In a similar way, one can also show that
$$
\lambda_U(\varphi)\leq\frac{c}{1-\epsilon}.
$$
Using the Markov evolution equations in \eqref{gibbs-tv} together with (\ref{lyap-eq-sink-o}) we find that, for any $n\geq n_0$,
\begin{eqnarray*}
\pi_{2n}(\psi)&=&\pi_{2(n-1)} (\Sa_{2n}^{\circ}(\psi))\\
&\leq& \epsilon~\pi_{2(n-1)}(\psi)+c\leq  \epsilon^{n-n_0}~\pi_{2n_0}(\psi) +\frac{c}{1-\epsilon}
\end{eqnarray*}
and, by the same argument,
\begin{eqnarray*}
\pi_{2n+1}(\varphi)&=&\pi_{2n-1} (\Sa_{2n+1}^{\circ}(\varphi))\\
&\leq& \epsilon~\pi_{2n-1}(\varphi)+c\leq  \epsilon^{n-n_0}~\pi_{2n_0+1}(\varphi) +\frac{c}{1-\epsilon}.
\end{eqnarray*}
\cqfd

\subsection*{Proof of Proposition~\ref{prop-gauss-model}}\label{prop-gauss-model-proof}

It is straightforward to show that 
\begin{eqnarray*}
2W_U(y)&=&\left\Vert \tau^{-1/2}(y-m_0)\right\Vert_F^2+\tr(\tau^{-1}\sigma_{\beta})+\log{\mbox{\rm det}(2\pi \tau)} \quad \text{and}\\
2W^V(x)&=&\left\Vert \tau^{-1/2}(\beta (x-m)+(m_0-\overline{m}))\right\Vert_F^2+~\tr(\tau^{-1}\overline{\sigma})+\log{\mbox{\rm det}(2\pi \tau)},
\end{eqnarray*}
with $m_0:=(\alpha+\beta m)$ and $\sigma_{\beta}:=(\beta\sigma
 \beta^{\prime})$, hence $W_U$ and $W^V$ are locally bounded. Moreover, for any $\delta\geq 0$ we obtain
\begin{eqnarray*}
2(\delta U-W^V) &=& -\tr(\tau^{-1}\overline{\sigma})+\delta~\log{\mbox{\rm det}(2\pi\sigma)}-\log{\mbox{\rm det}(2\pi\tau)} \\
&& +(x-m)^{\prime}\left(\delta~\sigma^{-1}-\beta^{\prime}\tau^{-1}\beta\right)(x-m)\\
&& -2(m_0-\overline{m})^{\prime}\tau^{-1}\beta (x-m)-(m_0-\overline{m})^{\prime}\tau^{-1}(m_0-\overline{m})
\end{eqnarray*}
and, in the same way,
\begin{eqnarray*}
2(\delta V-W_U) &=& -\tr(\tau^{-1}\sigma_{\beta})+\delta~\log{\mbox{\rm det}(2\pi\overline{\sigma})}-\log{\mbox{\rm det}(2\pi\tau)}\\
&& +(y-\overline{m})^{\prime}\left(\delta~\overline{\sigma}^{-1}-\tau^{-1}\right)(y-\overline{m})\\
&& -2(\overline{m}-m_0)^{\prime}\tau^{-1}(y-\overline{m}) -(\overline{m}-m_0)^{\prime}\tau^{-1}(\overline{m}-m_0).
\end{eqnarray*}
Therefore, 
$$
\Va_\delta = \delta V-W_U\in \Ba_{\infty}(\RR^d)\quad \mbox{\rm and}\quad \Ua_\delta=\delta U-W^V\in \Ba_{\infty}(\RR^d)
$$
if, and only if,
$$
\beta^{\prime}\tau^{-1}\beta<\delta~\sigma^{-1}\quad \mbox{\rm and}\quad
\tau^{-1}<\delta~\overline{\sigma}^{-1}.
$$

On the other hand, straightforward manipulations yield
\begin{eqnarray*}
2(W^V-(1-\delta) U) &=& \tr(\tau^{-1}\overline{\sigma})-(1-\delta)~\log{\mbox{\rm det}(2\pi\sigma)}+\log{\mbox{\rm det}(2\pi\tau)}\\
&&-(x-m)^{\prime}\left((1-\delta)~\sigma^{-1}-\beta^{\prime}\tau^{-1}\beta\right)(x-m)\\
&&+2(m_0-\overline{m})^{\prime}\tau^{-1}\beta (x-m) + (m_0-\overline{m})^{\prime}\tau^{-1}(m_0-\overline{m}),
\end{eqnarray*}
which shows that
$$
\lambda_{(1-\delta)U}(e^{W^V})<\infty\Longleftrightarrow
\beta^{\prime}\tau^{-1}\beta<(1-\delta)~\sigma^{-1}.
$$
In the same vein, we have
\begin{eqnarray*}
2(W_U-(1-\delta) V) &=& \tr(\tau^{-1}\sigma_{\beta})-(1-\delta)~\log{\mbox{\rm det}(2\pi\overline{\sigma})}-\log{\mbox{\rm det}(2\pi\tau)}\\
&&-(y-\overline{m})^{\prime}\left((1-\delta)~\overline{\sigma}^{-1}-\tau^{-1}\right)(y-\overline{m})\\
&&+2(\overline{m}-m_0)^{\prime}\tau^{-1}(y-\overline{m}) + (\overline{m}-m_0)^{\prime}\tau^{-1}(\overline{m}-m_0),
\end{eqnarray*}
which shows that
$$
\nu_{(1-\delta)V}(e^{W_U})<\infty\Longleftrightarrow
\tau^{-1}<(1-\delta)~\overline{\sigma}^{-1}.
$$
\cqfd

\subsection*{Proof of Proposition~\ref{prop-mixture-intro}}\label{prop-mixture-intro-proof}

Let us first recall that
$$
\lambda_U(dx)=\sum_{i\in E_0} \xi_0(i)~\lambda_{U_i}(dx)
\quad \mbox{\rm and}\quad
\nu_V(dy)=\sum_{j\in E_1} \xi_1(j)~\nu_{V_j}(dx),
$$
which implies that
\begin{equation}
W_U=\sum_{i\in E_0}~\xi_0(i)~W_{U_i}\quad \mbox{\rm and}\quad
W^V=\sum_{j\in E_1}~\xi_1(j)~W^{V_j}
\label{eqC5.5}
\end{equation}
and ensures that the functions $W_U$ and $W^V$ are locally bounded.

Moreover, we note that
$$
\sum_{i\in E_0} \xi_0(i)~e^{-U_i}\leq e^{-\inf_{i\in E_0}U_i},
$$
which implies that
$$
-U=\log{\sum_{i\in E_0} \xi_0(i)~e^{-U_i}}\leq -\inf_{i\in E_0}U_i
$$ 
and, hence,
$$
U\geq \inf_{i\in E_0}U_i.
$$
As a consequence, we readily obtain the inequalities
$$
\delta U-W^V\geq \inf_{i\in E_0}\sum_{j\in E_1}~\xi_1(j)~(\delta~U_i-~W^{V_j})\geq 
\inf_{(i,j)\in (E_0\times E_1)}~(\delta~U_i-~W^{V_j})
$$
which, in turn, enable us to show that
\begin{eqnarray*}
\left\{x~:~\delta U(x)-W^V(x)\leq r\right\} &\subset&
\left\{x~:~\inf_{(i,j)\in (E_0\times E_1)}~(\delta~U_i(x)-~W^{V_j}(x))\leq r\right\}\\
&=&\bigcup_{(i,j)\in (E_0\times E_1)}
\left\{
x~:~(\delta~U_i(x)-~W^{V_j}(x))\leq r
\right\}
\end{eqnarray*}
and this ensures that $\Ua_\delta=\delta U-W^{V}$ have compact sub-level sets. In the same vein, one can also prove that
$\Va_\delta=\delta V-W_{U}$ have compact sub-level sets.

On the other hand, using the fact that $(a+b)^{1-\delta}\leq a^{1-\delta}+b^{1-\delta}$ for any $a,b\geq 0$, we have
\begin{eqnarray}
\lambda_{(1-\delta)U}(dx)&=&\left(\sum_{i\in E_0} \xi_0(i)~e^{-U_i(x)}\right)^{1-\delta}~\lambda(dx) \nonumber\\
&\leq& \sum_{i\in E_0} \xi_0(i)^{1-\delta}~e^{-(1-\delta)U_i(x)}~\lambda(dx)=\sum_{i\in E_0} \xi_0(i)^{1-\delta}~\lambda_{(1-\delta)U_i}(dx) \label{eqA5.5}
\end{eqnarray}
and, by way of H\"older's inequality, we also have
\begin{equation}
\lambda_{(1-\delta)U_i}\left(\prod_{j\in E_1}~e^{\xi_1(j)~W^{V_j}}\right)\leq 
\prod_{j\in E_1}~\left(\lambda_{(1-\delta)U_i}\left(e^{W^{V_j}}\right)\right)^{\xi_1(j)}.
\label{eqB5.5}
\end{equation}
Combining \eqref{eqC5.5}, \eqref{eqA5.5} and \eqref{eqB5.5} we obtain the estimates
\begin{eqnarray*}
\lambda_{(1-\delta)U}\left(e^{W^V}\right)&\leq& 
\sum_{i\in E_0} \xi_0(i)^{1-\delta}~\lambda_{(1-\delta)U_i}\left(\prod_{j\in E_1}~e^{\xi_1(j)~W^{V_j}}\right)\\
&\leq &\sum_{i\in E_0} \xi_0(i)^{1-\delta}~\prod_{j\in E_1}~\left(\lambda_{(1-\delta)U_i}\left(e^{W^{V_j}}\right)\right)^{\xi_1(j)}
\end{eqnarray*}
from which we conclude that $\lambda_{(1-\delta)U}\left(e^{W^V}\right)<\infty$ (because $\lambda_{(1-\delta)U_i}\left(e^{W^{V_j}}\right)<\infty$ for almost every $(i,j)$). Using the same argument, we can also prove that $\nu_{(1-\delta)V}\left(e^{W_U}\right)
<\infty$. \cqfd

\bibliography{references}
\bibliographystyle{plain}

\end{document}